\documentclass[reqno,11pt]{amsart}
\usepackage{amssymb, amsmath,latexsym,amsfonts,amsbsy, amsthm}
\usepackage{color}
\usepackage{xcolor}
\usepackage{pgfplots}
\pgfplotsset{compat=1.18}

\allowdisplaybreaks[4]

\newcommand{\beq}{\begin{equation}}
	\newcommand{\eeq}{\end{equation}}
\newcommand{\ben}{\begin{eqnarray}}
	\newcommand{\een}{\end{eqnarray}}
\newcommand{\beno}{\begin{eqnarray*}}
	\newcommand{\eeno}{\end{eqnarray*}}

\newtheorem{remark}{Remark}[section]

\renewcommand{\theequation}{\thesection.\arabic{equation}}

\newtheorem{theorem}{Theorem}[section]

\newtheorem{lemma}[theorem]{Lemma}
\newtheorem{proposition}[theorem]{Proposition}

\numberwithin{equation}{section} 
\newtheorem{Theorem}{Theorem}[section]

\newtheorem{Proposition}[Theorem]{Proposition}
\newtheorem{Lemma}[Theorem]{Lemma}
\newtheorem{Corollary}[Theorem]{Corollary}

\begin{document}
	\title[Vanishing viscosity limit in a channel]
	{Prandtl-Batchelor flow
in a Periodic Channel with Counter-Moving Walls}
	\date{\today}

	\author{Zeyu Jia}
\address{School of  Mathematics, Shandong University, jinan, Shandong, 250100, China}
\email{jzy280012@126.com}

\author{Tao Tao}
\address{School of  Mathematics, Shandong University, jinan, Shandong, 250100, China}
\email{taotao@sdu.edu.cn}

	\date{\today}
\maketitle

\renewcommand{\theequation}{\thesection.\arabic{equation}}
\setcounter{equation}{0}
\begin{abstract}
We study the vanishing-viscosity limit of steady two-dimensional incompressible Navier–Stokes flows in a horizontally periodic channel with counter-moving walls. The prescribed wall velocities are small periodic perturbations of constant velocities with opposite signs. The Batchelor–Wood relations select an affine Euler shear with constant vorticity and a unique interior zero. We establish the existence of steady Navier–Stokes solutions that are uniformly approximated by the selected Euler shear together with two Prandtl boundary layers of thickness $O(\varepsilon)$, with an approximation error of order $O(\epsilon)$, where $\epsilon^2$ denotes the viscosity. Moreover, the vorticity converges locally uniformly to the selected constant vorticity at a rate of $O(\epsilon)$.

The principal analytical difficulty is the loss of global coercivity caused by the interior zero of the Euler shear. To overcome this difficulty, we develop a compactness–rigidity argument that combines signed near-wall positivity estimates, interior vorticity estimates, adjoint tests eliminating the affine limiting modes, and a streamline circulation equation ruling out concentration near the walls. These ingredients yield an $O(\epsilon^{-2})$ estimate for the linearized inverse, which, together with a high-order Euler–Prandtl approximation, enables us to solve the nonlinear remainder equations. Our result provides a rigorous realization of the Prandtl–Batchelor selection principle in a periodic channel with counter-moving walls.
\end{abstract}

\makeatletter
\begingroup
\renewcommand{\thefootnote}{}
\renewcommand{\@makefntext}[1]{\parindent=1.5em\indent#1}
\footnotetext[0]{%
  }
\endgroup
\makeatother

\numberwithin{equation}{section}

\indent

	\section{Introduction}
	
	The vanishing-viscosity limit of the steady incompressible Navier–Stokes equations is a fundamental and challenging problem in mathematical fluid mechanics. In the presence of solid boundaries, the discrepancy between the no-slip condition for viscous flows and the impermeability condition for inviscid flows generally leads to the formation of boundary layers. A further difficulty in the stationary setting arises from the nonuniqueness of steady Euler flows: the inviscid equations may admit infinitely many solutions, and a selection principle is needed to identify those that can arise as limits of viscous flows. In bounded domains, the classical Prandtl--Batchelor principle
	provides a framework for this selection problem under suitable
	assumptions. For two-dimensional steady flows with closed streamlines, the classical Prandtl–Batchelor principle provides such a selection mechanism under suitable assumptions, predicting that the limiting vorticity is constant in regions filled with closed streamlines; see \cite{B,prandtl} and
	\cite{edwards,feymann-lagerstrom,kim-formula,riley81,van wijngarden,W}.
	Related investigations can also be found in \cite{OK,OK1}
	
	The mathematical theory of Prandtl--Batchelor flows has been
	investigated in a series of works by Kim and collaborators,
	including studies on a disk 
	\cite{K1998,K2000,K2003,kim-childress,kim-free boundary,kim-thesis,kim-formula}
	and on an annulus \cite{kim-annular}.
	These studies address the well-posedness of the Prandtl equations
	under the Batchelor--Wood condition and the asymptotic structure
	of boundary layer expansions.
	A further step is to justify these expansions for exact
	Navier--Stokes solutions as the Reynolds number tends to infinity.
	For a disk with a prescribed wall velocity that is a small
	perturbation of rigid rotation, this was achieved in \cite{FGLT1}
	through the construction of higher-order approximate solutions
	and a rigorous analysis of the remainder.
	Related results were established for an annulus with an external
	force in \cite{FGLT2}, for an expanding annulus with degenerate inner velocity in \cite{Pan2026} and for a disk containing a point vortex
	in \cite{CFLZ}. We also refer to \cite{Guo2026} for recent progress on the validity of steady Prandtl expansions in three-dimensional axisymmetric flows.
	
	These developments are closely related to the mathematical justification of steady Prandtl boundary layer expansions. Significant progress has been made for flows over moving boundaries, favorable shear profiles, and more general background Euler flows; see [13, 14, 17--19, 24--30] and the references therein. A central ingredient in many of these constructions is a positivity estimate for the linearized remainder equations. Such estimates typically rely on the sign
and the nondegeneracy of the tangential background velocity, thereby providing the coercivity needed to
control the singular perturbation induced by the small viscosity; see, in particular, \cite{FGLT2, FPZ}.
	
	In contrast to the configurations considered in many previous works, we study a horizontally periodic channel whose two walls move in opposite directions. The selected Euler shear consequently changes sign and vanishes at a unique interior level. This interior degeneracy prevents the direct application of conventional global positivity estimates, even though the background tangential velocity remains nondegenerate near both walls.

The principal difficulty is therefore to connect the coercive estimates available in the two near-wall regions with the behavior of the remainder across the interior zero. Moreover, local compactness alone cannot exclude concentration of the horizontal velocity near the walls. To overcome these difficulties, we develop a compactness–rigidity argument combined with a streamline circulation method. This allows us to establish the existence of steady Prandtl–Batchelor flows and justify the associated two-sided Prandtl boundary layer expansions without imposing any symmetry condition on the wall perturbations.
	
	\subsection{Mathematical formulation and main result}
	We consider the steady Navier--Stokes system
	\begin{align}
		\begin{cases}
			\textbf{u}^{\epsilon}\cdot\nabla\textbf{u}^{\epsilon}+\nabla p^{\epsilon}-\epsilon^2\triangle\textbf{u}^\epsilon=0,\\[3pt]
			\nabla\cdot\textbf{u}^{\epsilon}=0,\\[3pt]
			\textbf{u}^{\epsilon}|_{y=1}=(\alpha+\delta f(x))\textbf{e}_x,\quad\textbf{u}^{\epsilon}|_{y=0}=(-\beta+\delta g(x))\textbf{e}_x.
		\end{cases}\label{eq:NS0}
	\end{align}
	on the channel $\Omega:=[0,2\pi]\times[0,1]$, with periodicity imposed in the horizontal variable. Here, $\epsilon>0$ is a small parameter, with $\epsilon^2$
	denoting the viscosity coefficient.
	The parameters $\alpha>0$ and $\beta>0$ are fixed,
	while $\delta>0$ measures the amplitude of the boundary
	perturbation.
	The functions $\textbf{u}^\epsilon$ and $p^\epsilon$
	are the velocity and pressure, respectively, and the pressure is normalized by
	\[
	\int_\Omega p^\epsilon\,dx\,dy=0.
	\]
	$\textbf{\textbf{e}}_x$ is the horizontal unit vector, and $f,g$ are two prescribed smooth $2\pi$-periodic functions.
	
	Let $u^\epsilon(x,y)$ and $v^\epsilon(x,y)$ denote
	the horizontal and vertical components of
	$\textbf{u}^\epsilon$, respectively. Write the equations \eqref{eq:NS0} in component form as follows
	\begin{align}
		\begin{cases}
			(u^\epsilon\partial_{x}+v^\epsilon\partial_{y})u^\epsilon+p_{x}^\epsilon-\epsilon^2\triangle u^\epsilon=0\text{ in }  \Omega,\\[3pt]
			(u^\epsilon\partial_{x}+v^\epsilon\partial_{y})v^\epsilon+p_{y}^\epsilon-\epsilon^2\triangle v^\epsilon=0\text{ in }  \Omega,\\[3pt]
			\partial_{x} u^\epsilon+\partial_{y} v^\epsilon=0\text{ in }\Omega,\\[3pt]
			u^\epsilon(x+2\pi,y)=u^\epsilon(x,y),v^\epsilon(x+2\pi,y)=v^\epsilon(x,y)\text{ in }\Omega,\\[3pt]
			u^\epsilon(x,1)=\alpha+\delta f(x),~v^\epsilon(x,1)=0,x\in[0,2\pi],\\[3pt]
			u^\epsilon(x,0)=-\beta+\delta g(x),~v^\epsilon(x,0)=0,x\in[0,2\pi].\label{eq:NS}
		\end{cases}
	\end{align}
	Formally letting $\epsilon\to0$ in \eqref{eq:NS}
	yields the steady Euler system
	\begin{align}
		\begin{cases}
			(u_e\partial_x+v_e\partial_y)u_e+\partial_{x}p_e=0\text{ in }\Omega,\\[3pt]
			(u_e\partial_x+v_e\partial_y)v_e+\partial_{y}p_e=0\text{ in }\Omega,\\[3pt]
			\partial_{x}u_e+\partial_{y}v_e=0\text{ in }\Omega,\\[3pt]
			u_e(x+2\pi,y)=u_e(x,y),v_e(x+2\pi,y)=v_e(x,y)\text{ in }\Omega,\\[3pt]
			v_e|_{y=0,1}=0,x\in[0,2\pi].\label{eq:Euler}
		\end{cases}
	\end{align}

	We choose the reference Euler flow to be the plane Couette-type flow
	\[
	(u_e(x,y),v_e(x,y),p_e(x,y))=(u_e(y),0,0),
	\]
	where the profile $u_e$ is determined by
	\begin{align}
		\begin{cases}
			u_e''(y)=0, \\[3pt]
u_e(1)=A_+,\\[3pt]
u_e(0)=-A_-,\label{eq:Couette flow}
	\end{cases}
	\end{align}
	where
	$$A_+=\left(\alpha^2+\frac{\alpha\delta}{\pi}\int_{0}^{2\pi}f(x)dx+\frac{\delta^2}{2\pi}\int_{0}^{2\pi}f^2(x)dx\right)^{\frac{1}{2}}$$
	and
	$$A_-=\left(\beta^2-\frac{\beta\delta}{\pi}\int_{0}^{2\pi}g(x)dx+\frac{\delta^2}{2\pi}\int_{0}^{2\pi}g^2(x)dx\right)^{\frac{1}{2}}.$$
	Solving this boundary value problem gives 
	\[
	u_e(y)=ay+b,
	\]  
	where 
	\begin{align}
		a=A_++A_-,~~b=-A_-.\label{eq:uecoefficient}
	\end{align}
	It has constant vorticity $a$ and vanishes at the unique interior level
	\[
	y_*=\frac{A_-}{A_++A_-}\in(0,1).
	\]
	So $u_e(y)=a(y-y_*)$ and the tangential Euler velocity changes sign across $y=y_*$.
	
	This choice is motivated by the 
	Prandtl--Batchelor principle.
	Under its hypotheses, an inviscid limit of the form
	$(u_e(y),0)$ must satisfy the ordinary differential equation
	in \eqref{eq:Couette flow}.
	The boundary value is selected by the Batchelor--Wood
	formula.
	Our aim is to construct steady Navier--Stokes solutions
	with this Euler flow as their leading-order interior profile
	and to justify the associated boundary layer approximation.

	In general, the tangential velocity of the Euler flow does
	not agree with the prescribed viscous boundary data.
	This discrepancy motivates Prandtl's boundary layer theory,
	introduced in 1904, in which a rapidly varying correction
	is added near the boundary to match the wall velocity.
	For studies of the steady Prandtl equations and their
	linearizations, we refer to
	\cite{DM,DIT,GZO,GI3,Iyer,JLY,OS,renardy,SWZ,WZ,WZ2}.
	For related Sobolev stability results concerning shear flows
	and Prandtl expansions, see \cite{CWZ,GMM}.
	
	To describe the leading boundary layer near the upper and lower wall, we introduce the stretched normal variables
	\[
	\zeta=\frac{y-1}{\epsilon},\qquad \eta=\frac{y}{\epsilon},
	\]
	near the upper and lower walls, respectively.
	
	Let $(u_p^{(0)},v_p^{(1)})$ solve
	\begin{align}
		\begin{cases}
			(u_e(1)+u_p^{(0)})\partial_{x}u_p^{(0)}+(v_p^{(1)}-v_p^{(1)}(x,0))\partial_{\zeta}{u_p^{(0)}}-\partial_{\zeta}^2{u_p^{(0)}}=0 \text{ in }\Omega_1,\\[3pt]
			\partial_{x}u_p^{(0)}+\partial_{\zeta}v_p^{(1)}=0\text{ in }\Omega_1,\\[3pt]
			u_p^{(0)}(x+2\pi,\zeta)=u_p^{(0)}(x,\zeta),~v_p^{(1)}(x+2\pi,\zeta)=v_p^{(1)}(x,\zeta)\text{ in }\Omega_1,\\[3pt]
			u_p^{(0)}|_{\zeta=0}=\alpha+\delta{f(x)}-u_e(1),~x\in{[0,2\pi]},\\[3pt]
			\lim\limits_{\zeta\rightarrow{-\infty}}u_p^{(0)}=\lim\limits_{\zeta\rightarrow{-\infty}}v_p^{(1)}=0,~x\in{[0,2\pi]}.\label{eq:up0}
		\end{cases}
	\end{align}
	where \(\Omega_1 = [0, 2\pi] \times (-\infty, 0]\), and let $(\hat{u}_p^{(0)},\hat{v}_p^{(1)})$ solve
	\begin{align}
		\begin{cases}
			(u_e(0)+\hat{u}_p^{(0)})\partial_{x}\hat{u}_p^{(0)}+(\hat{v}_p^{(1)}-\hat{v}_p^{(1)}(x,0))\partial_{\eta}{\hat{u}_p^{(0)}}-\partial_{\eta}^2{\hat{u}_p^{(0)}}=0 \text{ in }\Omega_2,\\[3pt]
			\partial_{x}\hat{u}_p^{(0)}+\partial_{\eta}\hat{v}_p^{(1)}=0\text{ in }\Omega_2,\\[3pt]
			\hat{u}_p^{(0)}(x+2\pi,\eta)=\hat{u}_p^{(0)}(x,\eta),~\hat{v}_p^{(1)}(x+2\pi,\eta)=\hat{v}_p^{(1)}(x,\eta)\text{ in }\Omega_2,\\[3pt]
			\hat{u}_p^{(0)}|_{\eta=0}=-\beta+\delta{g(x)}-u_e(0),~x\in{[0,2\pi]},\\[3pt]
			\lim\limits_{\eta\rightarrow{+\infty}}\hat{u}_p^{(0)}=\lim\limits_{\eta\rightarrow{+\infty}}\hat{v}_p^{(1)}=0,~x\in{[0,2\pi]}.\label{eq:lup0}
		\end{cases}
	\end{align}
	where \(\Omega_2 = [0, 2\pi] \times [0,+\infty)\).  
	
	Their derivation through matched asymptotic expansions and solvability is discussed in the next section.
	
	The main result of the paper can be stated as the following.
	\begin{theorem}\label{thm:main}
		Let $\alpha,\beta>0$ be fixed, and let $f,g$ be
		two smooth $2\pi$-periodic functions.
		There exists constant $\delta_0>0, \epsilon_0>0, C>0$ such that for any $\delta\in (0,\delta_0), \epsilon\in(0,\epsilon_0)$, the system \eqref{eq:NS} admits a solution with periodic velocity $(u^\epsilon,v^\epsilon)\in H^2(\Omega)$ and periodic pressure $p^{\epsilon}$ with zero mean, satisfying the following properties:
		
		\smallskip
		\noindent
		\textup{(i) Uniform velocity approximation:}
		\begin{equation}\label{estimate on velocity}
			\begin{aligned}
				\left\|
				u^\epsilon-u_e-u_p^{(0)}\left(x,\frac{y-1}{\epsilon}\right)-\hat{u}_p^{(0)}\left(x,\frac{y}{\epsilon}\right)
				\right\|_{L^\infty(\Omega)}
				+\|v^\epsilon\|_{L^\infty(\Omega)}
				\leq C\epsilon.
			\end{aligned}
		\end{equation}
		Here, $u_e$ is the Couette profile defined by
		\eqref{eq:Couette flow}, and $u_p^{(0)},\hat{u}_p^{(0)}$ are the horizontal
		component of the solution to \eqref{eq:up0},\eqref{eq:lup0}, respectively.
		
		\smallskip
		\noindent
		\textup{(ii) Interior vorticity convergence:}
		For every $0<y_1<y_2<1$, there exist a constant $C_{y_1,y_2}=C(\alpha,\beta,f,g,\delta,y_1,y_2)>0$ such that 
		\begin{equation}\label{estimate on vorticity}
			\|\omega^\epsilon-a\|_{L^{\infty}([0,2\pi]\times[y_1,y_2])}\leq C_{y_1,y_2}\epsilon,
		\end{equation}
		where $
		\omega^\epsilon
		:=u_y^{\epsilon}-v_x^{\epsilon}
		$
		is the vorticity and $a$ is defined in $\eqref{eq:uecoefficient}$.
	\end{theorem}
	
	\begin{remark}
		The estimate \eqref{estimate on velocity} justifies the
composite approximation consisting of the selected
Euler shear and the two leading-order Prandtl boundary
layers uniformly throughout the channel. In particular,
the interior zero of the Euler shear does not require
an additional leading-order internal boundary layer in
our construction.

Although the tangential Euler velocity changes sign in
the interior, it remains nondegenerate at both physical
boundaries. Consequently, both wall layers are of the
classical Prandtl type. The essential difficulty caused
by the interior degeneracy lies in the linear stability
analysis of the remainder equations.

Together with the interior vorticity estimate
\eqref{estimate on vorticity}, this provides a rigorous
realization of the Prandtl--Batchelor selection
principle in the present setting.
	\end{remark}

\begin{remark}
		We expect that the compactness--rigidity argument
and the streamline circulation method developed in
this paper can also be adapted to Prandtl--Batchelor
flows in an annulus with counter-rotating boundaries.
In this setting, the expected limiting Couette--Taylor
flow vanishes on an interior circle, giving rise to
a loss of coercivity analogous to that encountered
in the present paper. Although the curvature of the
annular domain introduces additional technical
difficulties, we believe that the main ideas of our
approach remain applicable. We leave the detailed
analysis for future work.
	\end{remark}
	
	\subsection{Main difficulties and ideas of the proof}
	
	The proof of Theorem 1.1 consists of the construction of high-order approximate solutions, the establishment of quantitative linear stability estimates, and the solution of the nonlinear remainder equations. Among these steps, the main analytical novelty lies in the linear stability analysis.
	
	\textbf{Step 1. Construction of the approximate solution.}
	We first construct an Euler--Prandtl approximate solution
	\[
	\mathbf u^a=\mathbf u_e^a+\mathbf u_{p,+}^a+\mathbf u_{p,-}^a+\mathbf h,
	\]
	where $\mathbf u_e^a$ is the outer Euler expansion, $\mathbf u_{p,+}^a$ and $\mathbf u_{p,-}^a$ are the localized upper and lower boundary layer corrections, and $\mathbf h$ accounts for the necessary divergence and matching corrections.
	
	A useful structural property of the outer expansion is that it can be chosen to be harmonic, with spatially constant vorticity. This property is preserved throughout the finite Euler expansion and plays a crucial role in the subsequent interior stability analysis.
	
	The approximate solution satisfies the prescribed
boundary conditions and the incompressibility constraint.
By carrying out the expansion to sufficiently high order,
we obtain a residual $R^a$ satisfying
\[
\|R^a\|_{L^2(\Omega)}
\le C\epsilon^8.
\]
This high-order smallness compensates for the polynomial
losses in the subsequent linear and nonlinear estimates.
	
	\textbf{Step 2. Linear stability in the presence of an interior zero.}
	Writing $\mathbf u^\epsilon=\mathbf u^a+\mathbf u$ and $p^\epsilon=p^a+p$, we obtain the linearized equation
	\[
	-\epsilon^2\Delta\mathbf u+(\mathbf u^a\cdot\nabla)\mathbf u+(\mathbf u\cdot\nabla)\mathbf u^a+\nabla p=\mathbf F.
	\]
	The main objective is to establish the estimate
	\begin{equation}\label{eq:intro-L2}
		\|\mathbf u\|_{L^2(\Omega)}\le C\epsilon^{-2}\|\mathbf F\|_{L^2(\Omega)}
	\end{equation}
	for divergence-free solutions with homogeneous no-slip boundary conditions.
	
	For a background shear of fixed sign, a suitable positivity argument often provides direct control of the tangential derivatives. In the present problem, however, the limiting shear changes sign at $y=y^*$. The associated quadratic form is therefore not globally positive, and the usual argument cannot be applied directly.
	
	We overcome this difficulty by developing a compactness--rigidity argument based on four complementary ingredients.
	
	First, we establish signed positivity estimates in the two wall regions. The signs of the multipliers are chosen in accordance with the directions of the background flow, so that coercivity is recovered separately near each wall. Combined with the basic energy estimate, these estimates control the tangential derivatives and the viscosity-weighted normal derivatives of the remainder.
	
	Second, we derive interior estimates from the vorticity equation. The outer approximate flow has constant vorticity, which eliminates the coupling with the gradient of the outer vorticity. By constructing cutoff functions that are constant along the outer streamlines, we also eliminate the leading transport contribution in the localized vorticity energy identity. This yields interior $H^2$ estimates without dividing by the tangential Euler velocity and, in particular, without introducing singular coefficients at the interior zero.
	
	Third, we argue by contradiction. If the desired $L^2$ estimate were false, there would exist a sequence of normalized approximate kernel elements satisfying
	\[
	\|\mathbf u_\epsilon\|_{L^2}=1,\qquad \|\mathbf F_\epsilon\|_{L^2}=o(\epsilon^2).
	\]
	The interior estimates yield sufficient compactness to identify the possible limiting modes. Passing to the limit in the vorticity equation shows that the limiting normal velocity vanishes and that the tangential velocity depends only on the normal variable. Suitable streamline-dependent tests then imply that the limiting tangential velocity is affine.
	
	The homogeneous no-slip condition cannot be transferred directly to this affine limit, since the available compactness is only local in the interior. To recover the missing boundary information, we construct adjoint test functions adapted to the Euler--Prandtl composite flow. Their leading boundary layer contributions cancel precisely through the Prandtl equations. The resulting limiting identities force the affine velocity to vanish at both walls, thereby eliminating all affine limiting modes.
	
	Finally, the vanishing of the interior limit does not, by itself, exclude concentration of the normalized $L^2$ mass near the walls. The remaining obstruction is associated with the horizontal mean velocity, which cannot be controlled by a tangential Poincar\'e inequality. We introduce the circulation of the remainder velocity along the streamlines of the approximate flow in each wall region. This circulation satisfies a scalar second-order ordinary differential equation. An explicit Green-kernel representation, together with weighted Hardy and Volterra estimates, rules out concentration of the horizontal mean mode near either wall.
	
	Combining the interior rigidity with the exclusion of wall concentration yields strong convergence of the normalized remainder to zero in the entire channel, contradicting its normalization. This establishes the desired $L^2$ estimate.
	
	The basic energy estimate and the Stokes estimates subsequently give the quantitative bound
	\[
	\|\mathbf u\|_2+\|\partial_x\mathbf u\|_2+\epsilon\|\partial_y\mathbf u\|_2+\epsilon^2\|\nabla^2\mathbf u\|_2\le C\epsilon^{-2}\|\mathbf F\|_2.
	\]
	The linear solvability follows from this estimate and a Fredholm argument.
	
	\textbf{Step 3. Solvability of the nonlinear remainder equations and the vanishing-viscosity limit.}
	With the linear inverse established, the nonlinear remainder equations are solved by a contraction mapping argument in a viscosity-weighted Sobolev space. The high-order smallness of the approximate residual compensates for the polynomial loss in the linear inverse and the derivative losses in the quadratic nonlinear terms.
	
	The resulting remainder is sufficiently small to justify the composite Euler--Prandtl approximation uniformly in the channel. Higher-order estimates and the rapid decay of the boundary layer profiles then yield the local uniform convergence of the vorticity.
	
	The essential feature of the proof is that global linear stability is obtained without requiring a globally positive tangential background velocity. Instead, the combination of local coercivity, interior compactness, adjoint rigidity, and circulation control provides a mechanism for treating the loss of coercivity caused by an interior zero of the selected Euler shear.
	
\subsection{Organization of the paper}
The remainder of the paper is organized as follows. In Section 2, we construct high-order Euler--Prandtl approximate solutions and establish the corresponding residual estimates. Section 3 is devoted to the linear stability analysis, including the signed near-wall estimates, the interior vorticity estimates, the compactness--rigidity argument, and the streamline circulation estimates. In Section 4, we solve the nonlinear remainder equations and complete the proof of Theorem 1.1. The Appendices contain a formal derivation of the Prandtl--Batchelor selection principle in a periodic channel and the construction of the divergence correction.

	\section{Construction of Approximate Solutions}
	In this section, we construct an approximate solution to $(\ref{eq:NS})$ by the matched asymptotic expansion. 
	
\subsection{Euler expansions}
	Away from the boundary, we make the following Euler expansions 
	\begin{align}
		\begin{cases}
			u^{\epsilon}(x,y)=u_e^{(0)}(x,y)+\sum\limits_{i=1}^{k}\epsilon^i u_e^{(i)}(x,y)+\cdot\cdot\cdot,\\
			v^{\epsilon}(x,y)=v_e^{(0)}(x,y)+\sum\limits_{i=1}^{k}\epsilon^i v_e^{(i)}(x,y)+\cdot\cdot\cdot,\\
			p^{\epsilon}(x,y)=p_e^{(0)}(x,y)+\sum\limits_{i=1}^{k}\epsilon^i p_e^{(i)}(x,y)+\cdot\cdot\cdot.\label{eq:EulerExpansions}
		\end{cases}
	\end{align}
\subsubsection{Equation for $(u_e^{(0)},v_e^{(0)},p_e^{(0)})$}
	Putting $(\ref{eq:EulerExpansions})$ into $\eqref{eq:Euler}$ and collecting the $\epsilon^0$ order terms, we obtain the following steady Euler equations
	$$
	\begin{cases}
		(u_e^{(0)}\partial_x+v_e^{(0)}\partial_y)u_e^{(0)}+\partial_{x}p_e^{(0)}=0\text{ in }\Omega,\\[3pt]
		(u_e^{(0)}\partial_x+v_e^{(0)}\partial_y)v_e^{(0)}+\partial_{y}p_e^{(0)}=0\text{ in }\Omega,\\[3pt]
		\partial_{x}u_e^{(0)}+\partial_{y}v_e^{(0)}=0\text{ in }\Omega,\\[3pt]
		u_e^{(0)}(x+2\pi,y)=u_e^{(0)}(x,y),v_e^{(0)}(x+2\pi,y)=v_e^{(0)}(x,y)\text{ in }\Omega.
	\end{cases} 
	$$
	In the present paper, we take $(u_e^{(0)},v_e^{(0)})$ to be shear flow $(u_e(y),0)$, the associated pressure $p_e^{(0)}=0$.

\subsubsection{Equations for $(u_e^{(k)},v_e^{(k)},p_e^{(k)}),k=1,2,3,4,5,6,7$}
	
 Collecting $\epsilon^k$ order terms, we obtain the following steady linearized Euler equations
	\begin{align}
		\begin{cases}
			u_e(y)\partial_xu_e^{(k)}+v_e^{(k)}u'_e(y)+\partial_xp_e^{(k)}=f_e^{(k)},\\[3pt]
			u_e(y)\partial_xv_e^{(k)}+\partial_yp_e^{(k)}=g_e^{(k)},\\[3pt]
			\partial_xu_e^{(k)}+\partial_yv_e^{(k)}=0,\\[3pt]
			 ~u_e^{(k)}(x+2\pi,y)=u_e^{(k)}(x,y),~v_e^{(k)}(x+2\pi,y)=v_e^{(k)}(x,y)\label{eq:uei}
		\end{cases}
	\end{align}
	where
	\begin{align*}
		f_e^{(1)}=&0,\\
		f_e^{(2)}=&-u_e^{(1)}\partial_xu_e^{(1)}-v_e^{(1)}\partial_yu_e^{(1)}+u_e''(y),\\
		f_e^{(k)}=&-\sum_{i=1}^{k-1}u_e^{(i)}\partial_xu_e^{(k-i)}-\sum_{i=1}^{k-1}v_e^{(i)}\partial_yu_e^{(k-i)}+\Delta u_e^{(k-2)},k\geq3.
	\end{align*}
	and
	\begin{align*}
		g_e^{(1)}=&0,\\
		g_e^{(2)}=&-u_e^{(1)}\partial_xv_e^{(1)}-v_e^{(1)}\partial_yv_e^{(1)},\\
		g_e^{(k)}=&-\sum_{i=1}^{k-1}u_e^{(i)}\partial_xv_e^{(k-i)}-\sum_{i=1}^{k-1}v_e^{(i)}\partial_yv_e^{(k-i)}+\Delta v_e^{(k-2)},k\geq3.
	\end{align*}
which are equipped with the boundary conditions
\begin{align}
 v_e^{(k)}|_{y=1}=-v_p^{(k)}|_{\zeta=0},~v_e^{(k)}|_{y=0}=-\hat{v}_p^{(k)}|_{\eta=0},\label{outer-1 order-bc}
\end{align}
where $v_p^{(k)}, \hat{v}_p^{(k)}$ is defined in \eqref{eq:upi} and  \eqref{eq:^upi} respectively in next subsection.
	
\subsection{Prandtl expansions near $y=1$}

We  introduce the scaled variable  $\zeta=\frac{y-1}{\epsilon}\in (-\infty,0]$ and make the upper-wall boundary layer expansions
	\begin{align}
		\begin{cases}
			u^{\epsilon}(x,y)=u_e^{(0)}(x,y)+u_p^{(0)}(x,\zeta)+\sum\limits_{i=1}^{k}\epsilon^i (u_e^{(i)}(x,y)+u_p^{(i)}(x,\zeta))+\cdot\cdot\cdot,\\
			v^{\epsilon}(x,y)=v_e^{(0)}(x,y)+v_p^{(0)}(x,\zeta)+\sum\limits_{i=1}^{k}\epsilon^i (v_e^{(i)}(x,y)+v_p^{(i)}(x,\zeta))+\cdot\cdot\cdot,\\
			p^{\epsilon}(x,y)=p_e^{(0)}(x,y)+p_p^{(0)}(x,\zeta)+\sum\limits_{i=1}^{k}\epsilon^i (p_e^{(i)}(x,y)+p_p^{(i)}(x,\zeta))+\cdot\cdot\cdot.\label{eq:UpperBoundaryExpansions}\\
		\end{cases}
	\end{align}
where as $\zeta\rightarrow -\infty$
\begin{align}
 \partial_x^l\partial_\zeta^mv_p^{(k)}(x,\zeta)\rightarrow 0, \ \partial_x^l\partial_\zeta^mp_p^{(k)}(x,\zeta)\rightarrow 0,\label{matching condition}
\end{align}
here $l,m\geq 0, k=0,1,\cdots,$ and satisfy the following boundary conditions
\begin{align*}
&u_e^{(0)}(x,1)+u_p^{(0)}(x,0)=\alpha+\delta f(x), \ u_e^{(k)}(x,1)+u_p^{(k)}(x,0)=0,\ k\geq 1,\\
&v_e^{(k)}(x,1)+v_p^{(k)}(x,0)=0, k\geq 0, \ \lim_{\zeta\rightarrow -\infty}(u_p^{(0)}(x,\zeta),\partial_{\zeta}u_p^{(k)}(x,\zeta))=(0,0), k\geq1.
\end{align*}

	\subsubsection{Equations for $(v_p^{(0)},p_p^{(0)})$}
	
	Putting (\ref{eq:UpperBoundaryExpansions})
	into $\eqref{eq:NS}$ and collecting $\epsilon^{-1}$ order terms, combining (\ref{matching condition}), we obtain the following equations
	$$
	\begin{cases}
		v_p^{(0)}\partial_{\zeta}v_p^{(0)}+\partial_{\zeta}p_p^{(0)}=0,\\[3pt]
		\partial_{\zeta}v_p^{(0)}=0,\\[3pt]
		v_p^{(0)}(x+2\pi,\zeta)=v_p^{(0)}(x,\zeta),~p_p^{(0)}(x+2\pi,\zeta)=p_p^{(0)}(x,\zeta),\\[3pt]
		\lim\limits_{\zeta\rightarrow-\infty}v_p^{(0)}=\lim\limits_{\zeta\rightarrow-\infty}p_p^{(0)}=0.
	\end{cases}
	$$
	From this, we obtain $v_p^{(0)}=p_p^{(0)}=0$.
	
	\subsubsection{Equations for $(u_p^{(0)},v_p^{(1)},p_p^{(1)})$}
	
	Collecting $\epsilon^0$ order terms, we obtain the following equations
\begin{align}
	\begin{cases}
		(u_e(1)+u_p^{(0)})\partial_x{u_p^{(0)}}+(v_p^{(1)}-v_p^{(1)}(x,0))\partial_{\zeta}u_p^{(0)}-\partial_{\zeta}^2u_p^{(0)}=0,\\[3pt]
		\partial_{\zeta}p_p^{(1)}=0,\\[3pt]
		\partial_xu_p^{(0)}+\partial_{\zeta}v_p^{(1)}=0,\\[3pt]
		u_p^{(0)}(x+2\pi,\zeta)=u_p^{(0)}(x,\zeta),v_p^{(1)}(x+2\pi,\zeta)=v_p^{(1)}(x,\zeta),\\[3pt] u_p^{(0)}|_{\zeta=0}=\alpha+\delta{f(x)}-u_e(1),\lim\limits_{\zeta\rightarrow-\infty}(u_p^{(0)}, v_p^{(1)}, p_p^{(1)})=(0,0,0).\nonumber
	\end{cases}
\end{align}
	From this, we otain $p_p^{(1)}=0$ and the steady Prandtl equation \eqref{eq:up0}. 

\subsubsection{Equations for $(u_p^{(k)},v_p^{(k+1)},p_p^{(k+1)}),k=1,2,3,4,5,6$}
	
Collecting the $\epsilon^k$ order terms, we obtain that
	\begin{align}
		\begin{cases} (u_e(1)+u_p^{(0)})\partial_xu_p^{(k)}+u_p^{(k)}\partial_xu_p^{(0)}+(v_e^{(1)}(x,1)+v_p^{(1)})\partial_{\zeta}u_p^{(k)}\\
			+(v_p^{(k+1)}-v_p^{(k+1)}(x,0))\partial_{\zeta}u_p^{(0)}-\partial_{\zeta}^2u_p^{(k)}=f_{\text{upper}}^{(k)},\\[3pt]
			\partial_xu_p^{(k)}+\partial_{\zeta}v_p^{(k+1)}=0,\\[3pt]
			u_p^{(k)}(x+2\pi,\zeta)=u_p^{(k)}(x,\zeta),v_p^{(k+1)}(x+2\pi,\zeta)=v_p^{(k+1)}(x,\zeta),\\[3pt]
u_p^{(k)}|_{\zeta=0}=-u_e^{(k)}|_{y=1},\lim\limits_{\zeta\rightarrow-\infty}\partial_{\zeta}u_p^{(k)}=\lim\limits_{\zeta\rightarrow-\infty}v_p^{(k+1)}=0, \label{eq:upi}
		\end{cases}
	\end{align}
where
\begin{align*}
	f_{\text{upper}}^{(1)}=&(\zeta\partial_{\zeta}u_p^{(0)}-u_p^{(0)})\partial_xu_e^{(1)}(x,1)-(u_e^{(1)}(x,1)+a\zeta)\partial_xu_p^{(0)}-av_p^{(1)},\\
	f^{(k)}_{\text{upper}}=&\partial_{xx}u_p^{(k-2)}-\partial_xp_p^{(k)}
	-u_p^{(0)}\partial_x\bar a_k-\bar a_k\partial_xu_p^{(0)}
	-v_p^{(1)}\partial_\zeta\bar a_k\\[5pt]
	&-\sum_{\substack{i+j=k\\ i,j\ge1}}\Big[u_p^{(i)}\partial_x\big(u_p^{(j)}+\bar a_j\big)
	+\bar a_i\,\partial_xu_p^{(j)}\Big]\\[5pt]
	&-\sum_{\substack{i+j=k+1\\ i\ge2,\ j\ge1}}\Big[v_p^{(i)}\partial_\zeta\big(u_p^{(j)}+\bar a_j\big)
	+\bar b_i\,\partial_\zeta u_p^{(j)}\Big]
	-\big(\bar b_{k+1}-v_e^{(k+1)}(x,1)\big)\partial_\zeta u_p^{(0)},k\geq2
\end{align*}
and $\bar a_1=a\zeta+u_e^{(1)}(x,1)$, for $k\ge2$,
\begin{align*}
	\bar a_k&=u_e^{(k)}(x,1)+\zeta\partial_yu_e^{(k-1)}(x,1)+\frac{\zeta^2}{2}\partial_{yy}u_e^{(k-2)}(x,1)
	+\cdots+\frac{\zeta^{k-1}}{(k-1)!}\partial_y^{k-1}u_e^{(1)}(x,1),\\
	\bar b_k&=v_e^{(k)}(x,1)+\zeta\partial_yv_e^{(k-1)}(x,1)+\frac{\zeta^2}{2}\partial_{yy}v_e^{(k-2)}(x,1)
	+\cdots+\frac{\zeta^{k-1}}{(k-1)!}\partial_y^{k-1}v_e^{(1)}(x,1).
\end{align*}
Moreover, 
$$\partial_{\zeta}p_p^{(k+1)}=g_{\operatorname{upper}}^{(k)},\lim\limits_{\zeta\to-\infty}p_p^{(k+1)}=0,$$
where
\begin{align*}
	g_{\operatorname{upper}}^{(1)}=&\partial_{\zeta}^2v_p^{(1)}-(u_e(1)+u_p^{(0)})\partial_xv_p^{(1)}-(v_e^{(1)}(x,1)+v_p^{(1)})\partial_{\zeta}v_p^{(1)}-u_p^{(0)}\partial_xv_e^{(1)}(x,1),\\
	g^{(k)}_{\text{upper}}=&\partial_{xx}v_p^{(k-2)}+\partial_{\zeta\zeta}v_p^{(k)}
	-\big(u_e(1)+u_p^{(0)}\big)\partial_xv_p^{(k)}
	-\big(v_e^{(1)}(x,1)+v_p^{(1)}\big)\partial_\zeta v_p^{(k)}\\[5pt]
	&-u_p^{(0)}\partial_x\bar b_k-v_p^{(1)}\partial_\zeta\bar b_k-\sum_{\substack{i+j=k\\ i,j\ge1}}\Big[u_p^{(i)}\partial_x\big(v_p^{(j)}+\bar b_j\big)
	+\bar a_i\,\partial_xv_p^{(j)}\Big]\\[5pt]
	&-\sum_{\substack{i+j=k+1\\ i\ge2,\ j\ge1}}\Big[v_p^{(i)}\partial_\zeta\big(v_p^{(j)}+\bar b_j\big)
	+\bar b_i\,\partial_\zeta v_p^{(j)}\Big],k\geq2.
\end{align*}

\subsubsection{Equations for $(u_p^{(7)},v_p^{(8)},p_p^{(8)})$}

Collecting the $\epsilon^7$ order terms, we obtain that
\begin{align}
	\begin{cases}
		(u_e(1)+u_p^{(0)})\partial_xu_p^{(7)}+u_p^{(7)}\partial_xu_p^{(0)}+(v_p^{(1)}+v_e^{(1)}(x,1))\partial_{\zeta}u_p^{(7)}
		\\[5pt]
		+(v_p^{(8)}-v_p^{(8)}(x,0))\partial_{\zeta}u_p^{(0)}-\partial_{\zeta}^2u_p^{(7)}=f_{\text{upper}}^{(7)},\\[5pt]
		\partial_xu_p^{(7)}+\partial_{\zeta}v_p^{(8)}=0,\\[5pt]
		u_p^{(7)}(x+2\pi,\zeta)=u_p^{(7)}(x,\zeta),~v_p^{(8)}(x+2\pi,\zeta)=v_p^{(8)}(x,\zeta),\\[5pt]			u_p^{(7)}|_{\zeta=0}=-u_e^{(7)}|_{y=1},~\lim\limits_{\zeta\rightarrow-\infty}\partial_{\zeta}u_p^{(7)}=v_p^{(8)}(x,0)=0\label{eq:up end}
	\end{cases}
\end{align}
where the expression for \(f^{(7)}_{\text{upper}}\) can be obtained by substituting \(k=7\) into \(f^{(k)}_{\text{upper}}\). Moreover,
$$\partial_{\zeta}p_p^{(8)}=g_{\operatorname{upper}}^{(7)},~p_p^{(8)}(x,0)=0,$$
where the expression for \(g^{(7)}_{\text{upper}}\) can be obtained by substituting \(k=7\) into \(g^{(k)}_{\text{upper}}\).

\subsection{Prandtl expansions near $y=0$}

We make the lower-wall boundary layer expansions
\begin{align}
		\begin{cases}
			u^{\epsilon}(x,y)=u_e^{(0)}(x,y)+\hat{u}_p^{(0)}(x,\eta)+\sum\limits_{i=1}^{k}\epsilon^i (u_e^{(i)}(x,y)+\hat{u}_p^{(i)}(x,\eta))+\cdot\cdot\cdot,\\
			v^{\epsilon}(x,y)=v_e^{(0)}(x,y)+\hat{v}_p^{(0)}(x,\eta)+\sum\limits_{i=1}^{k}\epsilon^i (v_e^{(i)}(x,y)+\hat{v}_p^{(i)}(x,\eta))+\cdot\cdot\cdot,\\
			p^{\epsilon}(x,y)=p_e^{(0)}(x,y)+\hat{p}_p^{(0)}(x,\eta)+\sum\limits_{i=1}^{k}\epsilon^i (p_e^{(i)}(x,y)+\hat{p}_p^{(i)}(x,\eta))+\cdot\cdot\cdot\label{eq:LowerBoudaryExpansions}\\
		\end{cases}
	\end{align}
where as $\eta\rightarrow +\infty$
\begin{align*}
 \partial_x^l\partial_\eta^m\widehat{v}_p^{(k)}(x,\eta)\rightarrow 0, \ \partial_x^l\partial_{\eta}^m\widehat{p}_p^{(k)}(x,\eta)\rightarrow 0,
\end{align*}
here $l,m\geq 0, k=0,1,\cdots,$ and  there hold the following boundary conditions
\begin{align*}
&u_e^{(0)}(x,0)+\hat{u}_p^{(0)}(x,0)=-\beta+\delta g(x), \ u_e^{(k)}(x,0)+\hat{u}_p^{(k)}(x,0)=0,k\geq 1,\\
&v_e^{(k)}(x,0)+\hat{v}_p^{(k)}(x,0)=0, k\geq 0, \ \lim_{\eta\rightarrow +\infty}(\hat{u}_p^{(0)}(x,\eta),\partial_\eta\hat{u}_p^{(k)}(x,\eta))=(0,0),k\geq1.
\end{align*}

\subsubsection{Equations for $(\hat{v}_p^{(0)}, \hat{p}_p^{(0)})$} 
Putting (\ref{eq:LowerBoudaryExpansions})
	into $\eqref{eq:NS}$ and collecting $\epsilon^{-1}$ order terms, we obtain $\hat{v}_p^{(0)}=0,\  \hat{p}_p^{(0)}=0$.

\subsubsection{Equations for $(\hat{u}_p^{(0)},\hat{v}_p^{(1)},\hat{p}_p^{(1)})$}
Collecting $\epsilon^{0}$ order terms, we obtain $\hat{p}_p^{(1)}=0$ and the steady Prandtl equation $(\ref{eq:lup0})$.

\subsubsection{Equations for $(\hat{u}_p^{(k)},\hat{v}_p^{(k+1)},\hat{p}_p^{(k+1)}),k=1,2,3,4,5,6$}
	Collecting $\epsilon^k$ order terms, we obtain that
	\begin{align}
		\begin{cases}		
		(u_e(0)+\hat{u}_p^{(0)})\partial_x\hat{u}_p^{(k)}+\hat{u}_p^{(k)}\partial_x\hat{u}_p^{(0)}+(v_e^{(k)}(x,0)
+\hat{v}_p^{(1)})\partial_{\eta}\hat{u}_p^{(k)}\\[5pt]
			+(\hat{v}_p^{(k+1)}-\hat{v}_p^{(k+1)}(x,0))\partial_{\eta}\hat{u}_p^{(0)}-\partial_{\eta}^2\hat{u}_p^{(k)}=f_{\text{lower}}^{(k)},\\[5pt]
			\partial_x\hat{u}_p^{(k)}+\partial_{\eta}\hat{v}_p^{(k+1)}=0,\\[5pt]
			\hat{u}_p^{(k)}(x+2\pi,\eta)=\hat{u}_p^{(k)}(x,\eta),\hat{v}_p^{(k+1)}(x+2\pi,\eta)=\hat{v}_p^{(k+1)}(x,\eta),\\[3pt] \hat{u}_p^{(k)}|_{\eta=0}=-u_e^{(k)}|_{y=0},\lim\limits_{\eta\rightarrow+\infty}\partial_{\eta}\hat{u}_p^{(k)}=\lim\limits_{\eta\rightarrow+\infty}\hat{v}_p^{(k+1)}=0.\label{eq:^upi}
		\end{cases}
	\end{align}
	Moreover,
	$$\partial_{\eta}\hat{p}_p^{(k+1)}=g_{\operatorname{lower}}^{(k)},\lim\limits_{\eta\to+\infty}\hat{p}_p^{(k+1)}=0,$$
	where $f_{\text{lower}}^{(k)},g_{\text{lower}}^{(k)}$ is similar to $f_{\text{upper}}^{(k)},g_{\text{upper}}^{(k)}$, respectively, so we omit the details.
	
\subsubsection{Equations for $(\hat{u}_p^{(7)},\hat{v}_p^{(8)},\hat{p}_p^{(8)})$}

Collecting the $\epsilon^7$ order terms, we obtain that
\begin{align}
	\begin{cases}
		(u_e(0)+\hat{u}_p^{(0)})\partial_x\hat{u}_p^{(7)}+\hat{u}_p^{(7)}\partial_xu_p^{(0)}+(\hat{v}_p^{(1)}+\hat{v}_e^{(1)}(x,0))\partial_{\eta}\hat{u}_p^{(7)}
		\\[5pt]
		+(\hat{v}_p^{(8)}-\hat{v}_p^{(8)}(x,0))\partial_{\eta}\hat{u}_p^{(0)}-\partial_{\eta}^2\hat{u}_p^{(7)}=f_{\text{lower}}^{(7)},\\[5pt]
		\partial_x\hat{u}_p^{(7)}+\partial_{\eta}\hat{v}_p^{(8)}=0,\\[5pt]
		\hat{u}_p^{(7)}(x+2\pi,\eta)=\hat{u}_p^{(7)}(x,\eta),~\hat{v}_p^{(8)}(x+2\pi,\eta)=\hat{v}_p^{(8)}(x,\eta),\\[5pt]
		\hat{u}_p^{(7)}|_{\eta=0}=-u_e^{(7)}|_{y=0},\lim\limits_{\eta\rightarrow+\infty}\partial_{\eta}\hat{u}_p^{(7)}=\hat{v}_p^{(8)}(x,0)=0\label{eq:^up end}
	\end{cases}
\end{align}
and 
$$\partial_{\eta}\hat{p}_p^{(8)}=g_{\text{lower}}^{(7)},~\hat{p}_p^{(8)}(x,0)=0,$$
where $f_{\text{lower}}^{(7)}$ and $g_{\text{lower}}^{(7)}$ are similar to $f_{\text{upper}}^{(7)}$ and $g_{\text{upper}}^{(7)}$,  we omit them.

We will give the solvability of the above Euler equations and Prandtl equations.

\subsection{Solvabilities of Euler equations and Prandtl equations}
\indent

The order in which we solve the equation is as follows

\begin{tikzpicture}
	\node (O) at (0,0) {$(u_e^{(0)},v_e^{(0)})$};
	\node (P0) at (0,1.5) {$(u_p^{(0)},v_p^{(1)})$};
	\node (p0) at (0,-1.5) {$(\hat{u}_p^{(0)},\hat{v}_p^{(1)})$};
	\node (E1) at (2,0) {$(u_e^{(1)},v_e^{(1)})$};
	\node (P1) at (2,1.5) {$(u_p^{(1)},v_p^{(2)})$};
	\node (p1) at (2,-1.5) {$(\hat{u}_p^{(1)},\hat{v}_p^{(2)})$};
	\node (E2) at (4,0) {$(u_e^{(2)},v_e^{(2)})$};
	\node (P2) at (4,1.5) {$(u_p^{(2)},v_p^{(3)})$};
	\node (p2) at (4,-1.5) {$(\hat{u}_p^{(2)},\hat{v}_p^{(3)})$};
	\node (E3) at (6,0) {$(u_e^{(3)},v_e^{(3)})$};
	\node (P3) at (6,1.5) {$(u_p^{(3)},v_p^{(4)})$};
	\node (p3) at (6,-1.5) {$(\hat{u}_p^{(3)},\hat{v}_p^{(4)})$};
	\node (E4) at (8,0) {$(u_e^{(4)},v_e^{(4)})$};
	\node (P4) at (8,1.5) {$(u_p^{(4)},v_p^{(5)})$};
	\node (p4) at (8,-1.5) {$(\hat{u}_p^{(4)},\hat{v}_p^{(5)})$};
	\node (E5) at (10,0) {$(u_e^{(5)},v_e^{(5)})$};
	\node (P5) at (10,1.5) {$(u_p^{(5)},v_p^{(6)})$};
	\node (p5) at (10,-1.5) {$(\hat{u}_p^{(5)},\hat{v}_p^{(6)})$};
	\node (E6) at (12,0) {$(u_e^{(6)},v_e^{(6)})$};
	\node (P6) at (12,1.5) {$(u_p^{(6)},v_p^{(7)})$};
	\node (p6) at (12,-1.5) {$(\hat{u}_p^{(6)},\hat{v}_p^{(7)})$};
	\node (E7) at (14,0) {$(u_e^{(7)},v_e^{(7)})$};
	\node (P7) at (14,1.5) {$(u_p^{(7)},v_p^{(8)})$};
	\node (p7) at (14,-1.5) {$(\hat{u}_p^{(7)},\hat{v}_p^{(8)})$};
	\draw[->] (O) -- node[above] {$upper\quad boundary$} (P0);
	\draw[->] (O) -- node[above] {$lower\quad boundary$} (p0);
	\draw[->] (P0) -- node[above] {} (E1);
	\draw[->] (p0) -- node[above] {} (E1);
	
	\draw[->] (E1) -- node[above] {} (P1);
	\draw[->] (E1) -- node[above] {} (p1);
	\draw[->] (P1) -- node[above] {} (E2);
	\draw[->] (p1) -- node[above] {} (E2);
	
	\draw[->] (E2) -- node[above] {} (P2);
	\draw[->] (E2) -- node[above] {} (p2);
	\draw[->] (P2) -- node[above] {} (E3);
	\draw[->] (p2) -- node[above] {} (E3);
	
	\draw[->] (E3) -- node[above] {} (P3);
	\draw[->] (E3) -- node[above] {} (p3);
	\draw[->] (P3) -- node[above] {} (E4);
	\draw[->] (p3) -- node[above] {} (E4);
	
	\draw[->] (E4) -- node[above] {} (P4);
	\draw[->] (E4) -- node[above] {} (p4);
	\draw[->] (P4) -- node[above] {} (E5);
	\draw[->] (p4) -- node[above] {} (E5);
	
	\draw[->] (E5) -- node[above] {} (P5);
	\draw[->] (E5) -- node[above] {} (p5);
	\draw[->] (P5) -- node[above] {} (E6);
	\draw[->] (p5) -- node[above] {} (E6);
	
	\draw[->] (E6) -- node[above] {} (P6);
	\draw[->] (E6) -- node[above] {} (p6);
	\draw[->] (P6) -- node[above] {} (E7);
	\draw[->] (p6) -- node[above] {} (E7);
	
	\draw[->] (E7) -- node[above] {} (P7);
	\draw[->] (E7) -- node[above] {} (p7);
	
\end{tikzpicture}

We give some explanations for this picture.
\begin{itemize}
	\item $(u_e, v_e)$ denotes the Euler asymptotic expansion away from the boundary, while $(u_p, v_p)$ denotes the asymptotic expansion near the upper boundary and $(\hat{u}_p, \hat{v}_p)$ denotes the asymptotic expansion near the lower boundary. 
	\item The superscript on the Euler and the boundary profiles represents the matching order of $\epsilon$.
	\item The vertical arrows mean the matched boundary condition between $u_e$ and $u_p,\hat{u}_p$,while the slanted arrows mean the matched boundary condition between $v_e$ and $v_p,\hat{v}_p$.
	\item The equations satisfied by $(u_e, v_e)$, $(u_p, v_p)$ and $(\hat{u}_p, \hat{v}_p)$ will be derived by matching the order of $\epsilon$ in (\ref{eq:NS}). The asymptotic expansion will be solved column by column starting from $(u_e^{(0)}, v_e^{(0)}):=(u_e(y), 0)$.
\end{itemize}

\subsubsection{The solvability of Prandtl equations}
\indent
The following Lemma give some necessary conditions for the solvability of Prandtl equations, one can refer to Lemma 2.1 in \cite{FGLT1}.
	\begin{lemma}\label{lem:ue1}
		If the equation $\eqref{eq:up0}$ has solution $(u_p^{(0)},v_p^{(1)})$ which satisfies $$u_e(1)+u_p^{(0)}>0,\|u_p^{(0)}\|_{L^\infty}<+\infty,\forall \zeta\leq0,$$then there holds$$u_e^2(1)=\alpha^2+\frac{\alpha\delta}{\pi}\int_{0}^{2\pi}f(x)dx+\frac{\delta^2}{2\pi}\int_{0}^{2\pi}f^2(x)dx.$$
	\end{lemma}

	One can also refer to Corollary 2.4 in \cite{FGLT1}, which gives the existence of the solution of \eqref{eq:up0}.
	\begin{proposition}\label{prop:up0vp1}
		There exists $\delta_0>0$ such that for any $\delta\in(0,\delta_0)$ and any $j,k,l\in\mathbb{N}$, the equation $\eqref{eq:up0}$ has a unique solution $(u_p^{(0)},v_p^{(1)})$ which satisfies
		\begin{equation}
			\begin{aligned}
				&\int_{-\infty}^{0}\int_{0}^{2\pi}\left|\partial_x^j\partial_{\zeta}^k(u_p^{(0)},v_p^{(1)})\right|^2\langle\zeta\rangle^{2l}dxd\zeta\leq{C_{j,k,l}\delta^2},\nonumber\\
				&\int_{0}^{2\pi}v_p^{(1)}(x,\zeta)dx=0,~\forall\zeta\leq0. \nonumber
			\end{aligned}
		\end{equation}
	\end{proposition}
	
The same results holds for $(\hat{u}_p^{(0)},\hat{v}_p^{(1)})$.
	\begin{lemma}\label{lem:ue0}
		If the equation $\eqref{eq:lup0}$ has solution $(\hat{u}_p^{(0)},\hat{v}_p^{(1)})$ which satisfies $$u_e(0)+\hat{u}_p^{(0)}<0,\|\hat{u}_p^{(0)}\|_{L^\infty}<+\infty,\forall \eta\geq0,$$then there holds$$u_e^2(0)=\beta^2-\frac{\beta\delta}{\pi}\int_{0}^{2\pi}g(x)dx+\frac{\delta^2}{2\pi}\int_{0}^{2\pi}g^2(x)dx.$$
	\end{lemma}
	\begin{proof}
		Let
		$$\bar{u}_p^{(0)}(x,\zeta):=-\hat{u}_p^{(0)}(-x,-\zeta),\bar{v}_p^{(1)}(x,\zeta):=-\hat{v}_p^{(1)}(-x,-\zeta),x\in[0,2\pi],\zeta\le0.$$
		Then $(\bar{u}_p^{(0)},\bar{v}_p^{(1)})$ satisfies
		\begin{align*}
			\begin{cases}
				(-u_e(0)+\bar{u}_p^{(0)})\partial_{x}\bar{u}_p^{(0)}+(\bar{v}_p^{(1)}-\bar{v}_p^{(1)}(x,0))\partial_{\zeta}{\bar{u}_p^{(0)}}-\partial_{\zeta}^2{\bar{u}_p^{(0)}}=0 \text{ in }\Omega_1,\\[3pt]
				\partial_{x}\bar{u}_p^{(0)}+\partial_{\zeta}\bar{v}_p^{(1)}=0\text{ in }\Omega_1,\\[3pt]
				\bar{u}_p^{(0)}(x+2\pi,\zeta)=\bar{u}_p^{(0)}(x,\zeta),~\bar{v}_p^{(1)}(x+2\pi,\zeta)=\bar{v}_p^{(1)}(x,\zeta)\text{ in }\Omega_1,\\[3pt]
				\bar{u}_p^{(0)}|_{\zeta=0}=\beta-\delta{g(-x)}+u_e(0),~x\in{[0,2\pi]},\\[3pt]
				\lim\limits_{\zeta\rightarrow{-\infty}}\bar{u}_p^{(0)}=\lim\limits_{\zeta\rightarrow{-\infty}}\bar{v}_p^{(1)}=0,~x\in{[0,2\pi]}.
			\end{cases}
		\end{align*}
		This equation has the same form of $\eqref{eq:up0},$so the conclusion follows from Lemma $\ref{lem:ue1}.$
	\end{proof}
	\begin{proposition}\label{prop:^up0vp1}
		There exists $\delta_0>0$ such that for any $\delta\in(0,\delta_0)$ and any $j,k,l\in\mathbb{N}$, the equation $\eqref{eq:lup0}$ has a unique solution $(\hat{u}_p^{(0)},\hat{v}_p^{(1)})$ which satisfies 
		\begin{equation}
			\begin{aligned}
&\int_{0}^{+\infty}\int_{0}^{2\pi}\left|\partial_x^j\partial_{\eta}^k(\hat{u}_p^{(0)},\hat{v}_p^{(1)})\right|^2\langle\eta\rangle^{2l}dxd\eta\leq{C_{j,k,l}\delta^2},\nonumber\\
				&\int_{0}^{2\pi}\hat{v}_p^{(1)}(x,\eta)dx=0,~\forall\eta\geq0.\nonumber
			\end{aligned}
		\end{equation}
	\end{proposition}
	
\subsubsection{The solvability of linearized Euler equation for $(u_e^{(1)}, v_e^{(1)}, p_e^{(1)})$ }
\begin{proposition}\label{prop:ue1ve1}
	There exists $\delta_0>0$ such that for any $\delta\in(0,\delta_0)$, the linearized Euler equation $\eqref{eq:uei}$ for $k=1$ has a solution $(u_e^{(1)},v_e^{(1)},p_e^{(1)})$ which satisfies
	\begin{equation}
		\begin{aligned}
			&\|\partial_x^j\partial_y^k(u_e^{(1)},v_e^{(1)})\|_2\leq{C_{j,k}}\delta,~\forall j,k\geq0.\\
			&\Delta u_e^{(1)}=\Delta v_e^{(1)}=0, (x,y)\in\Omega.\\
			&\int_{0}^{2\pi}v_e^{(1)}(x,y)dx=0,~\forall y\in[0,1]. 
		\end{aligned}\label{eq:ue1ve1Estimation}
	\end{equation}
	Besides, the vorticity $\omega_e^{(1)}:=\partial_yu_e^{(1)}-\partial_xv_e^{(1)}$ is a constant.
\end{proposition}
\begin{proof}
	Eliminating $p_e^{(1)}$ from $\eqref{eq:uei}$, we obtain that
	\begin{align}
		\begin{cases}
			-u_e(y)\Delta{v_e^{(1)}}=0,\\[3pt]
			v_e^{(1)}|_{y=1}=-v_p^{(1)}|_{\zeta=0},~v_e^{(1)}|_{y=0}=-\hat{v}_p^{(1)}|_{\eta=0}.\label{eq:ve1}
		\end{cases}
	\end{align}
	Expand \(-v_p^{(1)}(x, 0)\) and \(-\hat{v}_p^{(1)}(x, 0)\) as  
	\begin{align*}
		-v_p^{(1)}(x, 0) &= \sum_{n=1}^{\infty} \big( a_{n1} \cos nx + b_{n1} \sin nx \big), \\
		-\hat{v}_p^{(1)}(x, 0) &= \sum_{n=1}^{\infty} \big( c_{n1} \cos nx + d_{n1} \sin nx \big).
	\end{align*}
	We define  
	\begin{align*}
		u_e^{(1)}(x, y) &= \sum_{n=1}^{\infty} \left( \frac{\partial_yb_n^{(1)}(y)}{n} \cos nx - \frac{\partial_ya_n^{(1)}(y)}{n} \sin nx \right), \\
		v_e^{(1)}(x, y) &= \sum_{n=1}^{\infty} \big( a_n^{(1)}(y) \cos nx + b_n^{(1)}(y) \sin nx \big), \\
		p_e^{(1)}(x, y) &= -\int_{0}^{x} \big( u_e(y) \partial_x u_e^{(1)}(\bar x, y) + u_e'(y) v_e^{(1)}(\bar x, y) \big) d\bar x + \varphi_1(y),
	\end{align*}
	where  
	\begin{align*}
		\frac{d^2}{dy^2}a_n^{(1)}(y) -n^2a_n^{(1)} = 0,a_n(1)=a_{n1},a_n(0)=c_{n1},\\
		\frac{d^2}{dy^2}b_n^{(1)}(y) -n^2b_n^{(1)} = 0,b_n(1)=b_{n1},b_n(0)=d_{n1},
	\end{align*}
	and \(\varphi_1(y)\) satisfies  
	\begin{align*}
		\varphi_1'(y)=-u_e(y)\partial_xv_e^{(1)}(0,y).
	\end{align*}
	Hence, direct calculation shows that $(u_e^{(1)}, v_e^{(1)},p_e^{(1)})$ satisfies the equation \eqref{eq:uei} and $\eqref{eq:ue1ve1Estimation}_{2,3}$. $\eqref{eq:ue1ve1Estimation}_1$ follows from Proposition $\ref{prop:up0vp1}$ and Proposition $\ref{prop:^up0vp1}$. $\omega_e^{(1)}$ is a constant since $\nabla\omega_e^{(1)}=(-\Delta v_e^{(1)},\Delta u_e^{(1)})=(0,0).$
\end{proof}

	\subsubsection{The solvability of linearized Prandtl equations for  $(u_p^{(1)},v_p^{(2)})$ and $(\hat{u}_p^{(1)},\hat{v}_p^{(2)})$}
\indent
	
Now we present the solvability of linearized Prandtl equation which can be found in \cite{FGLT2}.
	\begin{proposition}\label{prop:up1vp2Estimation}
		There exists $\delta_0>0$ such that for any $\delta\in(0,\delta_0)$, the equation $\eqref{eq:upi}$ for $k=1$ has a solution $(u_p^{(1)},v_p^{(2)})$ which satisfies
		\begin{align*}
			&\int_{-\infty}^{0}\int_{0}^{2\pi}|\partial_x^j\partial_{\zeta}^k(u_p^{(1)}-A_1,v_p^{(2)})|^2\langle\zeta\rangle^{2l}dxd\zeta\leq{C_{j,k,l}}\delta^2,\\
			&\int_{0}^{2\pi}v_p^{(2)}(x,\zeta)dx=0,~\forall \zeta\leq0,
		\end{align*}
		where $A_1=\lim\limits_{\zeta\rightarrow-\infty}u_p^{(1)}(x,\zeta)$ is a constant which satisfies $|A_1|\leq{C\delta}.$
	\end{proposition}
	
	\begin{proposition}\label{prop:^up1vp2Estimation}
		There exists $\delta_0>0$ such that for any $\delta\in(0,\delta_0)$, the equation \eqref{eq:^upi} for $k=1$ has a solution $(\hat{u}_p^{(1)},\hat{v}_p^{(2)})$ which satisfies
		\begin{align*}
			&\int_{0}^{+\infty}\int_{0}^{2\pi}|\partial_x^j\partial_{\eta}^k(\hat{u}_p^{(1)}-\hat{A}_1,\hat{v}_p^{(2)})|^2\langle\eta\rangle^{2l}dxd\eta\leq{C_{j,k,l}}\delta^2,\\
			&\int_{0}^{2\pi}\hat{v}_p^{(2)}(x,\eta)dx=0,~\forall \eta\geq0,
		\end{align*}
		where $\hat{A}_1=\lim\limits_{\eta\rightarrow+\infty}\hat{u}_p^{(1)}(x,\eta)$ is a constant which satisfies $|\hat{A}_1|\leq{C\delta}.$
	\end{proposition}
	\begin{proof}
		Applying the same reflection in the proof of Lemma $\ref{lem:ue0}$,for $x\in[0,2\pi],\zeta\le0,$we denote
		$$(\bar{u}_p^{(0)},\bar{u}_p^{(1)})(x,\zeta):=-(\hat{u}_p^{(0)},\hat{u}_p^{(1)})(-x,-\zeta),(\bar{v}_p^{(1)},\bar{v}_p^{(2)})(x,\zeta):=-(\hat{v}_p^{(1)},\hat{v}_p^{(2)})(-x,-\zeta).$$
		Then the equation satisfied by $(\bar{u}_p^{(1)},\bar{v}_p^{(2)})$ has the same form of $\eqref{eq:upi}$ for $k=1$ and the conclusion follows from Proposition $\ref{prop:up1vp2Estimation}.$
	\end{proof}
	Next we do some modifications. Let
	\begin{equation}
		\begin{aligned}
			\tilde{u}_e^{(1)}&=u_e^{(1)}+yA_1+(1-y)\hat{A}_1,~\tilde{v}_e^{(1)}=v_e^{(1)},~\tilde{p}_e^{(1)}=p_e^{(1)},\\[3pt]
			\tilde{u}_p^{(1)}&=u_p^{(1)}-A_1,~\tilde{v}_p^{(2)}=v_p^{(2)},~\tilde{p}_p^{(2)}=p_p^{(2)},\\[3pt]
			\tilde{\hat{u}}_p^{(1)}&=\hat{u}_p^{(1)}-\hat{A}_1,~\tilde{\hat{v}}_p^{(2)}=\hat{v}_p^{(2)},~\tilde{\hat{p}}_p^{(2)}=\hat{p}_p^{(2)},\label{eq:ue1modification}
		\end{aligned}
	\end{equation}
	Then  $(\tilde{u}_p^{(1)},\tilde{v}_p^{(2)},\tilde{p}_p^{(2)})$ satisfies 
	\begin{align}
		\begin{cases}			(u_e(1)+u_p^{(0)})\partial_x\tilde{u}_p^{(1)}+\tilde{u}_p^{(1)}\partial_xu_p^{(0)}+(\tilde{v}_e^{(1)}(x,1)+v_p^{(1)})\partial_{\zeta}\tilde{u}_p^{(1)}\\[5pt]
			+(\tilde{v}_p^{(2)}
			-\tilde{v}_p^{(2)}(x,0))\partial_{\zeta}u_p^{(0)}-\partial_{\zeta}^2\tilde{u}_p^{(1)}=\tilde{f}_{\text{upper}}^{(1)},\\[5pt]
			\partial_x\tilde{u}_p^{(1)}+\partial_{\zeta}\tilde{v}_p^{(2)}=0,\\[5pt]			\tilde{u}_p^{(1)}|_{\zeta=0}=-\tilde{u}_e^{(1)}|_{y=1},\lim\limits_{\zeta\rightarrow-\infty}\partial_{\zeta}\tilde{u}_p^{(1)}
			=\lim\limits_{\zeta\rightarrow-\infty}\tilde{v}_p^{(2)}=0\label{eq:up1vp2modification}
		\end{cases}
	\end{align}
	and $$(u_e(1)+u_p^{(0)})\partial_xv_p^{(1)}+(\tilde{v}_e(x,1)+v_p^{(1)})\partial_{\zeta}v_p^{(1)}+u_p^{(0)}\partial_x\tilde{v}_e^{(1)}(x,1)+\partial_{\zeta}\tilde{p}_p^{(2)}-\partial_{\zeta}^2v_p^{(1)}=0,\lim\limits_{\zeta\rightarrow-\infty}\tilde{p}_p^{(2)}=0,$$
	where $\tilde{f}_{\text{upper}}^{(1)}=f_{\text{upper}}^{(1)}+A_1\partial_xu_p^{(0)}$. Similarly, $(\tilde{\hat{u}}_p^{(1)},\tilde{\hat{v}}_p^{(2)})$ satisfies equation of the same form and so we omit.Moreover, $(\tilde{u}_e^{(1)},\tilde{v}_e^{(1)},\tilde{p}_e^{(1)})$ satisfies
	\begin{align}
		\begin{cases}
			u_e(y)\partial_x\tilde{u}_e^{(1)}+\tilde{v}_e^{(1)}u_e'(y)+\partial_x\tilde{p}_e^{(1)}=0,\\[3pt]
			u_e(y)\partial_x\tilde{v}_e^{(1)}+\partial_y\tilde{p}_e^{(1)}=0,\\[3pt]
			\partial_x\tilde{u}_e^{(1)}+\partial_y\tilde{v}_e^{(1)}=0,\\[3pt]
			\tilde{u}_e^{(1)}(x+2\pi,y)=\tilde{u}_e^{(1)}(x,y),~\tilde{v}_e^{(1)}(x+2\pi,y)=\tilde{v}_e^{(1)}(x,y),\\[3pt]
			\tilde{v}_e^{(1)}|_{y=1}=-v_p^{(1)}|_{\zeta=0},~\tilde{v}_e^{(1)}|_{y=0}=-\hat{v}_p^{(1)}|_{\eta=0}.\label{eq:ue1modeq}
		\end{cases}
	\end{align} 
	and
	\begin{align*}
		&\|\partial_x^j\partial_y^k(\tilde{u}_e^{(1)},\tilde{v}_e^{(1)})\|_2\leq{C_{j,k}}\delta,~\forall j,k\geq0.\\
		&\Delta \tilde{u}_e^{(1)}=\Delta \tilde{v}_e^{(1)}=0, (x,y)\in\Omega.\\
		&\int_{0}^{2\pi}\tilde{v}_e^{(1)}(x,y)dx=0,~\forall y\in[0,1]. 
	\end{align*}
	and has a constant vorticity.
	Besides $\tilde{u}_p^{(1)},\tilde{\hat{u}}_p^{(1)}$ have zero limit: 
	$$\lim\limits_{\zeta\to-\infty}\tilde{u}_p^{(1)}=\lim\limits_{\eta\to+\infty}\tilde{\hat{u}}_p^{(1)}=0.$$
	For simplicity, we omit the tilde superscripts.
	
	\subsubsection{The solvability of linearized Euler equations for $(u_e^{(k)}, v_e^{(k)}, p_e^{(k)}),k=2,3,4,5,6$}

\indent
	
	By substituting $\eqref{eq:EulerExpansions}$ into $\eqref{eq:NS}$ and collecting $\epsilon^k$ order terms, we obtain $\eqref{eq:uei}$.
	The solvability of $(u_e^{(k)}, v_e^{(k)}, p_e^{(k)})$ is similar to $(u_e^{(1)},v_e^{(1)},p_e^{(1)})$. We have the following Proposition. 
	\begin{proposition}\label{prop:uekvek}
		There exists $\delta_0>0$ such that for any $\delta\in(0,\delta_0)$, the linearized Euler equation $\eqref{eq:uei}$ has a solution $(u_e^{(k)},v_e^{(k)},p_e^{(k)})$ which satisfies
		\begin{align}
			&\|\partial_x^j\partial_y^l(u_e^{(k)},v_e^{(k)})\|_2\leq{C_{j,k,l}}\delta,~\forall j,l\geq0.\nonumber\\
			&\Delta u_e^{(k)}=\Delta v_e^{(k)}=0, (x,y)\in\Omega.\\
			&\int_{0}^{2\pi}v_e^{(k)}(x,y)dx=0,~\forall y\in[0,1].\nonumber
		\end{align}
		Besides, the vorticity $\omega_e^{(k)}$ is a constant.
	\end{proposition}
\begin{proof}
	We only prove the case for $k=2$ and $k=3$. For $k=2$, we only need to prove $\int_{0}^{2\pi}f_e^{(2)}(x,y)dx=0$. In fact, due to $\int_0^{2\pi}v_e^{(1)}(x,y)dx=0,~\forall y\in [0,1]$, we deduce that
\begin{align*}
			\int_{0}^{2\pi}v_e^{(1)}\partial_yu_e^{(1)}(x,y)dx&=\int_{0}^{2\pi}v_e^{(1)}(x,y)\int_{0}^{x}\partial_{xy}u_e^{(1)}(\bar{x},y)d\bar{x}dx\\
			&=-\int_{0}^{2\pi}v_e^{(1)}(x,y)\int_{0}^{x}\partial_{yy}v_e^{(1)}(\bar{x},y)d\bar{x}dx\\
			&=\int_{0}^{2\pi}v_e^{(1)}(x,y)\int_{0}^{x}\partial_{xx}v_e^{(1)}(\bar{x},y)d\bar{x}dx\\
			&=\int_{0}^{2\pi}v_e^{(1)}(x,y)\partial_xv_e^{(1)}(x,y)dx=0.
		\end{align*}
From this and $u_e''(y)=0$, we get $\int_{0}^{2\pi}f_e^{(2)}(x,y)dx=0$.

For $k=3$, We only need to show that $\int_{0}^{2\pi}f_e^{(3)}(x,y)dx=0,~\forall y\in[0,1]$. Since $u_e^{(1)}$ is harmonic, so it becomes
$$\int_{0}^{2\pi}\left(v_e^{(1)}\partial_yu_e^{(2)}+v_e^{(2)}\partial_yu_e^{(1)}\right)dx=0,~\forall y\in[0,1].$$
Denote $\displaystyle V^{(i)}=\int_{0}^{x}v_e^{(i)}(\bar x,y)d\bar x,i=1,2,$ then by the incompressibility and the harmonicity of $(u_e^{(1)},v_e^{(1)}),(u_e^{(2)},v_e^{(2)})$, we have
\begin{align*}
	\int_{0}^{2\pi}\left(v_e^{(1)}\partial_yu_e^{(2)}+v_e^{(2)}\partial_yu_e^{(1)}\right)dx&=-\int_{0}^{2\pi}\left(V^{(1)}\partial_{xy}u_e^{(2)}+V^{(2)}\partial_{xy}u_e^{(1)}\right)dx\\
	&=-\int_{0}^{2\pi}\left(V^{(1)}\partial_{xx}v_e^{(2)}+V^{(2)}\partial_{xx}v_e^{(1)}\right)dx\\
	&=\int_{0}^{2\pi}\left(v_e^{(1)}\partial_xv_e^{(2)}+v_e^{(2)}\partial_xv_e^{(1)}\right)dx=0.
\end{align*}
The remainder of the proof of this proposition follows as in the proof of Proposition $\ref{prop:ue1ve1}$. 
\end{proof}
	
	\begin{remark}
		In this remark, we show that the solvability of (\ref{eq:uei}) for $k=2$ also determines the form of $u_e(y)$. In fact, integrating the first equation of (\ref{eq:uei}) for $k=2$ with respect to $x$ in $[0,2\pi]$, we deduce that
$$u_e''(y)=0,$$
this is the equation in (\ref{eq:Couette flow}).
	\end{remark}
	
\subsubsection{The solvability of equations for $(u_p^{(k)},v_p^{(k+1)},p_p^{(k+1)})$, $(\hat{u}_p^{(k)},\hat{v}_p^{(k+1)},\hat{p}_p^{(k+1)}),k=2,3,4,5,6$}Let
\begin{align}
	\begin{cases}
		u^{\epsilon}(x,y)=u_e(y)+u_p^{(0)}+\sum\limits_{i=1}^{k}\epsilon^i (u_e^{(i)}+u_p^{(i)})+\cdot\cdot\cdot,\\[3pt]
		v^{\epsilon}(x,y)=\sum\limits_{i=1}^{k}\epsilon^i (v_e^{(i)}+v_p^{(i)})+\cdot\cdot\cdot,\\[3pt]
		p^{\epsilon}(x,y)=\sum\limits_{i=1}^{k}\epsilon^i (p_e^{(i)}+p_p^{(i)})+\cdot\cdot\cdot.\label{eq:kthUpperExpansions}
	\end{cases}
\end{align}
Putting $\eqref{eq:kthUpperExpansions}$ into $\eqref{eq:NS}$ and collecting $\epsilon^k$ order terms, we obtain that $\eqref{eq:upi}$.
	Similarly, putting
\begin{align*}
	\begin{cases}
		u^{\epsilon}(x,y)=u_e(y)+\hat{u}_p^{(0)}+\sum\limits_{i=1}^{k}\epsilon^i (u_e^{(i)}+\hat{u}_p^{(i)})+\cdot\cdot\cdot,\\[3pt]
		v^{\epsilon}(x,y)=\sum\limits_{i=1}^{k+1}\epsilon^k (v_e^{(i)}+\hat{v}_p^{(i)})+\cdot\cdot\cdot,\\[3pt]
		p^{\epsilon}(x,y)=\sum\limits_{i=1}^{k+1}\epsilon^k (p_e^{(i)}+\hat{p}_p^{(i)})+\cdot\cdot\cdot
	\end{cases}
\end{align*}
into $\eqref{eq:NS}$ and collecting $\epsilon^k$ order terms, we obtain that $\eqref{eq:^upi}$.
	
	The solvability of $\eqref{eq:upi}$ for $k\geq2$ and $\eqref{eq:^upi}$ for $k\geq2$ are the same with $\eqref{eq:upi}$ for $k=1$. Noticing that Proposition $\ref{prop:uekvek}$ and the fact that $f_{\text{upper}}^{(k)}$ and $g_{\text{upper}}^{(k)}$ decay fast as $\zeta\rightarrow-\infty$, $f_{\text{lower}}^{(k)}$ and $g_{\text{lower}}^{(k)}$ decay fast as $\eta\rightarrow+\infty$, the Prandtl profiles have the following property. 
	\begin{proposition}\label{prop:upkvpk+1}
		There exists $\delta_0>0$ such that for any $\delta\in(0,\delta_0)$, the equation $\eqref{eq:upi}$ has a solution $(u_p^{(k)},v_p^{(k+1)})$ and equations $\eqref{eq:^upi}$ has a solution $(\hat{u}_p^{(k)},\hat{v}_p^{(k+1)})$ satisfying
		\begin{equation}
			\begin{aligned}
				&\int_{-\infty}^{0}\int_{0}^{2\pi}|\partial_x^j\partial_{\zeta}^s(u_p^{(k)}-A_k,v_p^{(k+1)})|^2\langle\zeta\rangle^{2l}dxd\zeta\leq{C_{j,k,l,s}}\delta^2,\\
				&\int_{0}^{+\infty}\int_{0}^{2\pi}|\partial_x^j\partial_{\eta}^s(\hat{u}_p^{(k)}-\hat{A}_k,\hat{v}_p^{(k+1)})|^2\langle\eta\rangle^{2l}dxd\eta\leq{C_{j,k,l,s}}\delta^2,\\
				&\int_{0}^{2\pi}v_p^{(k+1)}(x,\zeta)dx=\int_{0}^{2\pi}\hat{v}_p^{(k+1)}(x,\eta)dx=0,\forall\zeta\leq0,\forall\eta\geq0,\nonumber
			\end{aligned}
		\end{equation}
		where $A_k=\lim\limits_{\zeta\rightarrow-\infty}u_p^{(k)},~\hat{A}_k=\lim\limits_{\eta\rightarrow+\infty}\hat{u}_p^{(k)}$ are constants satisfying $|A_k|+|\hat{A}_k|\leq{C}_k\delta$.
	\end{proposition}
	As before, we make modifications. Let
	\begin{equation}
		\begin{aligned}
			\tilde{u}_e^{(k)}&=u_e^{(k)}+yA_k+(1-y)\hat{A}_k,~\tilde{v}_e^{(k)}=v_e^{(k)},~\tilde{p}_e^{(k)}=p_e^{(k)},\\[3pt]
			\tilde{u}_p^{(k)}&=u_p^{(k)}-A_k,~\tilde{v}_p^{(k+1)}=v_p^{(k+1)},~\tilde{p}_p^{(k+1)}=p_p^{(k+1)},\\[3pt]
			\tilde{\hat{u}}_p^{(k)}&=\hat{u}_p^{(k)}-\hat{A}_k,~\tilde{\hat{v}}_p^{(k+1)}=\hat{v}_p^{(k+1)},~\tilde{\hat{p}}_p^{(k+1)}=\hat{p}_p^{(k+1)},\label{eq:uekmodification}
		\end{aligned}
	\end{equation}
	Then  $(\tilde{u}_p^{(k)},\tilde{v}_p^{(k+1)},\tilde{p}_p^{(k+1)})$ satisfies 
	\begin{align}
		\begin{cases}			(u_e(1)+u_p^{(0)})\partial_x\tilde{u}_p^{(k)}+\tilde{u}_p^{(k)}\partial_xu_p^{(0)}+(\tilde{v}_e^{(1)}(x,1)+v_p^{(1)})\partial_{\zeta}\tilde{u}_p^{(k)}\\[5pt]
			+(\tilde{v}_p^{(k+1)}
			-\tilde{v}_p^{(k+1)}(x,0))\partial_{\zeta}u_p^{(0)}-\partial_{\zeta}^2\tilde{u}_p^{(k)}=\tilde{f}_{\text{upper}}^{(k)},\\[5pt]
			\partial_x\tilde{u}_p^{(k)}+\partial_{\zeta}\tilde{v}_p^{(k+1)}=0,\\[5pt]			\tilde{u}_p^{(k)}|_{\zeta=0}=-\tilde{u}_e^{(k)}|_{y=1},\lim\limits_{\zeta\rightarrow-\infty}\partial_{\zeta}\tilde{u}_p^{(k)}
			=\lim\limits_{\zeta\rightarrow-\infty}\tilde{v}_p^{(k+1)}=0\label{eq:upkmodification}
		\end{cases}
	\end{align}
	and $$\partial_{\zeta}\tilde{p}_p^{(k+1)}=\tilde{g}_{\text{upper}}^{(k)},$$
	where $\tilde{f}_{\text{upper}}^{(k)}=f_{\text{upper}}^{(k)}+A_k\partial_xu_p^{(0)}$ and $\tilde{g}_{\text{upper}}^{(k)}=g_{\text{upper}}^{(k)}$. Similarly, $(\tilde{\hat{u}}_p^{(k)},\tilde{\hat{v}}_p^{(k+1)})$ satisfies equation of the same form and so we omit. Moreover, $(\tilde{u}_e^{(k)},\tilde{v}_e^{(k)},\tilde{p}_e^{(k)})$ satisfies
	\begin{align}
		\begin{cases}
			u_e(y)\partial_x\tilde{u}_e^{(k)}+\tilde{v}_e^{(k)}u_e'(y)+\partial_x\tilde{p}_e^{(k)}=\tilde{f}_e^{(k)},\\[3pt]
			u_e(y)\partial_x\tilde{v}_e^{(k)}+\partial_y\tilde{p}_e^{(k)}=\tilde{g}_e^{(k)},\\[3pt]
			\partial_x\tilde{u}_e^{(k)}+\partial_y\tilde{v}_e^{(k)}=0,\\[3pt]
			\tilde{u}_e^{(k)}(x+2\pi,y)=\tilde{u}_e^{(k)}(x,y),~\tilde{v}_e^{(k)}(x+2\pi,y)=\tilde{v}_e^{(k)}(x,y),\\[3pt]
			\tilde{v}_e^{(k)}|_{y=1}=-\tilde{v}_p^{(k)}|_{\zeta=0},~\tilde{v}_e^{(k)}|_{y=0}=-\tilde{\hat{v}}_p^{(k)}|_{\eta=0}.\label{eq:uekmodeq}
		\end{cases}
	\end{align} 
	where $\tilde{f}_e^{(k)}=f_e^{(k)},\tilde{g}_e^{(k)}=g_e^{(k)}$	and
	\begin{align*}
		&\|\partial_x^j\partial_y^l(\tilde{u}_e^{(k)},\tilde{v}_e^{(k)})\|_2\leq{C_{j,k,l}}\delta,~\forall j,l\geq0.\\
		&\Delta \tilde{u}_e^{(k)}=\Delta \tilde{v}_e^{(k)}=0, (x,y)\in\Omega.\\
		&\int_{0}^{2\pi}\tilde{v}_e^{(k)}(x,y)dx=0,~\forall y\in[0,1]. 
	\end{align*}
	and has a constant vorticity.
	Besides $\tilde{u}_p^{(k)},\tilde{\hat{u}}_p^{(k)}$ have zero limit: 
	$$\lim\limits_{\zeta\to-\infty}\tilde{u}_p^{(k)}=\lim\limits_{\eta\to+\infty}\tilde{\hat{u}}_p^{(k)}=0.$$
	For simplicity, we omit the tilde superscripts. We continue this process until $(u_e^{(6)},v_e^{(6)})$ $(u_p^{(6)},v_p^{(7)})$ and $(\hat{u}_p^{(6)},\hat{v}_p^{(7)})$ are solved out.

	\subsubsection{The solvability of equation for $(u_e^{(7)},v_e^{(7)},p_e^{(7)})$}
	
	Putting
	\begin{align}
		\begin{cases}
			u^{\epsilon}(x,y)=u_e(y)+\sum\limits_{i=1}^{7}\epsilon^i u_e^{(i)}(x,y)+\cdot\cdot\cdot,\\
			v^{\epsilon}(x,y)=\sum\limits_{i=1}^{7}\epsilon^i v_e^{(i)}(x,y)+\cdot\cdot\cdot,\\
			p^{\epsilon}(x,y)=\sum\limits_{i=1}^{7}\epsilon^i p_e^{(i)}(x,y)+\cdot\cdot\cdot.\label{eq:endEulerExpansions}
		\end{cases}
	\end{align}
	into $\eqref{eq:NS}$ and collecting $\epsilon^7$ order terms, we obtain $\eqref{eq:uei}$ for $k=7$.
	By the similar argument of Proposition $\ref{prop:uekvek}$, we have
	\begin{proposition}\label{prop:ue end}
		There exists $\delta_0>0$ such that for any $\delta\in(0,\delta_0)$, the linearized Euler equation $\eqref{eq:uei}$ for $k=7$ has a solution $(u_e^{(7)},v_e^{(7)})$ which satisfies
		\begin{align}
			&\|\partial_x^j\partial_y^l(u_e^{(7)},v_e^{(7)})\|_2\leq{C_{j,l}}\delta,\nonumber\\
			&\Delta u_e^{(7)}=\Delta v_e^{(7)}=0, (x,y)\in\Omega.\\
			&\int_{0}^{2\pi}v_e^{(7)}(x,y)dx=0,~\forall y\in[0,1].\nonumber
		\end{align}
		Besides, the vorticity $\omega_e^{(7)}$ is a constant.
	\end{proposition}
	
	\subsubsection{The solvability of equations for $(u_p^{(7)},v_p^{(8)},p_p^{(8)})$ and $(\hat{u}_p^{(7)},\hat{v}_p^{(8)},\hat{p}_p^{(8)})$}
	
	Putting
	\begin{align*}
		\begin{cases}
			u^{\epsilon}(x,y)=u_e(y)+u_p^{(0)}+\sum\limits_{i=1}^{7}\epsilon^i (u_e^{(i)}+u_p^{(i)})+\cdot\cdot\cdot,\\[3pt]
			v^{\epsilon}(x,y)=\sum\limits_{i=1}^{8}\epsilon^i (v_e^{(i)}+v_p^{(i)})+\cdot\cdot\cdot,\\[3pt]
			p^{\epsilon}(x,y)=\sum\limits_{i=1}^{8}\epsilon^i (p_e^{(i)}+p_p^{(i)})+\cdot\cdot\cdot
		\end{cases}
	\end{align*}
	into $\eqref{eq:NS}$ and collecting $\epsilon^7$ order terms, we obtain $\eqref{eq:up end}$. Besides, putting
	\begin{align*}
		\begin{cases}
			u^{\epsilon}(x,y)=u_e(y)+\hat{u}_p^{(0)}+\sum\limits_{i=1}^{7}\epsilon^i (u_e^{(i)}+\hat{u}_p^{(i)})+\cdot\cdot\cdot,\\[3pt]
			v^{\epsilon}(x,y)=\sum\limits_{i=1}^{8}\epsilon^i (v_e^{(i)}+\hat{v}_p^{(i)})+\cdot\cdot\cdot,\\[3pt]
			p^{\epsilon}(x,y)=\sum\limits_{i=1}^{8}\epsilon^i (p_e^{(i)}+\hat{p}_p^{(i)})+\cdot\cdot\cdot
		\end{cases}
	\end{align*}
	into $\eqref{eq:NS}$ and collecting $\epsilon^7$ order terms, we obtain $\eqref{eq:^up end}$.
	We have the following conclusion similar to Proposition $\ref{prop:up1vp2Estimation}.$
	
	\begin{proposition}\label{prop:up end}
		There exists $\delta_0>0$ such that for any $\delta\in(0,\delta_0)$, the equation $\eqref{eq:up end}$ has a solution $(u_p^{(7)},v_p^{(8)})$ and the equation $\eqref{eq:^up end}$ has a solution $(\hat{u}_p^{(7)},\hat{v}_p^{(8)})$ satisfying
		\begin{align*}
			&\int_{-\infty}^{0}\int_{0}^{2\pi}|\partial_x^j\partial_{\zeta}^s(u_p^{(7)},v_p^{(8)})|^2\langle\zeta\rangle^{2l}dxd\zeta\leq{C_{j,l,s}}\delta^2,\\
			&\int_{0}^{+\infty}\int_{0}^{2\pi}|\partial_x^j\partial_{\eta}^s(\hat{u}_p^{(7)},\hat{v}_p^{(8)})|^2\langle\eta\rangle^{2l}dxd\eta\leq{C_{j,l,s}}\delta^2,\\
			&\int_{0}^{2\pi}v_p^{(8)}(x,\zeta)dx=\int_{0}^{2\pi}\hat{v}_p^{(8)}(x,\eta)dx=0,~\forall\zeta\leq0,\forall \eta\geq0,
		\end{align*}
	\end{proposition}

\subsection{Approximate solution}\label{subsection:Approximate solution}

	 In this subsection, we construct an approximate solution of Navier-Stokes equations (\ref{eq:NS}). 
	 First we choose a function $\chi(y)\in{C^{\infty}([0,1])}$ satisfying
	\begin{align}
		\chi(y)=\begin{cases}
			1,~y\in[3/4,1],\\[3pt]
			0,~y\in[0,1/4].\label{eq:chi}
		\end{cases}
	\end{align}
	Set
	\begin{align}
\bar{u}_p^a(x,y)&:=\chi(y)\sum_{i=0}^{7}\epsilon^iu_p^{(i)}+(1-\chi(y))\sum_{i=0}^{7}\epsilon^i\hat{u}_p^{(i)}
:=\chi(y)u_p^a+(1-\chi(y))\hat{u}_p^a,\nonumber\\
	\bar{v}_p^a(x,y)&:=\chi(y)\sum_{i=1}^{8}\epsilon^iv_p^{(i)}+(1-\chi(y))\sum_{i=1}^{8}\epsilon^i\hat{v}_p^{(i)}:=\chi(y)v_p^a+(1-\chi(y))\hat{v}_p^a,\nonumber\\
	\bar{p}_p^a(x,y)&:=\chi^2(y)\sum_{i=1}^{8}\epsilon^ip_p^{(i)}+(1-\chi(y))^2\sum_{i=1}^{8}\epsilon^i\hat{p}_p^{(i)}:=\chi^2(y)p_p^a+(1-\chi(y))^2\hat{p}_p^a\nonumber
\end{align}
	and
\begin{align}
u_e^a:=u_e(y)+\sum_{i=1}^{7}\epsilon^iu_e^{(i)},~
	v_e^a:=\sum_{i=1}^{7}\epsilon^iv_e^{(i)},~
	p_e^a:=\sum_{i=1}^{7}\epsilon^ip_e^{(i)}.\nonumber
\end{align}
	We construct an approximate solution $(u^a,v^a,p^a)$ by
	\begin{align}
u^a(x,y):=u_e^a+\bar{u}_p^a+h(x,y),~
	v^a(x,y):=v_e^a+\bar{v}_p^a,~
	p^a(x,y):=p_e^a+\bar{p}_p^a,\nonumber
\end{align}
	where the corrector $h(x,y)$ will be given in Appendix,which satisfies
	$$h(x,0)=h(x,1)=0,~~\|\partial_x^j\partial_y^kh\|_2\leq{C_{j,k,m}}\epsilon^m\delta\rho_m$$
	and makes $(u^a,v^a)$ be divergence-free
	$$\partial_xu^a+\partial_yv^a=0,$$
	where $\rho_m=(1+|\zeta|)^{-m}+(1+\eta)^{-m}$ with any positive integer $m$.
	Moreover, $(u^a,v^a)$ satisfies the following boundary conditions
	\begin{align}
u^a(x+2\pi,y)&=u^a(x,y),~v^a(x+2\pi,y)=v^a(x,y),\nonumber\\[3pt]
	u^a(x,1)&=\alpha+\delta f(x),~v^a(x,1)=0,\nonumber\\[3pt]
	u^a(x,0)&=-\beta+\delta g(x),~v^a(x,0)=0.\nonumber
\end{align}
	By collecting the estimates in Proposition $\ref{prop:ue1ve1}$, Proposition $\ref{prop:uekvek}$, Proposition $\ref{prop:ue end}$ and using the Sobolev embedding, we deduce that
	$$\|\partial_x^j\partial_y^k(u_e^a-u_e(y))\|_{\infty}\leq C_{j,k}\epsilon\delta,~~\|\partial_x^j\partial_y^kv_e^a\|_{\infty}\leq C_{j,k}\epsilon\delta.$$ 
	and
	$$\Delta u_e^a=\Delta v_e^a=0,$$
	which implies that the outer vorticity $\omega_e^a:=\partial_yu_e^a-\partial_xv_e^a$ is a constant which satisfies $$|\omega_e^a-a|\leq C\epsilon\delta.$$
Similarly, collecting the estimates in Proposition $\ref{prop:up0vp1}$, $\ref{prop:^up0vp1}$, $\ref{prop:up1vp2Estimation}$, $\ref{prop:^up1vp2Estimation}$, $\ref{prop:upkvpk+1}$ and Proposition $\ref{prop:up end}$, we have
	\begin{align}
&\|\zeta^l\partial_x^j\partial_{\zeta}^ku_p^a\|_{\infty}\leq{C_{j,k,l}}\delta,
~~\|\zeta^l\partial_x^j\partial_{\zeta}^{k}v_p^a\|_{\infty}\leq{C_{j,k,l}}\epsilon\delta,\nonumber\\[3pt] &\|\eta^l\partial_x^j\partial_{\eta}^k\hat{u}_p^a\|_{\infty}\leq{C_{j,k,l}}\delta,~~
\|\eta^l\partial_x^j\partial_{\eta}^{k}\hat{v}_p^a\|_{\infty}\leq{C_{j,k,l}}\epsilon\delta.\nonumber
\end{align}
	
	Finally, we set 
	$$R_u^a:=u^au_x^a+v^au_y^a+p_x^a-\epsilon^2\triangle u^a,$$
	$$R_v^a:=u^av_x^a+v^av_y^a+p_y^a-\epsilon^2\triangle v^a,$$
	then there holds
	$$\|\partial_x^j\partial_y^k(R_u^a,R_v^a)\|_2\leq C_{j,k}\epsilon^{8-k},\forall j,k\geq0$$
	and $(u^a,v^a,p^a)$ satisfies
	\begin{align}
		\begin{cases}
			u^au_x^a+v^au_y^a+p_x^a-\epsilon^2\triangle u^a=R_u^a \text{ in } \Omega,\\[3pt]
			u^av_x^a+v^av_y^a+p_y^a-\epsilon^2\triangle v^a=R_v^a\text{ in }\Omega,\\[3pt]
			u_x^a+v_y^a=0\text{ in }\Omega,\\[3pt]
			u^a(x+2\pi,y)=u^a(x,y),~v^a(x+2\pi,y)=v^a(x,y)\text{ in }\Omega,\\[3pt]
			u^a(x,1)=\alpha+\delta f(x),~v^a(x,1)=0,~x\in[0,2\pi],\\[3pt]
			u^a(x,0)=-\beta+\delta g(x),~v^a(x,0)=0,~x\in[0,2\pi].
		\end{cases}\label{eq:approximate equation}
	\end{align}
	
	We conclude the above as follows.
	\begin{proposition}\label{lem:approxEstimates}
		The approximate solution $(u^a,v^a)$ satisfying $\eqref{eq:approximate equation}$ is decomposed as
		\begin{align}
			u^a:=u_e^a+\chi(y)u_p^a+(1-\chi(y))\hat{u}_p^a+h(x,y),~v^a:=v_e^a+\chi(y)v_p^a+(1-\chi(y))\hat{v}_p^a, \nonumber
		\end{align} 
		The Euler part $(u_e^a,v_e^a)$ is harmonic with a constant vorticity $\omega_e^a=a+O(\epsilon\delta)$ and satisfies
		\begin{align}
			\|\partial_x^j\partial_y^k(u_e^a-u_e(y))\|_{\infty}\leq C_{j,k}\epsilon\delta,~\|\partial_x^j\partial_y^kv_e^a\|_{\infty}\leq C_{j,k}\epsilon\delta.\nonumber
		\end{align}
		The upper boundary layer part $(u_p^a,v_p^a)$ satisfies 
		\begin{align}
			\|\zeta^l\partial_x^j\partial_{\zeta}^ku_p^a\|_{\infty}\leq{C_{j,k,l}}\delta,~
\|\zeta^l\partial_x^j\partial_{\zeta}^{k}v_p^a\|_{\infty}\leq C_{j,k,l}\epsilon\delta.\nonumber
		\end{align}
		The lower boundary layer part $(\hat{u}_p^a,\hat{v}_p^a)$ satisfies 
		\begin{align}
			\|\eta^l\partial_x^j\partial_{\eta}^k\hat{u}_p^a\|_{\infty}\leq{C_{j,k,l}}\delta,~
\|\eta^l\partial_x^j\partial_{\eta}^{k}\hat{v}_p^a\|_{\infty}\leq{C_{j,k,l}}\epsilon\delta.\nonumber
		\end{align}
		While $h(x,y)$ satisfies
		$$h(x,0)=h(x,1)=0,~~|\partial_x^j\partial_y^kh|\leq C_{j,k,m}\epsilon^m\delta\rho_m,$$
		where $\rho_m=(1+|\zeta|)^{-m}+(1+\eta)^{-m}$ with any positive integer $m$.
	\end{proposition}
		
	\section{Linear stability estimates of the error equations}

\subsection{Error equation}
	Define the error functions by
	$$u:=u^{\epsilon}-u^a,~v:=v^{\epsilon}-v^a,~p:=p^{\epsilon}-p^a.$$
	Then $(u,v,p)$ satisfies 
	\begin{align}
		\begin{cases}
			-\epsilon^2\Delta u+p_x+S_u=R_u\text{ in }\Omega,\\[3pt]
			-\epsilon^2\Delta v+p_y+S_v=R_v\text{ in }\Omega,\\[3pt]
			u_x+v_y=0\text{ in }\Omega,\\[3pt]
			u(x+2\pi,y)=u(x,y),~v(x+2\pi,y)=v(x,y)\text{ in }\Omega,\\[3pt]
			u|_{y=0,1}=v|_{y=0,1}=0,x\in[0,2\pi],\label{eq:errEqu}
		\end{cases}
	\end{align}
	where
\begin{equation}
		\begin{aligned}
			S_u:&=u^au_x+v^au_y+uu_x^a+vu_y^a,\nonumber\\[3pt]
			S_v:&=u^av_x+v^av_y+uv_x^a+vv_y^a\nonumber
		\end{aligned}
	\end{equation}
and
	\begin{equation}
		\begin{aligned}
			R_u:=-R_u^a-uu_x-vu_y,~
			R_v:=-R_v^a-uv_x-vv_y.\nonumber
		\end{aligned}
	\end{equation}

\subsection{Linear stability estimate}
	To derive the error estimates, we first give the following Lemma which can be deduced from Proposition $\ref{lem:approxEstimates}$.
	\begin{lemma}\label{lem:uxvxEstimate}
		There exist $\epsilon_0>0, \delta_0>0$ such that for any $\epsilon\in (0,\epsilon_0), \delta\in (0,\delta_0)$, there holds
		\begin{align*}
			& \|\partial_x^j(u^a-u_e(y))\|_\infty\leq{C_{j}\delta},
			~~\|\partial_x^jv^a\|_\infty\leq C_j\epsilon\delta,~\forall j\in \mathbb{N},\\
			&|u^a-u_e^a|\leq C_m\delta\rho_m, |v^a-v_e^a|\leq C_m\epsilon\delta\rho_m,|\omega^a-\omega_e^a|\leq C_m\epsilon^{-1}\delta\rho_m.
		\end{align*}
	\end{lemma}
\begin{proof}
	From the bound estimates in Proposition $\ref{lem:approxEstimates}$, we have
	\begin{align*}
|\partial_x^j(u^a-u_e(y))|&\leq|\partial_x^j(u_e^a-u_e(y))|+\chi(y)|\partial_x^ju_p^a|+(1-\chi(y))|\partial_x^j\hat{u}_p^a|+|\partial_x^jh|\leq C_j\delta,\nonumber\\[3pt]
|\partial_x^jv^a|&\leq|\partial_x^jv_e^a|+\chi(y)|\partial_x^jv_p^a|+(1-\chi(y))|\partial_x^j\hat{v}_p^a|\leq C_j\epsilon\delta.
\end{align*}
It follows from fast decay of Prandtl profiles in Proposition $\ref{lem:approxEstimates}$ that
\begin{align*}
	|u^a-u_e^a|&\leq|u_p^a|+|\hat{u}_p^a|+|h|\leq C_m\delta\rho_m,\\
	|v^a-v_e^a|&\leq|v_p^a|+|\hat{v}_p^a|\leq C_m\epsilon\delta\rho_m,\\
	|\omega^a-\omega_e^a|&\leq|\partial_y(u^a-u_e^a)-\partial_x(v^a-v_e^a)|\\
	&\leq|\partial_y(\chi(y)u_p^a)|+|\partial_y[(1-\chi(y))\hat{u}_p^a]|+|h_y|+|\partial_xv_p^a|+|\partial_x\hat{v}_p^a|\\
	&\leq\epsilon^{-1}(|\partial_{\zeta}u_p^a|+|\partial_{\eta}\hat{u}_p^a|)+C(|u_p^a|+|\hat{u}_p^a|)+|\partial_xv_p^a|+|\partial_x\hat{v}_p^a|+|h_y|\\
	&\leq C_m\epsilon^{-1}\delta\rho_m.
\end{align*}
	\end{proof}
	
The following weighted Hardy inequality will be used. Let $c_0>0$ is small but independent of $\epsilon$ such that the zero set of $u^a$ is contained in $(y^*-\frac{c_0}{2},y^*+\frac{c_0}{2})$.
	\begin{lemma}[Hardy inequality]\label{lem:Hardy}
		Suppose that $\sigma(y)\in C^{\infty}([0,1])$ is a signed monotone function satisfying
		\begin{equation*}
			\sigma(y)=
			\begin{cases}
				1,y\in[y^*+c_0,1],\\
				0,y\in[y^*-\frac{c_0}{2},y^*+\frac{c_0}{2}],\\
				-1,y\in[0,y^*-c_0].
			\end{cases}
		\end{equation*}
		Then for $f\in C^1([0,1])$ with $f(0)=f(1)=0$ and $\alpha>0$, we have
		\begin{align}\label{weight Hardy 1}
\int_{0}^{1}|\sigma|\frac{f^2}{y^2}dy+\int_{0}^{1}|\sigma|\frac{f^2}{(1-y)^2}dy\leq C\int_{0}^{1}|\sigma|f_y^2dy+C\int_{T_{\sigma}}f^2dy,
\end{align}
and
\begin{align}\label{weight Hardy 2}
\int_{0}^{1}|\sigma|f^2y^{\kappa-1}dy\leq C_{\kappa}\int_{0}^{1}|\sigma|f_y^2y^{\kappa+1}dy+C_{\kappa}\int_{T_{\sigma}} f^2dy,
\end{align}
		where $T_{\sigma}$ is the support of $\sigma'.$
	\end{lemma}
	\begin{proof}
		Integrating by parts and noticing that $\sigma(y^*)=0$, we have
		\begin{align*}
			\int_{0}^{1}|\sigma|\frac{f^2}{y^2}dy=&\int_{0}^{y^*}\frac{-\sigma f^2}{y^2}dy+\int_{y^*}^{1}\frac{\sigma f^2}{y^2}dy\\
			=&-\int_{0}^{y^*}\frac{\sigma'f^2+2\sigma ff_y}{y}dy+\int_{y^*}^{1}\frac{\sigma'f^2+2\sigma ff_y}{y}dy\\
			\leq&C\int_{T_{\sigma}}f^2dy+C\left(\int_{0}^{1}|\sigma|\frac{f^2}{y^2}dy\right)^{\frac{1}{2}}\left(\int_{0}^{1}|\sigma|f_y^2dy\right)^{\frac{1}{2}}.
		\end{align*}
		Then we obtain \eqref{weight Hardy 1} after applying Young inequality. The proof of (\ref{weight Hardy 2}) is similar. 
	\end{proof}

Now we state the linear stability estimate of error equation $\eqref{eq:errEqu}$. Let $L_{\sigma}^2(\Omega)$ be the space of periodic,distributionally divergence-free vector fields with zero normal trace, and let $\mathbb{P}$ be its orthogonal Leray projection.
\begin{Proposition}\label{prop:estimate}
		Suppose that $F_1,F_2\in L^2(\Omega)$, then there exists $\delta_0>0$ such that for any fixed $\delta\in(0,\delta_0)$, there are two constants $\epsilon_0>0$ and $C>0$ such that for any $\epsilon\in(0,\epsilon_0)$, the system
		\begin{align}
			\begin{cases}
				-\epsilon^2\Delta u+p_x+S_u=F_1\text{ in }\Omega,\\[3pt]
				-\epsilon^2\Delta v+p_y+S_v=F_2\text{ in }\Omega,\label{eq:linearEstimate}
			\end{cases}
		\end{align}
		has a unique solution $(u,v)$ in $D:=L_{\sigma}^2(\Omega)\cap[H_0^1(\Omega)]^2\cap[H^2(\Omega)]^2$, which satisfies
		\begin{align}
			\|(u,v)\|_2+\epsilon\|\nabla(u,v)\|_2+\epsilon^2\|\nabla^2(u,v)\|_2\leq&C\epsilon^{-2}\|(F_1,F_2)\|_2.\label{ineq:linear stability}
		\end{align}
		Besides, if $F_1,F_2\in H^1(\Omega)$, then $(u,v)\in[H^3(\Omega)]^2$ and satisfies
		\begin{align}
			\epsilon^2\|\nabla^3(u,v)\|_2\leq&C\|\nabla(F_1,F_2)\|_2+C\epsilon^{-4}\|(F_1,F_2)\|_2.\label{ineq:H3}
		\end{align}
	\end{Proposition}

This proposition can be obtained by combining the Proposition $\ref{prop:wall estimate}$, $\ref{prop:interior estimate}$, $\ref{prop:uv estimate}$, Corollary $\ref{cor:v gradient estimate}$, and Stokes estimates together. 

\subsubsection{Wall positivity estimate}	We first give the following positivity estimate near the wall.
	\begin{lemma}\label{lem:positivity}
		Suppose that $(u,v)$ is the smooth solution of $\eqref{eq:linearEstimate}$, then there exists $\epsilon_0>0,\delta_0>0$ and a constant $C>0$ such that for any $ \epsilon\in(0,\epsilon_0),~ \delta\in(0,\delta_0)$, the following estimate holds
		\begin{align*}
			\int_{0}^{1}\int_{0}^{2\pi}|\sigma|(u_x^2+v_x^2)dxdy\leq&C\int_{0}^{1}\int_{0}^{2\pi}[\epsilon^2(u_x^2+u_y^2+v_x^2)+v^2]dxdy+C\int_{0}^{1}\int_{0}^{2\pi}(F_1^2+F_2^2)dxdy.
		\end{align*}
	\end{lemma}
\begin{proof} Multiplying the first equation of $\eqref{eq:linearEstimate}$ with $-(\sigma v)_y$ and the second equation with $(\sigma v)_x$, adding them together and integrating in $\Omega$, we obtain 
	\begin{align}
		-\epsilon^2\int_{0}^{1}\int_{0}^{2\pi}&\Big[(\sigma u_x-\sigma'v)\Delta u+\sigma v_x\Delta v\Big]dxdy+\int_{0}^{1}\int_{0}^{2\pi}[S_u(\sigma u_x-\sigma'v)+S_v\sigma v_x]dxdy\notag\\
		=&\int_{0}^{1}\int_{0}^{2\pi}[F_1(\sigma u_x-\sigma'v)+F_2\sigma v_x]dxdy,\label{eq:poseq1}
	\end{align}
	where we eliminate pressure terms by integration by parts. 
Next we deal with $\eqref{eq:poseq1}$ term by term. \\
{\bf The diffusion term:} By divergence-free and periodic condition and integration by parts, we obtain
	\begin{align}
		\epsilon^2&\left|\int_{0}^{1}\int_{0}^{2\pi}\left[(\sigma u_x-\sigma'v)\Delta u+\sigma v_x\Delta v\right]dxdy\right|\nonumber\\
=&\epsilon^2\left|\int_{0}^{1}\int_{0}^{2\pi}\left[u_{yy}\sigma u_x-\sigma'(u_{xx}+u_{yy})v+v_{yy}\sigma v_x\right]dxdy\right|\notag\\
=&\epsilon^2\left|\int_{0}^{1}\int_{0}^{2\pi}\left[u_{y}(\sigma'u_x+\sigma u_{xy})-\sigma'u_xv_x-u_y(\sigma''v+\sigma'v_y)+v_y(\sigma'v_x+\sigma v_{xy})\right]dxdy\right|\notag\\
=&\epsilon^2\left|\int_{0}^{1}\int_{0}^{2\pi}\left[2\sigma'u_xu_y-2\sigma'u_xv_x-\sigma''u_yv\right]dxdy\right|\notag\\
\leq&C\epsilon^2\int_{0}^{1}\int_{0}^{2\pi}(u_x^2+u_y^2+v_x^2)dxdy.\label{eq:poseq2}
	\end{align}
{\bf The convective term}: By direct computation, we get that
	\begin{align}
		\int_{0}^{1}\int_{0}^{2\pi}[S_u(\sigma u_x-\sigma'v)+S_v\sigma v_x]dxdy&=\sum_{k=1}^{5}I_k,\label{est:positivity5}
	\end{align}
	where
	\begin{align}
		I_1&=\int_{0}^{1}\int_{0}^{2\pi}\sigma u^a(u_x^2+v_x^2)dxdy,\nonumber\\
		I_2&=\int_{0}^{1}\int_{0}^{2\pi}vu_y^a(\sigma u_x-\sigma'v)dxdy,\nonumber\\
		I_3&=\int_{0}^{1}\int_{0}^{2\pi}v^au_y(\sigma u_x-\sigma'v)dxdy,\nonumber\\
		I_4&=\int_{0}^{1}\int_{0}^{2\pi}uu_x^a(\sigma u_x-\sigma'v)dxdy,\nonumber\\
		I_5&=\int_{0}^{1}\int_{0}^{2\pi}[(uv_x^a+vv_y^a)\sigma v_x+v^av_y\sigma v_x]dxdy.\nonumber
	\end{align}
Estimate of $I_1$:	By the fact that $\sigma u^a\geq c|\sigma|$, we deduce that 
	\begin{align}
		I_1\geq&c\int_{0}^{1}\int_{0}^{2\pi}|\sigma|(u_x^2+v_x^2)dxdy.\label{est:positivityI_1}
	\end{align}
Estimate of $I_2$:	
	By integrating by parts and divergence-free condition, we obtain that
	\begin{align*}
		I_2=&\frac{1}{2}\int_{0}^{1}\int_{0}^{2\pi}v^2(u_{yy}^a\sigma-u_y^a\sigma')dxdy\notag\\
		=&\frac{1}{2}\int_{0}^{1}\int_{0}^{2\pi}[(u^a-u_e(y))_{yy}\sigma v^2-(u^a-u_e(y))_y\sigma'v^2-u_e'(y)\sigma v^2]dxdy\notag\\
		=&\frac{1}{2}\int_{0}^{1}\int_{0}^{2\pi}\{[(u_e^a-u_e(y))_{yy}+(\chi u_p^a)_{yy}+((1-\chi)\hat{u}_p^a)_{yy}+h_{yy}]\sigma v^2\notag\\
		&-[(u_e^a-u_e(y))_y+(\chi u_p^a)_y+((1-\chi)\hat{u}_p^a)_y+h_y]\sigma'v^2-u_e'(y)\sigma'v^2\}dxdy.
	\end{align*}
	Then from Proposition $\ref{lem:approxEstimates}$, Lemma $\ref{lem:Hardy}$ and Poincar\'e inequality, we have
	\begin{align}
		|I_2|\leq&C\int_{0}^{1}\int_{0}^{2\pi}\left(\epsilon\delta|\sigma| v^2+|\zeta^2\partial_{\zeta}^2u_p^a|\frac{|\sigma|v^2}{(1-y)^2}+|\eta^2\partial_{\eta}^2\hat{u}_p^a|\frac{|\sigma|v^2}{y^2}\right)dxdy\notag\\
		&+C\int_{0}^{1}\int_{0}^{2\pi}\left(v^2+|\zeta\partial_{\zeta}u_p^a|\frac{|\sigma'|v^2}{1-y}+|\eta\partial_{\eta}\hat{u}_p^a|\frac{|\sigma'|v^2}{y}\right)dxdy\notag\\
		\leq&C\delta\int_{0}^{1}\int_{0}^{2\pi}|\sigma|(v_x^2+v_y^2)dxdy+C\int_{0}^{1}\int_{0}^{2\pi}v^2dxdy.\label{est:positivityI_2}
	\end{align}
	Estimate of $I_3$: By Lemma $\ref{lem:uxvxEstimate}$,we obtain that
	\begin{align}
		|I_3|\leq&C\epsilon\delta\int_{0}^{1}\int_{0}^{2\pi}(|\sigma||u_xu_y|+|\sigma'||u_yv|)dxdy\notag\\
		\leq&C\delta\int_{0}^{1}\int_{0}^{2\pi}|\sigma|u_x^2dxdy+C\epsilon^2\delta\int_{0}^{1}\int_{0}^{2\pi}u_y^2dxdy+C\delta\int_{0}^{1}\int_{0}^{2\pi}v^2dxdy.\label{est:positivityI_3}
	\end{align}
	Estimate of $I_4$: Similarly,from Lemma $\ref{lem:uxvxEstimate}$ and Hardy inequality,there holds
	\begin{align}
		|I_4|\leq&\int_{0}^{1}\int_{0}^{2\pi}(|(u_e^a-u_e(y))_x+\chi u_{px}^a+(1-\chi)\hat{u}_{px}^a+h_x||\sigma u_xu|+C\epsilon\delta|\sigma'uv|)dxdy\notag\\
		\leq&C\int_{0}^{1}\int_{0}^{2\pi}\left(\epsilon\delta|\sigma u_xu|+\epsilon|\zeta u_p^a|\frac{|\sigma u_xu|}{1-y}+\epsilon|\eta\hat{u}_p^a|\frac{|\sigma u_xu|}{y}\right)dxdy+C\epsilon\delta\int_{0}^{1}\int_{0}^{2\pi}|uv|dxdy\notag\\
		\leq&C\epsilon^2\delta\int_{0}^{1}\int_{0}^{2\pi}u_y^2dxdy+C\delta\int_{0}^{1}\int_{0}^{2\pi}|\sigma|u_x^2dxdy+C\delta\int_{0}^{1}\int_{0}^{2\pi}v^2dxdy.\label{est:positivityI_4}
	\end{align}
	Estimate of $I_5$: It follows from Lemma $\ref{lem:uxvxEstimate}$ and Hardy inequality that
	\begin{align}
		|I_5|\leq&C\int_{0}^{1}\int_{0}^{2\pi}(\epsilon\delta|\sigma uv_x|+\delta|\sigma v_xv|+\epsilon\delta|\sigma u_xv_x|)dxdy\notag\\
		\leq&C\epsilon^2\delta\int_{0}^{1}\int_{0}^{2\pi}u_y^2dxdy+C\delta\int_{0}^{1}\int_{0}^{2\pi}|\sigma|(u_x^2+v_x^2)dxdy.
	\end{align}
{\bf The source term:}
Finally, by Young inequality, it's easy to get 
	\begin{align}
		\int_{0}^{1}\int_{0}^{2\pi}&[F_1(\sigma u_x-\sigma'v)+F_2\sigma v_x]dxdy\notag\\
		\leq&\frac{c}{2}\int_{0}^{1}\int_{0}^{2\pi}|\sigma|(u_x^2+v_x^2)dxdy+C\int_{0}^{1}\int_{0}^{2\pi}(F_1^2+F_2^2)dxdy+C\int_{0}^{1}\int_{0}^{2\pi}v^2dxdy. \label{e:source term}
	\end{align}
Collecting the estimate (\ref{eq:poseq1})-(\ref{e:source term}), we deduce that there exist $\epsilon_0>0, \delta_0>0$ such that for any $\epsilon\in (0,\epsilon_0), \delta\in (0,\delta_0)$, the desired inequality holds, this completes the proof of this Lemma.

\end{proof}
	
\subsubsection{Energy estimate} Next we give the following energy estimate.
	\begin{lemma}\label{lem:energy}
		Suppose that $(u,v)\in D$ is the solution of $\eqref{eq:linearEstimate}$, then there exists $\epsilon_0>0,\delta_0>0$ such that for any $\epsilon\in(0,\epsilon_0),~ \delta\in(0,\delta_0)$, there holds
		\begin{align}
			\epsilon^2\int_{0}^{1}\int_{0}^{2\pi}(u_y^2+u_x^2+v_x^2)dxdy\leq&C\delta\int_{0}^{1}\int_{0}^{2\pi}|\sigma|u_x^2dxdy\notag\\
			&+C\int_{0}^{1}\int_{0}^{2\pi}(u^2+v^2)dxdy+C\int_{0}^{1}\int_{0}^{2\pi}(F_1^2+F_2^2)dxdy.\nonumber
		\end{align}
	\end{lemma}
\begin{proof} Multiplying $\eqref{eq:linearEstimate}_1$ by $u$, $\eqref{eq:linearEstimate}_2$ by $v$, adding them together, integrating on $\Omega$ and using the divergence-free condition, we obtain that
	\begin{align}
\epsilon^2\int_{0}^{1}\int_{0}^{2\pi}(u_x^2+u_y^2+v_x^2+v_y^2)dxdy+\int_{0}^{1}\int_{0}^{2\pi}(uS_u+vS_v)dxdy=\int_{0}^{1}\int_{0}^{2\pi}(uF_1+vF_2)dxdy.\label{eq:energy1}
	\end{align}
 Next, we deal with these terms related to $S_u$ and $S_v$.\\
\textbf{The term of $S_u$}: By integration by parts and divergence-free condition, we have
	\begin{align}
		\int_{0}^{1}\int_{0}^{2\pi}uS_udxdy&=\int_{0}^{1}\int_{0}^{2\pi}\left[\frac{1}{2}u^a(u^2)_x+\frac{1}{2}v^a(u^2)_y+u^2u_x^a+uvu_y^a\right]dxdy\notag\\
		&=\int_{0}^{1}\int_{0}^{2\pi}\left[u^2u_x^a+uvu_y^a\right]dxdy.\label{eq:energy2}
	\end{align}
	By Lemma $\ref{lem:approxEstimates}$, there holds
	\begin{align}
		\left|\int_{0}^{1}\int_{0}^{2\pi}u^2u_x^adxdy\right|\leq&C\delta\int_{0}^{1}\int_{0}^{2\pi}u^2dxdy.\nonumber
	\end{align}
	Besides, we decompose $\displaystyle\int_{0}^{1}\int_{0}^{2\pi} uvu_y^adxdy$ as follows
	\begin{align}
		\int_{0}^{1}\int_{0}^{2\pi}uvu_y^adxdy=\int_{0}^{1}\int_{0}^{2\pi}[(u^a-u_e(y))_yuv+u_e'(y)uv]dxdy.\nonumber
	\end{align}
	By Proposition $\ref{lem:approxEstimates}$ and Hardy inequality, we deduce that
	\begin{align}
		\Big|\int_{0}^{1}\int_{0}^{2\pi}&uv(u^a-u_e(y))_ydxdy\Big|\notag\\
		\leq&\left|\int_{0}^{1}\int_{0}^{2\pi}uv\partial_y(u_e^a-u_e(y))dxdy\right|
+\left|\int_{0}^{1}\int_{0}^{2\pi}uvh_ydxdy\right|\notag\\
		&+\left|\int_{0}^{1}\int_{0}^{2\pi}uv(\chi(y)u_p^a)_ydxdy\right|+\left|\int_{0}^{1}\int_{0}^{2\pi}uv((1-\chi(y))\hat{u}_p^a)_ydxdy\right|\notag\\
\leq& C\epsilon\delta\|u\|_2\|v\|_2+C\epsilon\int_{0}^{1}\int_{0}^{2\pi}|\sigma|\left(\left|\frac{uv}{(y-1)^2}\right||\zeta^2\partial_{\zeta}u_p^a|+\left|\frac{uv}{y^2}\right|\left|\eta^2\partial_{\eta}\hat{u}_p^a\right|\right)dxdy\notag\\
&+C\epsilon\int_{0}^{1}\int_{0}^{2\pi}(1-|\sigma|)\left(\left|\frac{uv}{(y-1)^2}\right||\zeta^2\partial_{\zeta}u_p^a|+\left|\frac{uv}{y^2}\right|\left|\eta^2\partial_{\eta}\hat{u}_p^a\right|\right)dxdy\notag\\
\leq& C\epsilon\delta\|u\|_2\|v\|_2+C\epsilon\delta\left(\int_{0}^{1}\int_{0}^{2\pi}\frac{u^2}{(1-y)^2}dxdy\right)^{\frac{1}{2}}\left(\int_{0}^{1}\int_{0}^{2\pi}|\sigma|\frac{v^2}{(1-y)^2}dxdy\right)^{\frac{1}{2}}\notag\\
&+C\epsilon\delta\left(\int_{0}^{1}\int_{0}^{2\pi}\frac{u^2}{y^2}dxdy\right)^{\frac{1}{2}}\left(\int_{0}^{1}\int_{0}^{2\pi}|\sigma|\frac{v^2}{y^2}dxdy\right)^{\frac{1}{2}}\notag\\
\leq&C\epsilon^2\delta\int_{0}^{1}\int_{0}^{2\pi}u_y^2dxdy+C\delta\int_{0}^{1}\int_{0}^{2\pi}|\sigma|v_y^2dxdy+C\delta\int_{0}^{1}\int_{0}^{2\pi}v^2dxdy.\notag
	\end{align}
	Thus, there holds
\begin{align}
\Big|\int_{0}^{1}\int_{0}^{2\pi}uvu_y^adxdy\Big|\leq&C\epsilon^2\delta\int_{0}^{1}\int_{0}^{2\pi}u_y^2dxdy\notag\\
&+C\delta\int_{0}^{1}\int_{0}^{2\pi}|\sigma|u_x^2dxdy+C\int_{0}^{1}\int_{0}^{2\pi}(u^2+v^2)dxdy. \label{second part energy}
\end{align}
 Collecting \eqref{eq:energy2}-\eqref{second part energy}, we deduce that
	\begin{align}
		\left|\int_{0}^{1}\int_{0}^{2\pi}uS_udxdy\right|\leq&C\epsilon^2\delta\int_{0}^{1}\int_{0}^{2\pi}u_y^2dxdy\notag\\
		&+C\delta\int_{0}^{1}\int_{0}^{2\pi}|\sigma|u_x^2dxdy+C\int_{0}^{1}\int_{0}^{2\pi}(u^2+v^2)dxdy.\label{est:Su}
	\end{align}
	\textbf{The term of $S_v$}: By integration by parts and divergence-free condition, we have
	\begin{align}
		\int_{0}^{1}\int_{0}^{2\pi}vS_vdxdy&=\int_{0}^{1}\int_{0}^{2\pi}\left[\frac{1}{2}u^a(v^2)_x+\frac{1}{2}v^a(v^2)_y+v^2v_y^a+uvv_x^a\right]dxdy\notag\\
		&=\int_{0}^{1}\int_{0}^{2\pi}\left[-v^2u_x^a+uvv_x^a\right]dxdy.\nonumber
	\end{align}
	By Lemma $\ref{lem:uxvxEstimate}$ and Hardy inequality, we have
	\begin{align}
		\left|\int_{0}^{1}\int_{0}^{2\pi}\left[-v^2u_x^a+uvv_x^a\right]dxdy\right|\leq C\epsilon^2\delta\int_{0}^{1}\int_{0}^{2\pi}u_y^2dxdy+C\delta\int_{0}^{1}\int_{0}^{2\pi}v^2dxdy.\label{est:Sv}
	\end{align}
	The force term can be estimated directly by Cauchy inequality, thus after combining $\eqref{eq:energy1}$, $\eqref{est:Su}$ and $\eqref{est:Sv}$, we complete the proof of this Lemma.
\end{proof}
	Combining the Lemma $\ref{lem:positivity}$ and Lemma $\ref{lem:energy}$, we obtain the following estimate:
	\begin{Proposition}\label{prop:wall estimate}
		Suppose that $(u,v)\in D$ is the solution of $\eqref{eq:linearEstimate}$, then there exists $\epsilon_0>0,\delta_0>0$ and a constant $C>0$ such that for any $ \epsilon\in(0,\epsilon_0),~ \delta\in(0,\delta_0)$, the following estimate holds
		\begin{align*}
			\int_{0}^{1}\int_{0}^{2\pi}&|\sigma|(u_x^2+v_x^2)dxdy+\epsilon^2\int_{0}^{1}\int_{0}^{2\pi}u_y^2dxdy\leq C\int_{0}^{1}\int_{0}^{2\pi}(u^2+v^2)dxdy+C\int_{0}^{1}\int_{0}^{2\pi}(F_1^2+F_2^2)dxdy.
		\end{align*}
	\end{Proposition}
	
\subsubsection{Interior $H^2$ estimate}
	We consider the following vorticity equation of $\eqref{eq:linearEstimate}$ 
	\begin{align}
		-\epsilon^2\Delta\omega+u^a\omega_x+v^a\omega_y=\partial_yF_1-\partial_xF_2-u\omega_x^a-v\omega_y^a,\label{eq:vorticity equation}
	\end{align}
	where $\omega=u_y-v_x,\omega^a=u_y^a-v_x^a.$

\begin{lemma}\label{gamma cutoff function}
	For any  $0<y_1<y^*<y_2<1$, there exist $\epsilon_0>0$ such that the following holds:

for any $\epsilon\in(0,\epsilon_0),\delta\in(0,1)$, there exists a cutoff $\gamma\in C_c^\infty(\Omega)$ satisfying $1\leq\gamma(x,y)\leq 3,\forall(x,y)\in[0,2\pi]\times [y_1,y_2]$,
		\begin{align}
			u_e^a\gamma_x+v_e^a\gamma_y=0,\label{eq:gamma equation}
		\end{align}
		and a constant $C>0$ exists such that
		\begin{align*}
			\|\gamma\|_{C^2(\Omega)}\leq C.
		\end{align*}
\end{lemma}
\begin{proof}
we choose a sufficient small constant $\bar{c}>0$ which is dependent on $y_1,y_2$
		and denote $\displaystyle\phi_e(y)=\int_{y^*}^{y}(az+b)dz$ and $\phi_e^a(x,y)=\displaystyle\int_{y^*}^{y}u_e^a(x,z)dz-\int_{0}^{x}v_e^a(\bar{x},y^*)d\bar{x}$. Define $\displaystyle\gamma_j=g_j(\phi_e^a)\chi_j,j=1,2,3$ and
		\begin{align*}
			&g_1(t)=\begin{cases}
				1,t\in[\phi_e(y^*-\frac{3}{2}\bar{c}),\phi_e(\frac{3}{2}\bar{c})],\\
				0,t\in[\phi_e(\bar{c}),+\infty)\cup(-\infty,\phi_e(y^*-\bar{c})].
			\end{cases}
			\quad\hspace{0.5em}\quad\chi_1(y)=\begin{cases}
				1,y\in[\frac{3}{4}\bar{c},y^*-\frac{3}{4}\bar{c}],\\
				0,y\in[0,\frac{\bar{c}}{2}]\cup[y^*-\frac{\bar{c}}{2},1].
			\end{cases}\\
			&g_2(t)=\begin{cases}
				1,t\in[\phi_e(y^*+\frac{3}{2}\bar{c}),\phi_e(1-\frac{3}{2}\bar{c})],\\
				0,t\in[\phi_e(1-\bar{c}),+\infty)\cup(-\infty,\phi_e(y^*+\bar{c})].
			\end{cases}
			\quad\chi_2(y)=\begin{cases}
				1,y\in[y^*+\frac{3}{4}\bar{c},1-\frac{3}{4}\bar{c}],\\
				0,y\in[1-\frac{\bar{c}}{2},1]\cup[0,y^*+\frac{\bar{c}}{2}].
			\end{cases}\\
			&g_3(t)=\begin{cases}
				1,t\leq\phi_e(y^*\pm2\bar{c}),\\
				0,t\geq\phi_e(y^*\pm\frac{5}{2}\bar{c}).
			\end{cases}
			\quad\hspace{10em}\chi_3(y)=1.
		\end{align*}
		are smooth functions with values in $[0,1]$.We define $\gamma=\gamma_1+\gamma_2+\gamma_3$.
		
		Given $y\in[y_1,y_2]$, then $y\in[2\bar{c},1-2\bar{c}]$. If $y\in[2\bar{c},y^*-\frac{7}{4}\bar{c}]$, then $\phi_e^a=\phi_e(y)+O(\epsilon\delta)\in[\phi_e(y^*-\frac{3}{2}\bar{c}),\phi_e(\frac{3}{2}\bar{c})]$, so $\gamma_1=1$. Similarly, if $y\in[y^*-\frac{7}{4}\bar{c},1-2\bar{c}]$, then $\phi_e^a\in[\phi_e(y^*+\frac{3}{2}\bar{c}),\phi_e(1-\frac{3}{2}\bar{c})]$, so $\gamma_2=1$. If $y\in[y^*-\frac{7}{4}\bar{c},y^*+\frac{7}{4}\bar{c}]$, then $\phi_e^a\leq\phi_e(y^*\pm2\bar{c})$, so $\gamma_3=1$. This proves $1\leq\gamma(x,y)\leq 3,\forall(x,y)\in[0,2\pi]\times [y_1,y_2]$.
		
		For any $y\in[0,1]$,if $y\in[0,\frac{7}{8}\bar{c}]$, then $\phi_e^a\geq\phi_e(\bar{c})$, so $\gamma=0$. Similarly, if $y\in[1-\frac{7}{8}\bar{c},1]$, then $\phi_e^a\geq\phi_e(1-\bar{c})$, so $\gamma=0$. If $y\in[\frac{7}{8}\bar{c},y^*-\frac{7}{8}\bar{c}]$, then $\phi_e^a\in[\phi_e(y^*-\frac{3}{4}\bar{c}),\phi_e(\frac{3}{4}\bar{c})]$, so $\chi_1=1,\chi_2=0,\gamma=g_1(\phi_e^a)+g_3(\phi_e^a)$. Similarly if $y\in[y^*+\frac{7}{8}\bar{c},1-\frac{7}{8}\bar{c}]$, then $\phi_e^a\in[\phi_e(y^*+\frac{3}{4}\bar{c}),\phi_e(1-\frac{3}{4}\bar{c})]$, so $\chi_1=0,\chi_2=1,\gamma=g_2(\phi_e^a)+g_3(\phi_e^a)$. If $y\in[y^*-\frac{7}{8}\bar{c},y^*+\frac{7}{8}\bar{c}]$, then $\phi_e^a\leq\phi_e(y^*\pm\bar{c})$, so $g_1(\phi_e^a)=g_2(\phi_e^a)=0,\gamma=g_3(\phi_e^a)$. Hence, $\eqref{eq:gamma equation}$ follows from 
		$$u_e^a\partial_x(g_j(\phi_e^a))+v_e^a\partial_y(g_j(\phi_e^a))=0,j=1,2,3.$$
\end{proof}

	\begin{Proposition}\label{prop:interior estimate}
		Let $(u,v)\in D$ be the solution of $\eqref{eq:linearEstimate}$. There exists a constant \(\delta_0>0\) such that for any \(\delta\in(0,\delta_0)\) and any compact subset \(K\subseteq [0,2\pi]\times (0,1)\), there exist constants \(\epsilon_K>0\) and \(C_K>0\) such that for any \(\epsilon\in(0,\epsilon_K)\), the following estimate holds:
		\begin{align*}
			\|(u,v)\|_{H^2(K)} \leq C_K \left( \|(u,v)\|_2 + \epsilon^{-2}\|(F_1, F_2)\|_2 \right).
		\end{align*}
	\end{Proposition}
	\begin{proof}
		Suppose that $K\subseteq[0,2\pi]\times[y_1,y_2]$ where $y_1,y_2$ are independent of $\epsilon$ and satisfy $0<y_1<y^*<y_2<1$. By Lemma \ref{gamma cutoff function}, we can choose a smooth cutoff function $\gamma$ such that $1\leq\gamma\leq3\text{ on }K$, $\|\gamma\|_{C^2(\Omega)}\leq C$ and 
		\begin{align}
			u_e^a\gamma_x+v_e^a\gamma_y=0.\label{eq:gamma}
		\end{align}
	 Noticing that outer Euler flow has a constant vorticity, we can rewrite $\eqref{eq:vorticity equation}$ as the form 
		\begin{align}
			-\epsilon^2\Delta\omega+u_e^a\omega_x+v_e^a\omega_y=\partial_yF_1-\partial_xF_2+R_1,\label{eq:interior vorticity equation}
		\end{align}
		where 
		\begin{align*}
			R_1=-(u^a-u_e^a)\omega_x-(v^a-v_e^a)\omega_y-u(\omega^a-\omega_e^a)_x-v(\omega^a-\omega_e^a)_y.
		\end{align*}
		Multiplying the vorticity equation $\eqref{eq:interior vorticity equation}$ with $\gamma^4\omega$, integrating in $\Omega$ and using $\eqref{eq:gamma}$, we obtain that
		\begin{align}
			\epsilon^2\int_{0}^{1}\int_{0}^{2\pi}\gamma^4|\nabla\omega|^2dxdy=&\int_{0}^{1}\int_{0}^{2\pi}(\partial_yF_1-\partial_xF_2)\gamma^4\omega dxdy\notag\\
			&+\frac{\epsilon^2}{2}\int_{0}^{1}\int_{0}^{2\pi}\Delta(\gamma^4)\omega^2dxdy+\int_{0}^{1}\int_{0}^{2\pi}R_1\gamma^4\omega dxdy,\label{eq:interior estimate}
		\end{align}
		We estimate them term by term. The forcing term can be estimated as
		\begin{align}
			\int_{0}^{1}\int_{0}^{2\pi}(\partial_yF_1-\partial_xF_2)\gamma^4\omega dxdy=&\int_{0}^{1}\int_{0}^{2\pi}[(4\gamma^3\gamma_x\omega+\gamma^4\omega_x)F_2-(4\gamma^3\gamma_y\omega+\gamma^4\omega_y)F_1]dxdy\notag\\
			\leq&C\|(F_1,F_2)\|_2(\|\gamma\omega\|_2+\|\gamma^2\nabla\omega\|_2)\label{ineq:force interior estimate}
		\end{align}
		By Lemma $\ref{lem:approxEstimates}$, the term can by bounded by
		\begin{align}
			\int_{0}^{1}\int_{0}^{2\pi}R_1\gamma^4\omega dxdy\leq&C_M\epsilon^M(\|\gamma\omega\|_2^2+\|\gamma^2\nabla\omega\|_2^2+\|(u,v)\|_2^2),\label{ineq:R interior estimate}
		\end{align}
		where $M$ can be any positive integer.
		For the last term, by direct calculation, we arrive at
		\begin{align}
			\int_{0}^{1}\int_{0}^{2\pi}\Delta(\gamma^4)\omega^2dxdy=&\int_{0}^{1}\int_{0}^{2\pi}(12\gamma^2|\nabla\gamma|^2+4\gamma^3\Delta\gamma)\omega^2dxdy\leq C\|\gamma\omega\|_2^2.\label{ineq:Deltagamma estimate}
		\end{align} 
		Combining $\eqref{eq:interior estimate}-\eqref{ineq:Deltagamma estimate}$ and using Young inequality,we have
		\begin{align}
			\epsilon\|\gamma^2\nabla\omega\|_2\leq&C\epsilon^{-1}\|(F_1,F_2)\|_2+C\epsilon\|\gamma\omega\|_2+C\epsilon\|(u,v)\|_2.\label{ineq:vorticity interior estimate}
		\end{align}
	Moreover, there holds
		\begin{align*}
			\|\gamma\omega\|_2^2=&\int_{0}^{1}\int_{0}^{2\pi}\gamma^2\omega(u_y-v_x)dxdy\\
			=&\int_{0}^{1}\int_{0}^{2\pi}[v(2\gamma\gamma_x\omega+\gamma^2\omega_x)-u(2\gamma\gamma_y\omega+\gamma^2\omega_y)]dxdy\\
			\leq&C\|(u,v)\|_2(\|\gamma\omega\|_2+\|\gamma^2\nabla\omega\|_2).
		\end{align*}
		Hence by Young inequality,
		\begin{align*}
			\|\gamma\omega\|_2^2\leq\nu\|\gamma^2\nabla\omega\|_2^2+C_{\nu}\|(u,v)\|_2^2.
		\end{align*}
		Inserting this into $\eqref{ineq:vorticity interior estimate}$ and choosing small $\nu>0$, we obtain that
		\begin{align*}
			\|\gamma^2\nabla\omega\|_2\leq&C\epsilon^{-2}\|(F_1,F_2)\|_2+C\|(u,v)\|_2.
		\end{align*}
		Finally, notice that $\omega_x=-\Delta v,\omega_y=\Delta u$ and use the classical interior estimate for elliptic equation, we get that
		\begin{align*}
			\|(u,v)\|_{H^2(K)}\leq&C_K(\|\nabla\omega\|_{L^2(\tilde{K})}+\|(u,v)\|_2)\leq C_K(\|(u,v)\|_2+\epsilon^{-2}\|(F_1,F_2)\|_2),
		\end{align*}
		where $\tilde{K}$ is a compact subset of $\Omega$ containing $K$.
	\end{proof}
	Combining Proposition $\ref{prop:wall estimate}$ and Proposition $\ref{prop:interior estimate}$, we have the following result.
	\begin{Corollary}\label{cor:v gradient estimate}
		There exists $\epsilon_0>0,\delta_0>0$ and a constant $C>0$ such that for any $ \epsilon\in(0,\epsilon_0),~ \delta\in(0,\delta_0)$, then the solution $(u,v)\in D$ of $\eqref{eq:linearEstimate}$ satisfies
		\begin{align*}
			\int_{0}^{1}\int_{0}^{2\pi}(u_x^2+v_x^2)dxdy\leq&C\int_{0}^{1}\int_{0}^{2\pi}(u^2+v^2)dxdy+C\epsilon^{-4}\int_{0}^{1}\int_{0}^{2\pi}(F_1^2+F_2^2)dxdy.
		\end{align*}
	\end{Corollary}

\subsubsection{$L^2$ estimate}

\begin{Proposition}\label{prop:uv estimate}
		Suppose that $(u,v)\in D$ is the solution of $\eqref{eq:linearEstimate}$, $\delta_0$ is the constant in Corollary $\ref{cor:v gradient estimate}$. Fix $\delta\in(0,\delta_0)$, then there exist $\epsilon_0>0$ and a constant $C>0$ such that for any $ \epsilon\in(0,\epsilon_0)$, the following estimate holds
		\begin{align*}
			\int_{0}^{1}\int_{0}^{2\pi}(u^2+v^2)dxdy\leq C\epsilon^{-4}\int_{0}^{1}\int_{0}^{2\pi}(F_1^2+F_2^2)dxdy.
		\end{align*}
	\end{Proposition}

	The proof of Proposition $\ref{prop:uv estimate}$ is very complicated, we first give a outline of idea.  Assume that Proposition $\ref{prop:uv estimate}$ fails as $\epsilon\to0$, then there would be a sequence, indexed implicitly below, such that
	\begin{align}
		\epsilon\to0,~~\|(u,v)\|_2=1,~~\|(F_1,F_2)\|_2=o(\epsilon^2).\label{blow up assumption}
	\end{align}
Proposition $\ref{prop:wall estimate}$, Proposition $\ref{prop:interior estimate}$ and Corollary $\ref{cor:v gradient estimate}$ give that for any compact subset $K$ of $\Omega$,
	\begin{align}
		\|(u_x,v_x)\|_2\leq&C,~~\|(u,v)\|_{H^2(K)}\leq C_K\label{ineq:limit bound}
	\end{align}
	if the index is large enough. After extraction, $(u,v)\to(u_0,v_0)$ weakly in $H_{loc}^2(\Omega)$,  $(u,v)\to(u_0,v_0)$ strongly in $H^1(K)$ on every interior compact set, and the corresponding vorticities $w$ converge to $\omega_0=\partial_yu_0-\partial_xv_0$ weakly in $H_{\operatorname{loc}}^1$. Thus, by $\eqref{eq:vorticity equation}$, we deduce 
	\begin{align}
		\int_{0}^{1}\int_{0}^{2\pi}u_e(y)\partial_x\omega_0\psi dxdy=0,\label{eq:limiting vorticity equation}
	\end{align}
	for any periodic $\psi\in C_c^{\infty}([0,2\pi]\times(0,1))$. Since $u_e(y)$ vanishes on just one horizontal line and $\partial_x\omega_0\in L_{\operatorname{loc}}^2$, it follows that $\partial_x\omega_0=0$ in $[0,2\pi]\times (0,1)$. Moreover, by Corollary $\ref{cor:v gradient estimate}$, we deduce that $v_0\in H_0^1$. Since $\Delta v_0=-\partial_x\omega_0=0$, uniqueness for the periodic problem implies $v_0=0$. Thus $u_0$ depends only on $y$. The convergence $v\to0$ is strong in global $L^2$ in Rellich compactness.
	
	The strong convergence of $u\to0$ in $L^2$ can be proved by the following two steps.

{\bf Step 1:} We prove that $u_0(y)=0$. \\

{\bf Step 2:} An exact collar circulation equation then rules out the remaining possibility that its $L^2$ mass concentrates near the wall. Together, they constitute a contradiction and Proposition $\ref{prop:uv estimate}$ follows consequently.\\

	Let's first carry out the first step. We first show that $u_0(y)$ is a linear function, that is $u_0(y)=cy+d$. Then, we show that $c=d=0$ by choosing suitable wall adjoint tests.

	\begin{Lemma}\label{lem:affine}
		The interior limit has the form
		$$(u_0,v_0)=(cy+d, 0).$$
	\end{Lemma}
	\begin{proof}
		We show that for any $h\in C_c^{\infty}((y_\ast,1))$. Then there holds
		\begin{align}
			\int_{y^*}^{1}u_0(y)h''(y)dy=0.\label{eq:affine test}
		\end{align}
		Let 
		\begin{align*}
			f'(\phi_e(y))=\frac{h(y)}{u_e(y)},~~ y\in(y^*,1],
		\end{align*}
		then extend it to the entire channel by zero. It's easy to see that for any compact subset $K\subseteq [0,2\pi]\times(y^*,1)$, there holds in $C^k(K)$
		\begin{align}
			(Z_1,Z_2):=(\partial_y(f(\phi_e^a)),-\partial_x(f(\phi_e^a)))\rightarrow(h(y),0).\label{eq:multiplier convergence}
		\end{align}

Multiplying $\eqref{eq:linearEstimate}_1$ by $Z_1$, $\eqref{eq:linearEstimate}_2$ by $Z_2$, adding them together, integrating on $\Omega$ and using the divergence-free condition, we obtain that
		\begin{align}
			\langle-\epsilon^2\Delta u,Z_1\rangle+\langle-\epsilon^2\Delta v,Z_2\rangle+E=\langle F_1,Z_1\rangle+\langle F_2,Z_2\rangle+R_2,\label{eq:affine equation}
		\end{align}
		where $\left\langle f,g\right\rangle$ denotes $\int_{y^*}^{1}\int_{0}^{2\pi}fgdxdy$ and
		\begin{align*}
			E=&\langle u_e^au_x+v_e^au_y,Z_1\rangle+\langle u_e^av_x+v_e^av_y,Z_2\rangle+\langle uu_{ex}^a+vu_{ey}^a,Z_1\rangle+\langle uv_{ex}^a+vv_{ey}^a,Z_2\rangle,\\[3pt]
			R_2=&\langle-(u^a-u_e^a)u_x-(v^a-v_e^a)u_y-u(u^a-u_e^a)_x-v(u^a-u_e^a)_y,Z_1\rangle\\[3pt]
			&+\langle-(u^a-u_e^a)v_x-(v^a-v_e^a)v_y-u(v^a-v_e^a)_x-v(v^a-v_e^a)_y,Z_2\rangle.
		\end{align*}
		By integration by parts, the constant outer vorticity and divergence-free condition, we have
		\begin{align}
			E=&-[\langle u_e^au,\partial_xZ_1\rangle+\langle v_e^au,\partial_yZ_1\rangle+\langle u_e^av,\partial_xZ_2\rangle+\langle v_e^av,\partial_yZ_2\rangle]\notag\\
			&+\langle uu_{ex}^a+vu_{ey}^a,Z_1\rangle+\langle uv_{ex}^a+vv_{ey}^a,Z_2\rangle\notag\\
			=&-[\langle u_e^au,f'u_{ex}^a\rangle+\langle v_e^au,f'u_{ey}^a\rangle+\langle u_e^av,f'v_{ex}^a\rangle+\langle v_e^av,f'v_{ey}^a\rangle]\notag\\
			&+\langle uu_{ex}^a+vu_{ey}^a,f'u_e^a\rangle+\langle uv_{ex}^a+vv_{ey}^a,f'v_e^a\rangle\notag\\
			&+\langle u_e^au,f''u_e^av_e^a\rangle-\langle v_e^au,f''(u_e^a)^2\rangle+\langle u_e^av,f''(v_e^a)^2\rangle-\langle v_e^av,f''u_e^av_e^a\rangle\notag\\
			=&\omega_e^a\langle v,f'u_e^a\rangle-\omega_e^a\langle u,f'v_e^a\rangle
			=-\omega_e^a\langle f(\phi_e^a),u_x+v_y\rangle=0.\label{eq:affine Euler term}
		\end{align}
		By the properties of approximate solutions, $\eqref{blow up assumption}$ and $\eqref{ineq:limit bound}$, we deduce 
		\begin{align}
			\langle F_1,Z_1\rangle+\langle F_2,Z_2\rangle+R_2=o(\epsilon^2). \label{source term estimate}
		\end{align}
		Combining $\eqref{eq:affine equation}$-$\eqref{source term estimate}$, we deduce that
		\begin{align*}
			\langle u,\Delta Z_1\rangle+\langle v,\Delta Z_2\rangle=o(1).
		\end{align*}
		Therefore, noticing $\eqref{eq:multiplier convergence}$, we obtain $\eqref{eq:affine test}$ by taking the limit $\epsilon\to0$. Hence $u_0(y)$ is a linear function in $(y^*,1)$.

 By the same argument,  $u_0(y)$ is also a linear function in $(0,y^*)$. The local $H^2$ bound across the critical level $y=y^*$ matches both the value and first derivative, so $u_0(y)=cy+d, ~~\forall y\in (0,1)$. \\
	\end{proof}
	
	We then choose a delicate wall adjoint test to eliminate the affine mode. We carry out the upper-wall construction and introduce a variable
	$$t(x,y)=\phi^a(x,1)-\phi_e^a(x,y),y\in(y^*,1)$$
	where $\phi^a(x,y)=\displaystyle\int_{y^*}^{y}u^a(x,z)dz-\int_{0}^{x}v^a(\bar{x},y^*)d\bar{x}$ is the stream function of $(u^a,v^a)$. Obviously,
	\begin{align*}
		t=(\phi^a(x,1)-\phi_e^a(x,1))+(\phi_e^a(x,1)-\phi_e^a(x,y))=\frac{a}{2}(1-y^2)+b(1-y)+O(\epsilon\delta).
	\end{align*}
	Let $f_\ast \in C^\infty((-1,1))$ and satisfies 
\[
f_*(t)=
\begin{cases}
\displaystyle\frac{t^2}{2u_e(1)}, ~~|t|\leq\bar{t},\\[10pt]
0,\quad\quad\hspace{0.5em} |t|\geq2\cdot\bar{t}.
\end{cases}
\]
We choose sufficiently small $\bar{t}>0$ such that the support of $f_*(t(x,y))$ is contained in $[0,2\pi]\times[1-10^{-3},1]$. Extend it by zero through a region where $f_*(t(x,y))$ already vanishes. This ensures both smoothness and
	\begin{align}
		u_e^a\partial_x(f_*(t))+v_e^a\partial_y(f_*(t))=0.\label{eq:f* equality}
	\end{align}
	In the upper inner variable, we have
	\begin{align*}
		&u^a(x,1+\epsilon\zeta)=u_p^{(0)}+u_e(1)+O(\epsilon(1+|\zeta|)),\\
		&v^a(x,1+\epsilon\zeta)=\epsilon(v_p^{(1)}+v_e^{(1)}(x,1))+O(\epsilon^2\delta(1+|\zeta|)).
	\end{align*}
	More importantly, there holds
	\begin{align*}
		t(x,1+\epsilon\zeta)=&\int_{y^*}^{1}(u^a-u_e^a)(x,z)dz+\int_{y}^{1}u_e^a(x,z)dz-\int_{0}^{x}(v^a(\bar{x},y^*)-v_e^a(\bar{x},y^*))d\bar{x}\notag\\
		=&\int_{y^*}^{1}\chi(z)u_p^a(x,\frac{z-1}{\epsilon})dz+\int_{y^*}^{1}(1-\chi(z))\hat{u}_p^a(x,\frac{z}{\epsilon})dz+\int_{y^*}^{1}h(x,z)dz+\int_{y}^{1}u_e^a(x,z)dz\notag\\
		&-\int_{0}^{x}(\chi v_p^a+(1-\chi)\hat{v}_p^a)(\bar{x},y^*)d\bar{x}\notag\\
		=&\int_{y^*}^{1}\partial_z\left(\epsilon\chi(z)\int_{-\infty}^{\frac{z-1}{\epsilon}}u_p^a(x,\bar{\zeta})d\bar{\zeta}\right)dz-\int_{y^*}^{1}\partial_z\left(\epsilon(1-\chi(z))\int_{\frac{z}{\epsilon}}^{+\infty}\hat{u}_p^a(x,\bar{\eta})d\bar{\eta}\right)dz\notag\\
		&+\int_{y}^{1}(az+b)dz+O(\epsilon^2(1+|\zeta|))\notag\\
		=&\epsilon\left(\int_{-\infty}^{0}u_p^{(0)}(x,\bar{\zeta})d\bar{\zeta}-u_e(1)\zeta\right)+O(\epsilon^2(1+\zeta^2)).
	\end{align*}
	Denote $\Phi_p(x,\zeta)=\displaystyle\int_{-\infty}^{\zeta}u_p^{(0)}(x,\bar{\zeta})d\bar{\zeta}$, then
	\begin{align}
		f_*(t(x,1+\epsilon\zeta))=&\frac{1}{2u_e(1)}\left[\epsilon\left(\Phi_p(x,0)-u_e(1)\zeta\right)+O(\epsilon^2(1+\zeta^2))\right]^2\notag\\
		=&\epsilon^2\left(\frac{u_e(1)}{2}\zeta^2-\zeta\Phi_p(x,0)+\frac{\Phi_p^2(x,0)}{2u_e(1)}\right)+O(\epsilon^3(1+|\zeta|^3)).\label{eq:f* upper expansion}
	\end{align}
	Let $\chi_b(y)\in C^\infty([0,1])$ satisfies 
\[
\chi_b(y)=
\begin{cases}
1, ~~y\in[1-10^{-3},1],\\[5pt]
0, y\in [0, 1-2\cdot10^{-3}],
\end{cases}
\]
then we define
	\begin{align}
		\theta_{\operatorname{pre}}(x,y)=&f_*(t)+\epsilon^2\chi_b(y)\left(-\frac{\Phi_p^2(x,0)}{2u_e(1)}
-\int_{-\infty}^{0}\Phi_p(x,\zeta)d\zeta\right)+\epsilon^2\chi_b(y)\int_{-\infty}^{\zeta}\Phi_p(x,\bar{\zeta})d\bar{\zeta}.\label{theta pre}
	\end{align}  
	By $\eqref{eq:f* upper expansion}$, we deduce that 
	\begin{align}
		\theta_{\operatorname{pre}}=&\epsilon^2\left(\int_{\zeta}^{0}\int_{\bar{\zeta}}^{0}u_e(1)dZd\bar{\zeta}+\int_{\zeta}^{0}\int_{-\infty}^{0}u_p^{(0)}(x,Z)dZd\bar{\zeta}-\int_{\zeta}^{0}\int_{-\infty}^{\bar{\zeta}}u_p^{(0)}(x,Z)dZd\bar{\zeta}\right)+O(\epsilon^3(1+|\zeta|^3))\notag\\
		=&\epsilon^2\int_{\zeta}^{0}\int_{\bar{\zeta}}^{0}(u_e(1)+u_p^{(0)}(x,Z))dZd\bar{\zeta}+O(\epsilon^3(1+|\zeta|^3)).\label{theta pre expansion}
	\end{align}
	Later we will see that the contribution of the principal term of $\theta_{\operatorname{pre}}$ vanishes due to the Prandtl equation. 
	
	Choose compactly supported smooth functions $r_0,r_1$ on $(-\infty,0]$ with
	\begin{align*}
		r_0(0)=1,~~r_0'(0)=0;\qquad  r_1(0)=0,~~r_1'(0)=1,
	\end{align*}
	and define the exact test
	\begin{align}
		\theta(x,y)=\theta_{\operatorname{pre}}-\theta_{\operatorname{pre}}(x,1)r_0(\zeta)-\epsilon\partial_y\theta_{\operatorname{pre}}(x,1)r_1(\zeta).\label{exact test}
	\end{align}
	Then $\theta(x,1)=\theta_y(x,1)=0$. Moreover, it follows from $\eqref{theta pre expansion}$ that for any $j\in\mathbb{N}$,
	\begin{align}
		\|\partial_x^j\theta_{\operatorname{pre}}(x,1)\|_{L^{\infty}([0,2\pi])}\leq C_j\epsilon^3,~~\|\partial_x^j\partial_y\theta_{\operatorname{pre}}(x,1)\|_{L^{\infty}([0,2\pi])}\leq C_j\epsilon^2,~~\|\nabla\theta\|_2\leq C.\label{ineq:exact test}
	\end{align}
	Let $\phi=\displaystyle\int_{1}^{y}u(x,z)dz$ be the stream function of $(u,v)$, then we have the following Lemma. 
	\begin{Lemma}\label{lem:matched adjoint residual}
		The test in $\eqref{exact test}$ satisfies
		\begin{align}
			\int_{0}^{1}\int_{0}^{2\pi}\phi S_{\epsilon}\theta dxdy=\int_{0}^{1}\int_{0}^{2\pi}(F_1\theta_y-F_2\theta_x)dxdy,\label{eq:matched adjoint equality}
		\end{align}
		where
		$$S_{\epsilon}\theta=\epsilon^2\Delta^2\theta+\Delta(u^a\theta_x+v^a\theta_y)+\theta_y\omega_x^a-\theta_x\omega_y^a,$$
		and satisfies
		$$S_{\epsilon}\theta=\epsilon^2 G_{\epsilon}+R_{\epsilon}.$$
Moreover, there hold 
\begin{align}\label{convergence of G}
G_{\epsilon}\to G_0:=\Delta^2\theta_0+\Delta\left(u_e(y)\chi_b(y)H_x\right),
\end{align}
		where 
		$$H(x)=-\frac{\Phi_p^2(x,0)}{2u_e(1)}-\int_{-\infty}^{0}\Phi_p(x,\zeta)d\zeta,$$
		and $\theta_0=f_*\left(\frac{a}{2}(1-y^2)+b(1-y)\right)$ on the selected branch. The residual term satisfies
		\begin{align}
			|R_{\epsilon}(x,1+\epsilon\zeta)|\leq&C\epsilon(1+|\zeta|)^{-m}.\label{ineq:exact test residual estimate}
		\end{align}
		Here $m>3$ can be fixed sufficiently large and $R_{\epsilon}$ is supported in the upper collar,i.e.$\operatorname{supp}R_{\epsilon}\subseteq[0,2\pi]\times[1-2\cdot10^{-3},1]$.
	\end{Lemma}
	\begin{proof}
		Multiplying $\eqref{eq:linearEstimate}_1$ by $\theta_y$, $\eqref{eq:linearEstimate}_2$ by $-\theta_x$, adding them together, integrating on $\Omega$ and using $\theta(x,1)=\theta_y(x,1)=0$, we obtain $\eqref{eq:matched adjoint equality}$ after integration by parts.
		
		For simplicity, we denote
		\begin{align*}
			K=&\int_{-\infty}^{\zeta}\int_{-\infty}^{\bar{\zeta}}u_p^{(0)}(x,Z)dZd\bar{\zeta},\\
			\Theta=&\int_{\zeta}^{0}\int_{\bar{\zeta}}^{0}(u_e(1)+u_p^{(0)}(x,Z))dZd\bar{\zeta}.
		\end{align*}
		Then  $\eqref{exact test}$ and $\eqref{eq:f* equality}$ give
		\begin{align}
			S_{\epsilon}\theta=&\epsilon^2\Delta^2\theta_{\operatorname{pre}}+\Delta(u^a\partial_x\theta_{\operatorname{pre}}+v^a\partial_y\theta_{\operatorname{pre}})+(\partial_y\theta_{\operatorname{pre}}\omega_x^a-\partial_x\theta_{\operatorname{pre}}\omega_y^a)+O(\epsilon\rho_m)\notag\\
			=&\underbrace{\epsilon^2\Delta^2f_*(t)+\epsilon^4\Delta^2(\chi_bH)+\epsilon^2\Delta(u_e^a\chi_bH_x+v_e^a\chi_b'H)}_{:=\epsilon^2G_{\epsilon}}
+\underbrace{\sum_{i=3}^{5}R_i+O(\epsilon\rho_m)}_{R_\epsilon},\label{eq:Sepsilon decomposition}
		\end{align}
		where $\rho_m=(1+|\zeta|)^{-m}$ and 
		\begin{align*}
			R_3=&\epsilon^4\Delta^2(\chi_bK)+\epsilon^2\Delta(u_e^a\chi_bK_x+v_e^a\chi_bK_y),\\
			R_4=&\Delta((u^a-u_e^a)\partial_x\theta_{\operatorname{pre}}+(v^a-v_e^a)\partial_y\theta_{\operatorname{pre}}),\\
			R_5=&\partial_y\theta_{\operatorname{pre}}(\omega^a-\omega_e^a)_x-\partial_x\theta_{\operatorname{pre}}(\omega^a-\omega_e^a)_y.
		\end{align*}
Since the functions in $G_{\epsilon}$ are smooth and 
		$$t\to\frac{a}{2}(1-y^2)+b(1-y),~~u_e^a\to u_e(y),~~v_e^a\to0(\epsilon\to0),$$
		thus (\ref{convergence of G}) holds.

Moreover, direct computation gives that 
		\begin{align*}
			R_3=&\epsilon^4\Delta^2K+\epsilon^2\Delta(u_e^aK_x+v_e^aK_y)+O(\epsilon\rho_m)\\[3pt]
			=&\partial_{\zeta}^4K+\partial_{\zeta}^2(u_e^aK_x+\epsilon^{-1}v_e^aK_{\zeta})+O(\epsilon\rho_m)\\[3pt]
			=&\partial_{\zeta}^2u_p^{(0)}+u_e(1)\partial_xu_p^{(0)}+v_e^{(1)}(x,1)\partial_{\zeta}u_p^{(0)}+O(\epsilon\rho_m).
		\end{align*}
		By $\eqref{theta pre expansion}$, we deduce that 
		\begin{align*}
			R_4=&\Delta[\epsilon^2(u^a-u_e^a)\Theta_x+\epsilon(v^a-v_e^a)\Theta_{\zeta}]+O(\epsilon\rho_m)\\[3pt]
			=&\partial_{\zeta}^2((u^a-u_e^a)\Theta_x)+\partial_{\zeta}^2(\epsilon^{-1}(v^a-v_e^a)\Theta_{\zeta})+O(\epsilon\rho_m)\\[3pt]
			=&u_p^{(0)}\partial_xu_p^{(0)}+v_p^{(1)}\partial_{\zeta}u_p^{(0)}-2\partial_{\zeta}u_p^{(0)}(v_p^{(1)}+v_e^{(1)}(x,1))\\[3pt]
			&+\partial_{\zeta}^2u_p^{(0)}\Theta_x+\partial_{\zeta}^2v_p^{(1)}\Theta_{\zeta}+2\partial_{\zeta}v_p^{(1)}(u_e(1)+u_p^{(0)})+O(\epsilon\rho_m),\\[3pt]
			R_5=&\epsilon^{-1}\epsilon^2\Theta_{\zeta}\epsilon^{-1}\partial_{x\zeta}u_p^{(0)}-\epsilon^2\Theta_x\epsilon^{-2}\partial_{\zeta}^2u_p^{(0)}+O(\epsilon\rho_m)\\[3pt]
			=&\partial_{x\zeta}u_p^{(0)}\Theta_{\zeta}-\Theta_x\partial_{\zeta}^2u_p^{(0)}+O(\epsilon\rho_m).
		\end{align*}
		Thus, using $\eqref{eq:up0}$ and $v_p^{(1)}(x,0)=-v_e^{(1)}(x,1)$, we deduce 
		\begin{align*}
R_{\epsilon}=&\partial_{\zeta}^2u_p^{(0)}-(u_e(1)+u_p^{(0)})\partial_xu_p^{(0)}-(v_p^{(1)}+v_e^{(1)}(x,1))\partial_{\zeta}u_p^{(0)}+O(\epsilon\rho_m)=O(\epsilon\rho_m),
		\end{align*}
		this is (\ref{ineq:exact test residual estimate}).
	\end{proof}
	\begin{Lemma}\label{lem:weighted Hardy infty}
		Let  $\phi\in C^1([0,1])$ with $\phi(1)=0$, then for any $m\geq 2$, there holds 
		\begin{align*}
			\int_{0}^{1}\left(1+\left|\frac{y-1}{\epsilon}\right|\right)^{-m}\phi^2(y)dy\leq&C_m\epsilon^2\int_{0}^{1}(\phi'(y))^2dy,\quad\int_{0}^{1}\left(1+\left|\frac{y-1}{\epsilon}\right|\right)^{-m}dy\leq C_m\epsilon.
		\end{align*}
		\begin{proof}
			Direct calculation gives
			\begin{align*}
				\int_{0}^{1}\left(1+\left|\frac{y-1}{\epsilon}\right|\right)^{-m}dy=&\epsilon\int_{-\infty}^{-\frac{1}{\epsilon}}\left(1+\left|\zeta\right|\right)^{-m}d\zeta\leq C_m\epsilon.
			\end{align*}
			Similarly, since $|\phi(y)|^2=\left|\displaystyle\int_{1}^{y}\phi'(t)dt\right|^2\leq(1-y)\displaystyle\int_{0}^{1}(\phi'(t))^2dt$, we have
			\begin{align*}
				\int_{0}^{1}\left(1+\left|\frac{y-1}{\epsilon}\right|\right)^{-m}\phi^2(y)dy\leq&
\left[\int_{0}^{1}\left(1+\left|\frac{y-1}{\epsilon}\right|\right)^{-m}(1-y)dy\right]\int_{0}^{1}(\phi'(y))^2dy\\
				\leq&C_m\epsilon^2\int_{0}^{1}(\phi'(y))^2dy.
			\end{align*}
		\end{proof}
	\end{Lemma}
	\begin{Proposition}\label{prop:u0 is 0}
		The limiting affine $u_0$ satisfies $u_0(1)=u_0(0)=0$, thus $u_0(y)\equiv0.$
	\end{Proposition}
	\begin{proof}
		By (\ref{ineq:exact test residual estimate}) and Lemma $\ref{lem:weighted Hardy infty}$, we deduce 
		\begin{align}
			\left|\int_{0}^{1}\int_{0}^{2\pi}R_{\epsilon}\phi dxdy\right|\leq&C_m\epsilon\left(\int_{0}^{1}\int_{0}^{2\pi}\rho_m\phi^2dxdy\right)^{\frac{1}{2}}\left(\int_{0}^{1}\int_{0}^{2\pi}\rho_mdxdy\right)^{\frac{1}{2}}\notag\\
			\leq&C_m\epsilon^{\frac{5}{2}}\left(\int_{0}^{1}\int_{0}^{2\pi}\phi_y^2dxdy\right)^{\frac{1}{2}}\leq C_m\epsilon^{\frac{5}{2}},\label{ineq:limiting residual}
		\end{align}
		where we have used $\|(u,v)\|_2=1$. Since $\|\phi\|_{H^1}\leq C\|(u,v)\|_2=C$, there exists a subsequence converges strongly in $L^2(\Omega)$ and weakly in $H^1(\Omega)$ to $\phi_0$. The local affine limit and the inherited selected-wall trace identify
		\begin{align}
			\phi_0(y)=\int_{1}^{y}u_0(s)ds.\label{eq:phi0}
		\end{align}
	Dividing by $\epsilon^2$ and taking $\epsilon\rightarrow 0$ in $\eqref{eq:matched adjoint equality}$,  using $\eqref{ineq:exact test}$, $\eqref{ineq:limiting residual}$  and (\ref{convergence of G}), we deduce that 
		\begin{align*}
			\int_{0}^{1}\int_{0}^{2\pi}\phi_0(y)G_0dxdy=0.
		\end{align*}
		The $H$ contribution vanishes because its horizontal average is zero, so there holds
		\begin{align*}
			\int_{0}^{1}\phi_0(y)\partial_y^4\Big[f_*\Big(\frac{a}{2}(1-y^2)+b(1-y)\Big)\Big]dy=0.
		\end{align*}
		By $f_*(0)=f_*'(0)=0,~~f_*''(0)=u_e(1)^{-1}$ and $\eqref{eq:phi0}$, four integration by parts in $y$ give
		\begin{align*}
			\int_{0}^{1}\phi_0(y)\partial_y^4\Big[f_*\Big(\frac{a}{2}(1-y^2)+b(1-y)\Big)\Big]dy=-u_e(1)u_0(1).
		\end{align*}
		Thus $u_0(1)=0$. Similar argument or the reflection $(x',y')=(-x,1-y),(u',v')=(-u,-v)$ gives $u_0(0)=0$.
	\end{proof}
	\begin{remark}
		For now, we have local convergence $(u,v)\to(0,0)$ in $H^1(K)$ for any compact set $K$ of $\Omega$ and $v\to0$ in global $L^2$. Returning $\eqref{eq:energy1}$ in Lemma $\ref{lem:energy}$, noticing that $\|(F_1,F_2)\|_2=o(\epsilon^2)$ and using Hardy inequalities, we have
		\begin{align*}
			\epsilon^2\int_{0}^{1}\int_{0}^{2\pi}u_y^2dxdy\leq&\int_{0}^{1}\int_{0}^{2\pi}(F_1u+F_2v)dxdy-\int_{0}^{1}\int_{0}^{2\pi}u_x^a(u^2-v^2)dxdy-\int_{0}^{1}\int_{0}^{2\pi}(u_y^a+v_x^a)uvdxdy\notag\\
			=&-\int_{0}^{1}\int_{0}^{2\pi}u_x^au^2dxdy-\int_{0}^{1}\int_{0}^{2\pi}u_y^auvdxdy+o(1)\\
			=&-\int_{0}^{1}\int_{0}^{2\pi}((\chi u_{px}^a)+(1-\chi)\hat{u}_{px}^a)u^2dxdy-\epsilon^{-1}\int_{0}^{1}\int_{0}^{2\pi}(\chi u_{p\zeta}^a+(1-\chi)\hat{u}_{p\eta}^a)uvdxdy+o(1)\\
			\leq&C\epsilon^2\delta\int_{0}^{1}\int_{0}^{2\pi}\left(\frac{u^2}{(y-1)^2}+\frac{u^2}{y^2}\right)dxdy+C\epsilon\delta\int_{0}^{1}\int_{0}^{2\pi}|\sigma|\left(\frac{|uv|}{(y-1)^2}+\frac{|uv|}{y^2}\right)dxdy+o(1)\\
			\leq&C\epsilon^2\delta\int_{0}^{1}\int_{0}^{2\pi}u_y^2dxdy+C\delta\int_{0}^{1}\int_{0}^{2\pi}|\sigma|v_y^2dxdy+o(1).
		\end{align*}
		It follows from Lemma $\ref{lem:energy}$ that
		\begin{align*}
			\int_{0}^{1}\int_{0}^{2\pi}|\sigma|(u_x^2+v_x^2)dxdy\leq&C\epsilon^2\delta\int_{0}^{1}\int_{0}^{2\pi}u_y^2dxdy+o(1),
		\end{align*}
		which implies that
		$$\int_{0}^{1}\int_{0}^{2\pi}|\sigma|(u_x^2+v_x^2)dxdy+\epsilon^2\int_{0}^{1}\int_{0}^{2\pi}u_y^2dxdy\to0.$$
		Combining the first convergence result above with the strong $H^1$ convergence on a fixed overlapping core, we have
		\begin{align}
			\int_{0}^{1}\int_{0}^{2\pi}(u_x^2+v_x^2)dxdy+\epsilon^2\int_{0}^{1}\int_{0}^{2\pi}u_y^2dxdy\to0.\label{ineq:v H1 convergence}
		\end{align}
	\end{remark}

Next, we carry out the step 2. 
	The next argument excludes the possibility that a mean horizontal velocity concentrates in a layer of thickness between $\epsilon$ and a fixed collar width, i.e.  we will prove that $u\to0$ in $L^2(W)$ for the upper collar $W:=[0,2\pi]\times[1-d_0,1]$, where $d_0< \frac12(1-y_\ast)$.
	
	In the upper collar, $u^a\geq c>0$, so we can define $t(x,y)=\displaystyle\int_{y}^{1}u^a(x,z)dz$, which is different from the one used previously. Write $y=Y(x,t)$, then
	\begin{align*}
		Y_x=\frac{v^a}{u^a},~~Y_t=-\frac{1}{u^a}.
	\end{align*}
	The Jacobian and its reciprocal are uniformly bounded. Choose a fixed $T>0$ such that the interface $t=T$ lies in the exact outer region for every sufficiently small $\epsilon$. 

Denote $\tilde{W}=[0,2\pi]\times[0,T]$ and define the circulation $C(t)$ by
	\begin{align*}
		C(t):=\frac{1}{2\pi}\int_{\int_{y}^{1}u^a(x,z)dz=t}udx+vdy=P_0\left(\hat{u}+\frac{v^a}{u^a}\hat{v}\right),
	\end{align*}
	where hat denotes pullback to $(x,t)$ and
	\begin{align*}
		P_0g:=\frac{1}{2\pi}\int_{0}^{2\pi}g(x)dx
	\end{align*}
	for any integrable periodic function $g$. All coefficient quotients in this section are pull back to the stream line coordinates. 
	\begin{Lemma}\label{lem:C equation}
		The circulation $C(t)$ in $[0,T]$ satisfies
		\begin{align}
			C'(t)=-P_0\left(\frac{\hat{\omega}}{u^a}\right)\label{eq:circulation first derivative}
		\end{align}
		and the exact scalar ODE
		\begin{align}
			(\omega_e^ay_0+b)^2C''-\omega_e^aC'=P+S'+\epsilon^{-2}(R_N-R_F),\label{eq:circulation equation}
		\end{align}
		where
		\begin{align*}
			P&=P_0\left(\frac{(\omega^a-\omega_e^a)\hat{\omega}}{u^a}\right),~~
			S=P_0\left(\frac{(u^a)^2+(v^a)^2-(\omega_e^ay_0+b)^2}{u^a}\hat{\omega}\right),\\
			R_N&=P_0\displaystyle\left((\omega^a-\omega_e^a)\left(\hat{v}-\frac{v^a}{u^a}\hat{u}\right)\right),~~
			R_F=P_0\displaystyle\left(\hat{F_1}+\frac{v^a}{u^a}\hat{F_2}\right)
		\end{align*}
		and $y_0=y_0(t)$ satisfies
		\begin{align*}
			\int_{y_0}^{1}(\omega_e^az+b)dz=t.
		\end{align*}
	\end{Lemma}
	\begin{proof}
		Direct differentiation give
		\begin{align*}
			C'(t)=&\frac{1}{2\pi}\int_{0}^{2\pi}Y_t\partial_y\left(u+\frac{v^a}{u^a}v\right)(x,Y(x,t))dx\\
			=&\frac{1}{2\pi}\int_{0}^{2\pi}\frac{-1}{u^a}\left(\frac{v^a}{u^a}v_y+\frac{v_y^au^a-v^au_y^a}{(u^a)^2}v+u_y\right)(x,Y(x,t))dx\\
			=&\frac{1}{2\pi}\int_{0}^{2\pi}\frac{-1}{u^a}\left(\frac{dv}{dx}-v_x+\frac{v_y^au^a-v^au_y^a}{(u^a)^2}v+u_y\right)(x,Y(x,t))dx\\
			=&\frac{1}{2\pi}\int_{0}^{2\pi}\left[v\frac{d}{dx}\left(\frac{1}{u^a}\right)+\frac{-1}{u^a}\left(-v_x+\frac{v_y^au^a-v^au_y^a}{(u^a)^2}v+u_y\right)\right](x,Y(x,t))dx\\
			=&\frac{1}{2\pi}\int_{0}^{2\pi}v\left(\partial_x\left(\frac{1}{u^a}\right)+\frac{v^a}{u^a}\partial_y\left(\frac{1}{u^a}\right)\right)(x,Y(x,t))dx\\
			&+\frac{1}{2\pi}\int_{0}^{2\pi}\frac{-1}{u^a}\left(-v_x+\frac{v_y^au^a-v^au_y^a}{(u^a)^2}v+u_y\right)(x,Y(x,t))dx\\
			=&-\frac{1}{2\pi}\int_{0}^{2\pi}\frac{u_y-v_x}{u^a}(x,Y(x,t))dx
			=-P_0\left(\frac{\hat{\omega}}{u^a}\right),
		\end{align*}
		where $\displaystyle\frac{d}{dx}=\frac{\partial}{\partial x}+\frac{v^a}{u^a}\frac{\partial}{\partial y}$ is the total derivative along the streamline, this completes the proof of (\ref{eq:circulation first derivative}).

Multiplying $\eqref{eq:linearEstimate}_2$ by $\frac{v^a}{u^a}$, adding it to $\eqref{eq:linearEstimate}_1$, we obtain
\begin{align*}
	\epsilon^2\left(\omega_y-\frac{v^a}{u^a}\omega_x\right)+N-\left(p_x+\frac{v^a}{u^a}p_y\right)=-\left(F_1+\frac{v^a}{u^a}F_2\right),
\end{align*}
where
\begin{align*}
	N=-u^a\left(u_x+\frac{v^a}{u^a}v_x\right)-v^a\left(u_y+\frac{v^a}{u^a}v_y\right)-u\left(u_x^a+\frac{v^a}{u^a}v_x^a\right)-v\left(u_y^a+\frac{v^a}{u^a}v_y^a\right).
\end{align*}
Along the streamline, we have $\displaystyle p_x+\frac{v^a}{u^a}p_y=\frac{dp}{dx}$ and
\begin{align*}
	N=&-u^a\left(u_x+\frac{v^a}{u^a}u_y\right)-v^a\left(v_x+\frac{v^a}{u^a}v_y\right)-u\left(u_x^a+\frac{v^a}{u^a}u_y^a\right)\\
	&+u\frac{v^a}{u^a}(u_y^a-v_x^a)-v(u_y^a-v_x^a)-v\left(v_x^a+\frac{v^a}{u^a}v_y^a\right)\\
	=&-u^a\frac{du}{dx}-u\frac{du^a}{dx}-v^a\frac{dv}{dx}-v\frac{dv^a}{dx}+\omega^a\left(\frac{v^a}{u^a}u-v\right)\\
	=&-\frac{d}{dx}(u^au+v^av)+\omega^a\left(\frac{v^a}{u^a}u-v\right),\\
	\frac{v^a}{u^a}u-v=&\frac{d\phi}{dx}.
\end{align*}
Thus, integrating along the streamline and using the tact that $w^a_e$ is a constant, we obtain 
		\begin{align}
			\epsilon^2P_0\left(\omega_y-\frac{v^a}{u^a}\omega_x\right)=R_N-R_F.\label{eq:circulation vorticity equation}
		\end{align}
		Moreover, by integrating by part and $\eqref{eq:circulation first derivative}$, we deduce that 
		\begin{align}
			P_0\left(\omega_y-\frac{v^a}{u^a}\omega_x\right)=&P_0\left(-u^a\hat{\omega}_t-\frac{v^a}{u^a}(\hat{\omega}_x+v^a\hat{\omega}_t)\right)\notag\\
			=&P_0\left(-\frac{(u^a)^2+(v^a)^2}{u^a}\hat{\omega}_t+\frac{d}{dx}\left(\frac{v^a}{u^a}\right)\hat{\omega}\right)\notag\\
			=&P_0\left(-\frac{(u^a)^2+(v^a)^2}{u^a}\hat{\omega}_t-\left[\left(\frac{(u^a)^2+(v^a)^2}{u^a}\right)_t+\frac{\omega^a}{u^a}\right]\hat{\omega}\right)\notag\\
			=&-\frac{d}{dt}P_0\left(\frac{(u^a)^2+(v^a)^2}{u^a}\hat{\omega}\right)-P_0\left(\frac{\omega^a}{u^a}\hat{\omega}\right)\notag\\
			=&-S'+\frac{d}{dt}((\omega_e^ay_0+b)^2C')-P+\omega_e^aC'\notag\\[5pt]
			=&(\omega_e^ay_0+b)^2C''+2\omega_e^ay_0'(\omega_e^ay_0+b)C'+\omega_e^aC'-S'-P,\label{eq:circulation second derivative}
		\end{align}
		where we use the equality
		\begin{align*}
			\frac{d}{dx}\left(\frac{v^a}{u^a}\right)=-\left[\left(\frac{(u^a)^2+(v^a)^2}{u^a}\right)_t+\frac{\omega^a}{u^a}\right]
		\end{align*}
		which can be directly computed as follows,
		\begin{align*}
			\left(\frac{(u^a)^2+(v^a)^2}{u^a}\right)_t+\frac{d}{dx}\left(\frac{v^a}{u^a}\right)=&\frac{-1}{u^a}\left(\frac{(u^a)^2+(v^a)^2}{u^a}\right)_y+\partial_x\left(\frac{v^a}{u^a}\right)+\frac{v^a}{u^a}\partial_y\left(\frac{v^a}{u^a}\right)=-\frac{\omega^a}{u^a}.
		\end{align*}
		 Thus, combining $\eqref{eq:circulation vorticity equation}$, $\eqref{eq:circulation second derivative}$ and $y_0'=-\displaystyle\frac{1}{\omega_e^ay_0+b}$,
 we obtain $\eqref{eq:circulation equation}$.
	\end{proof}
	\begin{Lemma}\label{lem:circulation equation solution}
		The solution of $\eqref{eq:circulation equation}$ can be expressed as
		\begin{align}
			C(t)=&H_{\epsilon}(t)+\int_{t}^{T}G(t,s)P(s)ds-\int_{t}^{T}G_s(t,s)S(s)ds\nonumber\\
			&\qquad \qquad \qquad \qquad +\epsilon^{-2}\int_{t}^{T}G(t,s)(R_N(s)-R_F(s))ds,\label{eq:C expression}
		\end{align}
		where
		\begin{align*}
			H_{\epsilon}(t)=C(T)+\frac{-(\omega_e^ay_0(T)+b)^2C'(T)+S(T)}{\omega_e^ay_0(T)+b}(y_0(t)-y_0(T))
		\end{align*}
		and
		\begin{align*}
			G(t,s)=\frac{y_0(t)-y_0(s)}{\omega_e^ay_0(s)+b},~~0\leq t\leq s\leq T.
		\end{align*}
		 Moreover, there holds 
\begin{align}\label{estimate on Green and homogenes}
|G(t,s)|\leq C(s-t),~~|G_s(t,s)|\leq C,~~ \|H_{\epsilon}\|_{L^\infty([0,T])}\to0(\epsilon\to0).
\end{align}
	\end{Lemma}
	\begin{proof}
		Let $Q$ denote $P+\epsilon^{-2}(R_N-R_F)$ and $J(t)$ denote $\omega_e^ay_0(t)+b$, then $\eqref{eq:circulation equation}$ can be rewritten as 
		\begin{align}
			(J(t)C')'=\frac{S'+Q}{J(t)},~~0\leq t\leq T.\label{eq:circulation equation2}
		\end{align}
		Integrating the equation in $[t,T]$, we obtain that
		\begin{align*}
			C'(t)=\frac{J(T)C'(T)}{J(t)}-\frac{1}{J(t)}\int_{t}^{T}\frac{S'(s)+Q(s)}{J(s)}ds.
		\end{align*}
		Integrating in $[t,T]$ again, we have
		\begin{align*}
			C(t)=C(T)-J(T)C'(T)\int_{t}^{T}\frac{ds}{J(s)}+\int_{t}^{T}\frac{dr}{J(r)}\int_{r}^{T}\frac{S'(s)+Q(s)}{J(s)}ds.
		\end{align*}
		It follows from $y_0'(t)=-J(t)^{-1}$ that
		\begin{align*}
			C(t)=C(T)-J(T)C'(T)(y_0(t)-y_0(T))+\int_{t}^{T}G(t,s)(S'(s)+Q(s))ds,
		\end{align*}
		where we interchange the order of integration for the second integral. After integration by parts for $S$ term, we obtain $\eqref{eq:C expression}$. 

Since
		\begin{align*}
			\int_{y_0}^{1}(\omega_e^az+b)dz=&t=\int_{y}^{1}u^a(x,z)dz\\
			=&\int_{y}^{1}(\omega_e^az+b+u_p^{(0)}+\hat{u}_p^{(0)})dz+O(\epsilon\delta)\\
			=&\int_{y}^{1}(\omega_e^az+b)dz+O(\epsilon\delta),
		\end{align*}
		so 
\begin{align}\label{estimate on y0}
|y_0(t)-y|\leq C\epsilon\delta,~~0<c_1\leq J(t)\leq C, ~~\forall t\in [0,T].
\end{align}
Consequently, $|G(t,s)|\leq C(s-t)$ and $|G_s(t,s)|\leq C$ for $0\leq t\leq s\leq T$. 
		
		To prove $C(T)\to0(\epsilon\to0)$. We choose a number $\tilde{T}>T$ such that the streamline $\displaystyle\int_{y}^{1}u^a(x,z)dz=t$ lies in the outer region for any $t\in[T,\tilde{T}]$. Then by trace theorem, $\eqref{eq:circulation first derivative}$ and interior $H^2$ estimates, we have
		\begin{align*}
			|C(T)|\lesssim&\|C'\|_{L^2(T,\tilde{T})}+\|C\|_{L^2(T,\tilde{T})}\lesssim\|\omega\|_{L^2(K)}+\|u\|_{L^2(K)}+\|v\|_{L^2(K)}\to0(\epsilon\to0),
		\end{align*}
		where we have suppressed the constants in the calculations and $K$ is the compact set of $\Omega$ containing these considered streamlines. Similarly, by trace theorem and interpolation inequality, we have
		\begin{align*}
			|C'(T)|+|S(T)|&\lesssim\|\hat{\omega}(T)\|_{L^2([0,2\pi])}\lesssim\|\omega\|_{H^\frac{3}{4}(K)}\lesssim\|\omega\|_{L^2(K)}^{\frac{1}{4}}\|\omega\|_{H^1(K)}^{\frac{3}{4}}\to0(\epsilon\to0).
		\end{align*}
Thus, we obtain $\|H_{\epsilon}\|_{L^\infty([0,T])}\to0(\epsilon\to0)$, this completes the proof of this Lemma. 
	\end{proof}

	In order to estimate non-homogeneous terms, we need the following Lemma.
	
	\begin{Lemma}[Weighted Volterra bounds]\label{lem:weighted Volterra bounds}
		For $j=0,1$ and a positive weight $\rho$ satisfying $\rho(|\zeta|)\leq\frac{C}{(1+|\zeta|)^m},m>2$, the operators
		\begin{align*}
			(T_{j,\epsilon}f)(t)=\int_{t}^{T}(s-t)^j\rho\left(\frac{s}{\epsilon}\right)f(s)ds
		\end{align*} 
		satisfy
		\begin{align*}
			\|T_{0,\epsilon}f\|_{L^2(0,T)}\leq C\epsilon\|f\|_{L^2(0,T)},\quad\|T_{1,\epsilon}f\|_{L^2(0,T)}\leq C\epsilon^2\|f\|_{L^2(0,T)}.
		\end{align*}
	\end{Lemma}
	\begin{proof}
		By Cauchy inequality, we deduce 
		\begin{align*}
			|T_{0,\epsilon}f(t)|\leq&\left(\int_{t}^{T}\rho\left(\frac{s}{\epsilon}\right)ds\right)^{\frac{1}{2}}\left(\int_{t}^{T}\rho\left(\frac{s}{\epsilon}\right)f^2(s)ds\right)^{\frac{1}{2}},\\
			|T_{1,\epsilon}f(t)|\leq&\left(\int_{t}^{T}(s-t)\rho\left(\frac{s}{\epsilon}\right)ds\right)^{\frac{1}{2}}\left(\int_{t}^{T}(s-t)\rho\left(\frac{s}{\epsilon}\right)f^2(s)ds\right)^{\frac{1}{2}}.
		\end{align*}
		Hence,
		\begin{align*}
			\|T_{0,\epsilon}f\|_{L^2(0,T)}^2\leq&C\epsilon\int_{0}^{T}\int_{t}^{T}\rho\left(\frac{s}{\epsilon}\right)f^2(s)dsdt\\
			=&C\epsilon\int_{0}^{T}s\rho\left(\frac{s}{\epsilon}\right)f^2(s)ds
			\leq C\epsilon^2\int_{0}^{T}f^2(s)ds,\\
			\|T_{1,\epsilon}f\|_{L^2(0,T)}^2\leq&C\epsilon^2\int_{0}^{T}\int_{t}^{T}(s-t)\rho\left(\frac{s}{\epsilon}\right)f^2(s)dsdt\\
			\leq&C\epsilon^2\int_{0}^{T}s^2\rho\left(\frac{s}{\epsilon}\right)f^2(s)ds
			\leq C\epsilon^4\int_{0}^{T}f^2(s)ds.
		\end{align*}
	\end{proof}
	\textbf{Proof of Proposition $\ref{prop:uv estimate}$:} 
\begin{proof}First, by Lemma $\ref{lem:uxvxEstimate}$, we deduce that for $t\in [0,T]$ and $y\in[1-2\cdot10^{-3},1]$, there holds
	\begin{align*}
		|\omega^a-\omega_e^a|\leq C\epsilon^{-1}\delta\rho_m,~~\Big|\frac{v^a}{u^a}\Big|\leq C\epsilon\delta.
	\end{align*}
	By (\ref{estimate on y0}), we obtain 
	\begin{align*}
		\left|\frac{(u^a)^2+(v^a)^2-(\omega_e^ay_0+b)^2}{u^a}\right|\leq C|u^a-(\omega_e^ay_0+b)|+C|v^a|
		\leq C\delta\rho_m+C\epsilon\delta.
	\end{align*}
By Lemma $\ref{lem:weighted Volterra bounds}$, the $P$ term in $\eqref{eq:C expression}$ is bounded by
	\begin{align*}
		C\epsilon^{-1}\delta\left\|T_{1,\epsilon}\left(\int_{0}^{2\pi}|\hat{\omega}(x,s)|dx\right)\right\|_{L^2(0,T)}\leq&C\epsilon\delta\|\hat{\omega}\|_{L^2(\tilde{W})}\leq C\epsilon\delta(\|u_y\|_2+\|v_x\|_2).
	\end{align*}
	Similarly, the $S$ term is bounded by
	\begin{align*}
		C\delta\left\|T_{0,\epsilon}\left(\int_{0}^{2\pi}|\hat{\omega}(x,s)|dx\right)\right\|_{L^2([0,T])}+C\epsilon\delta\|\hat{\omega}\|_{L^2(\tilde{W})}\leq&C\epsilon\delta\|\hat{\omega}\|_{L^2(\tilde{W})}\leq C\epsilon\delta(\|u_y\|_2+\|v_x\|_2).
	\end{align*}
	The forcing term is bounded by $C\epsilon^{-2}\|(F_1,F_2)\|_2$. It remains to estimate $R_N$. 

First, we decompose $R_N$ into two parts:
$$R_N=P_0\displaystyle\left((\omega^a-\omega_e^a)\hat{v}\right)+P_0\displaystyle\left(-(\omega^a-\omega_e^a)\frac{v^a}{u^a}\hat{u}\right):=R_{N1}+R_{N2}.$$
By Lemma $\ref{lem:weighted Hardy infty}$ and Hardy inequality in $t$, we deduce that 
	\begin{align*}
		\|R_{N1}\|_{L^2(0,T)}\leq C\epsilon^{-3}\delta&\left\|\int_{t}^{T}(s-t)\rho_m\left(\frac{s}{\epsilon}\right)\left(\int_{0}^{2\pi}|\hat{v}(x,s)|dx\right)ds\right\|_{L^2(0,T)}\\
		\leq&C\epsilon^{-3}\delta\left\|\left(\int_{t}^{T}(s-t)^2\rho_m\left(\frac{s}{\epsilon}\right)ds\right)^{\frac{1}{2}}\right\|_{L^2(0,T)}\left(\int_{0}^{T}\int_{0}^{2\pi}\rho_m\left(\frac{s}{\epsilon}\right)\hat{v}^2(x,s)dxds\right)^{\frac{1}{2}}\\
		\leq&C\epsilon^{-1}\delta\left(\int_{0}^{T}\int_{0}^{2\pi}\rho_m\left(\frac{s}{\epsilon}\right)\hat{v}^2(x,s)dxds\right)^{\frac{1}{2}}
		\leq C\delta\|u_x\|_2.
	\end{align*}
where we use $\hat{v}(0)=0$ and $\displaystyle\hat{v}_t=\frac{u_x}{u^a}$.

Similarly, 
	\begin{align*}
		\|R_{N2}\|_{L^2(0,T)}\leq C\epsilon^{-2}\delta^2&\left\|\int_{t}^{T}(s-t)\rho_m\left(\frac{s}{\epsilon}\right)\left(\int_{0}^{2\pi}|\hat{u}(x,s)|dx\right)ds\right\|_{L^2(0,T)}\\
		\leq&C\epsilon^{-2}\delta^2\left\|\left(\int_{t}^{T}(s-t)^2\rho_m\left(\frac{s}{\epsilon}\right)ds\right)^{\frac{1}{2}}\right\|_{L^2(0,T)}\left(\int_{0}^{T}\int_{0}^{2\pi}\rho_m\left(\frac{s}{\epsilon}\right)\hat{u}^2(x,s)dxds\right)^{\frac{1}{2}}\\
		\leq&C\delta^2\left(\int_{0}^{T}\int_{0}^{2\pi}\rho_m\left(\frac{s}{\epsilon}\right)\hat{u}^2(x,s)dxds\right)^{\frac{1}{2}}
		\leq C\epsilon\delta^2\|u_y\|_2.
	\end{align*}
	Combining the above estimates, we obtain
	\begin{align*}
		\|C(\cdot)-H_{\epsilon}\|_{L^2(0,T)}\leq&C\delta(\epsilon\|(u_y,v_x)\|_2+\|u_x\|_2)+C\epsilon^{-2}\|(F_1.F_2)\|_2.
	\end{align*}
	Therefore, by Poincar\'e inequality and $\|\hat{v}\|_{L^2(\tilde{W})}\leq C\|u_x\|_2$ by its zero trace at $t=0$, we have
	\begin{align}
		\|\hat{u}-H_{\epsilon}\|_{L^2(\tilde{W})}\leq&\|\hat{u}-P_0(\hat{u})\|_{L^2(\tilde{W})}+\|P_0(\hat{u})-C(\cdot)\|_{L^2(\tilde{W})}
+\|C(\cdot)-H_{\epsilon}\|_{L^2(\tilde{W})}\nonumber\\
		\leq&C\left\|\frac{d}{dx}\hat{u}\right\|_{L^2(\tilde{W})}+C\left\|P_0\left(\frac{v^a}{u^a}\hat{v}\right)\right\|_{L^2(\tilde{W})}
+\|C(\cdot)-H_{\epsilon}\|_{L^2(\tilde{W})}\nonumber\\
		\leq&C(\|u_x\|_2+\epsilon\delta\|(u_y,v_x)\|_2+\epsilon^{-2}\|(F_1,F_2)\|_2).\label{ineq:wall convergence}
	\end{align}
Thus, by (\ref{ineq:v H1 convergence}) and (\ref{estimate on Green and homogenes}), we deduce $\|u\|_{L^2(W)}\rightarrow 0(\epsilon\rightarrow 0).$
Similar argument gives the convergence result in the lower collar.
	
Finally, combining  $\|v\|_{L^2(\Omega)}\to0$ and the interior convergence of $u$, we show that $\|(u,v)\|_2\to0$, contradicting $\eqref{blow up assumption}$. This complete the proof of this Proposition.
\end{proof}
	
	\textbf{Proof of Proposition $\ref{prop:estimate}:$}
	\begin{proof}
	Combining Proposition $\ref{prop:wall estimate}$, Proposition $\ref{prop:uv estimate}$ and Corollary $\ref{cor:v gradient estimate}$, we have
	\begin{align}
		\int_{0}^{1}\int_{0}^{2\pi}(u_x^2+v_x^2)dxdy+\epsilon^2\int_{0}^{1}\int_{0}^{2\pi}u_y^2dxdy\leq&C\epsilon^{-4}\int_{0}^{1}\int_{0}^{2\pi}(F_1^2+F_2^2)dxdy.\label{ineq:linear stability estimate}
	\end{align}
Direct computation gives
	\begin{align*}
		\|S_u\|_2\leq&C\|u_x\|_2+C\epsilon\delta\|u_y\|_2+C\delta\|u\|_2+C\|v(u^a-u_e^a)_y\|_2,
	\end{align*}
	 By $\eqref{ineq:linear stability estimate}$ and
	\begin{align*}
		\|v(u^a-u_e^a)_y\|_2\leq&C\epsilon^{-1}\|\rho_mv\|_2\leq C\|v_y\|_2=C\|u_x\|_2,
	\end{align*}
we deduce 
$$\|S_u\|_2\leq C\epsilon^{-2}\|(F_1,F_2)\|_2.$$
	Similarly, we have
	\begin{align*}
		\|S_v\|_2\leq&C\|(u_x,v_x)\|_2+C\epsilon\delta\|u\|_2\leq C\epsilon^{-2}\|(F_1,F_2)\|_2.
	\end{align*}
	Thus, by the periodic-channel Stokes estimate, we have
	\begin{align*}
		\epsilon^2\|\nabla^2(u,v)\|_2\leq C\epsilon^{-2}\|(F_1,F_2)\|_2.
	\end{align*}
	So $\eqref{ineq:linear stability}$ holds. 

Furthermore, there holds
	\begin{align*}
		\|\nabla S_u\|_2
		\leq&C\|\nabla u^a\|_{\infty}\|u_x\|_2+C\|u^a\|_{\infty}\|\nabla u_x\|_2+C\|\nabla v^a\|_{\infty}\|u_y\|_2+C\|v^a\|_{\infty}\|\nabla u_y\|_2\\[3pt]
		&+C\|\nabla u\|_2\|u_x^a\|_{\infty}+C\|u\|_2\|\nabla u_x^a\|_{\infty}+C\|\nabla v\|_2\|u_y^a\|_{\infty}+C\|v\|_2\|\nabla u_y^a\|_{\infty}\\[3pt]
		\leq&C\epsilon^{-1}\|(u,v)\|_{H^1}+C\|u\|_{H^2}+C\epsilon^{-2}\|v\|_2\leq C\epsilon^{-4}\|(F_1,F_2)\|_2
	\end{align*}
	and
	\begin{align*}
		\|\nabla S_v\|_2
		\leq& C\|\nabla u^a\|_{\infty}\|v_x\|_2+C\|u^a\|_{\infty}\|\nabla v_x\|_2+C\|\nabla v^a\|_{\infty}\|v_y\|_2+C\|v^a\|_{\infty}\|\nabla v_y\|_2\\[3pt]
		&+C\|\nabla u\|_2\|v_x^a\|_{\infty}+C\|u\|_2\|\nabla v_x^a\|_{\infty}+C\|\nabla v\|_2\|v_y^a\|_{\infty}+C\|v\|_2\|\nabla v_y^a\|_{\infty}\\[3pt]
		\leq&C\epsilon^{-1}\|(u,v)\|_{H^1}+C\|v\|_{H^2}\leq C\epsilon^{-4}\|(F_1,F_2)\|_2.
	\end{align*}
	By Stokes estimate again, we deduce 
	\begin{align*}
		\epsilon^2\|\nabla^3(u,v)\|_2\leq C\|(F_1,F_2)\|_{H^1}+C\|S_u\|_{H^1}+C\|S_v\|_{H^1}
		\leq C\|(F_1,F_2)\|_{H^1}+C\epsilon^{-4}\|(F_1,F_2)\|_2.
	\end{align*}
	So $\eqref{ineq:H3}$ holds.
	
	The existence of $\eqref{eq:linearEstimate}$ is equivalent to the existence of the system
	$$L_{\epsilon}\textbf{u}:=\mathbb{P}[-\epsilon^2\Delta\textbf{u}+\textbf{u}^a\cdot\nabla\textbf{u}+\textbf{u}\cdot\nabla\textbf{u}^a]=\mathbb{P}\textbf{F}\text{ in }D,$$
	where $\textbf{u}^a=(u^a,v^a)^T,\textbf{u}=(u,v)^T,\textbf{F}=(F_1,F_2)^T$. For fixed $\epsilon>0$ and $\boldsymbol{\tilde{F}}\in L_{\sigma}^2(\Omega)$, by the classical argument in elliptic equation theory,  $L_{\lambda}\textbf{u}:=L_{\epsilon}\textbf{u}+\lambda_{\epsilon}\textbf{u}=\boldsymbol{\tilde{F}}$ has s solution in $D$ and the energy estimate holds if $\lambda_{\epsilon}>0$ is large enough. Therefore,$L_{\lambda}^{-1}$ is compact in $L_{\sigma}^2(\Omega)$. 
	
	By Fredholm alternative, for any $\textbf{F}\in L^2(\Omega)$, $L_{\epsilon}\textbf{u}=\mathbb{P}\textbf{F}$ has a solution in $D$ iff 
	$$(I-\lambda_{\epsilon}L_{\lambda}^{-1})\textbf{u}=L_{\lambda}^{-1}\mathbb{P}\textbf{F}$$ 
	has a solution in $D$, the compactness of $L_{\lambda}^{-1}$ shows that it is equivalent to the uniqueness of the system $(I-\lambda_{\epsilon}L_{\lambda}^{-1})\textbf{u}=\textbf{0}$ in $L_{\sigma}^2(\Omega)$, which means $L_{\epsilon}\textbf{u}=\textbf{0}$ has a unique solution in $D$. So the existence follows from $\eqref{ineq:linear stability}$.
\end{proof}
	
	\section{Existence for the error equations}
	We now can apply the contraction mapping theorem to prove the existence for the error equations $\eqref{eq:errEqu}.$
	\begin{proposition}\label{prop:errExist}
		There exists $\delta_0>0$ such that for any fixed $\delta\in(0,\delta_0)$, two constants $\epsilon_0>0$ and $C>0$ depending on $\alpha,\beta,f,g,\delta$ exist such that for any $\epsilon\in(0,\epsilon_0)$, the error system $\eqref{eq:errEqu}$ has a solution $(u,v)\in D\cap[H^3(\Omega)]^2$ satisfying
		$$\|(u,v)\|_{\infty}\leq C\epsilon^4.$$
		and
		\begin{align}
			\|\nabla^3(u,v)\|_2\leq C\epsilon^2.\label{ineq:error H3 estimate}
		\end{align}
	\end{proposition}
	\begin{proof} 
For each function $(u,v)\in D$, we set $(\tilde{u},\tilde{v},\tilde{p})$ to be the solution of the following equations
	\begin{align}
		\begin{cases}
			-\epsilon^2\triangle\tilde{u}+\tilde{p}_x+S_{\tilde{u}}=R_u\text{ in }\Omega,\\[3pt]
			-\epsilon^2\triangle\tilde{v}+\tilde{p}_y+S_{\tilde{v}}=R_v\text{ in }\Omega,\\[3pt]
			\tilde{u}_x+\tilde{v}_y=0\text{ in }\Omega,\\[3pt]
			\tilde{u}(x+2\pi,y)=\tilde{u}(x,y),~\tilde{v}(x+2\pi,y)=\tilde{v}(x,y)\text{ in }\Omega,\\[3pt]
			\tilde{u}|_{y=0,1}=\tilde{v}|_{y=0,1}=0,~x\in[0,2\pi],\nonumber
		\end{cases}
	\end{align}
	where
	$$S_{\tilde{u}}:=u^a\tilde{u}_x+v^a\tilde{u}_y+\tilde{u}u_x^a+\tilde{v}u_y^a,$$
	$$S_{\tilde{v}}:=u^a\tilde{v}_x+v^a\tilde{v}_y+\tilde{u}v_x^a+\tilde{v}v_y^a.$$
	By Proposition $\ref{prop:estimate}$ and $\eqref{ineq:linear stability estimate}$, we have $(\tilde{u},\tilde{v})\in D$ and
	\begin{align}
		\|(\tilde{u},\tilde{v})\|_2+\|(\tilde{u}_x,\tilde{v}_x)\|_2+\epsilon\|\tilde{u}_y\|_2+\epsilon^2\|\nabla^2(\tilde{u},\tilde{v})\|_2\leq&C\epsilon^{-2}\|(R_u,R_v)\|_2.\label{est:exist equation}
	\end{align}
	Direct calculation gives
	\begin{align}
		\|(R_u,R_v)\|_2\leq&C\|(R_u^a,R_v^a)\|_2+C\|(u,v)\|_{\infty}\|(u_x,u_y,v_x)\|_2.\label{ineq:residual}
	\end{align}
	We define the norm
	\begin{align*}
		\|(u,v)\|_X:=\|(u,v)\|_2+\|(u_x,v_x)\|_2+\epsilon\|u_y\|_2+\epsilon^2\|\nabla^2(u,v)\|_2.
	\end{align*}
	Then $\eqref{est:exist equation}$ and $\eqref{ineq:residual}$ show that
	\begin{align}
		\|(\tilde{u},\tilde{v})\|_X\leq&C\epsilon^6+C\epsilon^{-4}\|(u,v)\|_X^2,\label{ineq:contraction}
	\end{align}
	where we have used $\|(R_u^a,R_v^a)\|_2=O(\epsilon^8)$ and the fact that
	\begin{align*}
		\|(u,v)\|_{\infty}\leq&C\|(u,v)\|_2^{\frac{1}{2}}\|(u,v)\|_{H^2}^{\frac{1}{2}}\leq C\epsilon^{-1}\|(u,v)\|_X.
	\end{align*}
	Let $X:=\{(u,v)\in D:\|(u,v)\|_X<+\infty\}$,
	then for sufficiently small $\delta_0>0$ and any fixed $\delta\in(0,\delta_0)$, there exists an $\epsilon_0>0$ such that the operator
	$$(u,v)\longmapsto(\tilde{u},\tilde{v})$$
	maps $B:=\{(u,v)\in X:\|(u,v)\|_X\leq\epsilon^{5}\}$ into itself. To prove the operator is a contraction mapping, assume that $(u_1,v_1),(u_2,v_2)$ are any two pairs in $B$ and we set
	\begin{align*}
		\tilde{U}:=\tilde{u}_1-\tilde{u}_2,~\tilde{V}:=\tilde{v}_1-\tilde{v}_2,~\tilde{P}:=\tilde{p}_1-\tilde{p}_2,
	\end{align*}
	then $(\tilde{U},\tilde{V})$ satisfies
	\begin{align}
		\begin{cases}
			-\epsilon^2\triangle\tilde{U}+S_{\tilde{U}}+\tilde{P}_x=R_U,\\[3pt]
			-\epsilon^2\triangle\tilde{V}+S_{\tilde{V}}+\tilde{P}_y=R_V,\\[3pt]
			\partial_x\tilde{U}+\partial_y\tilde{V}=0,\\[3pt]
			\tilde{U}(x+2\pi,y)=\tilde{U}(x,y),\tilde{V}(x+2\pi,y)=\tilde{V}(x,y),\\[3pt]
			\tilde{U}|_{y=0,1}=\tilde{V}|_{y=0,1}=0,\nonumber
		\end{cases}
	\end{align}
	where 
	\begin{align*}
		R_U&\overset{\triangle}{=}R_{u_1}-R_{u_2}=(u_2-u_1)\partial_xu_2+u_1\partial_x(u_2-u_1)+(v_2-v_1)\partial_yu_2+v_1\partial_y(u_2-u_1),\\
		R_V&\overset{\triangle}{=}R_{v_1}-R_{v_2}=(u_2-u_1)\partial_xv_2+u_1\partial_x(v_2-v_1)+(v_2-v_1)\partial_yv_2+v_1\partial_y(v_2-v_1).
	\end{align*}
	Then following the estimate $\eqref{ineq:contraction}$ line by line, we obtain
	$$\|(\tilde{u}_1-\tilde{u}_2,\tilde{v}_1-\tilde{v}_2)\|_X\leq C\epsilon^{-4}(\|(u_1,v_1)\|_X+\|(u_2,v_2)\|_X)\|(u_1-u_2,v_1-v_2)\|_X,$$
	which implies that this operator is a contraction mapping in $B$. 
	
	Finally, we prove $\eqref{ineq:error H3 estimate}$, by $\eqref{ineq:H3}$ in Proposition $\ref{prop:estimate}$ and Sobolev embedding $H^1\subseteq L^4$, we have
	\begin{align*}
		\epsilon^2\|\nabla^3(u,v)\|_2\leq&C\|\nabla(R_u,R_v)\|_2+C\epsilon^{-4}\|(R_u,R_v)\|_2\\[3pt]
		\leq&C\|\nabla(R_u^a,R_v^a)\|_2+C\epsilon^{-4}\|(R_u^a,R_v^a)\|_2+C\epsilon^{-4}\|(u,v)\|_{\infty}\|(u_x,u_y,v_x)\|_2\\[3pt]
		&+\|(u,v)\|_{\infty}\|\nabla(u_x,u_y,v_x)\|_2+C\|\nabla(u,v)\|_4^2\\[3pt]
		\leq&C\epsilon^4+C\epsilon^{-6}\|(u,v)\|_X^2+C\|(u,v)\|_{H^2}^2\\[3pt]
		\leq&C\epsilon^4+C\epsilon^{-6}\|(u,v)\|_X^2.
	\end{align*}
	
	This completes the proof.\\
\end{proof}

Now we give the proof of Theorem $\ref{thm:main}$.

{\bf Proof of Theorem \ref{thm:main}}:

\textit{Proof.} First, the constructed solution takes the following form: 
\[
u^{\epsilon}=u_e^a+\bar{u}_p^a+h(x,y)+u,~~v^{\epsilon}=v_e^a+\bar{v}_p^a+v.
\]
Thus, combining Proposition $\ref{prop:errExist}$ with the Euler estimates and Prandtl estimates in Proposition \ref{lem:approxEstimates}, we obtain  (\ref{estimate on velocity}).

For $\eqref{estimate on vorticity}$, we note that $u_e'(y)=a$, which yields
\begin{align*}
	\omega^{\epsilon}-a=\omega+\partial_y(u^a-u_e(y))-v_x^a.
\end{align*}
It follows from Sobolev embedding $H^2\subseteq L^{\infty}$ and $\eqref{ineq:error H3 estimate}$ in Proposition $ \ref{prop:errExist}$ that
\begin{align*}
	\|\omega\|_{L^{\infty}}\le C\|(u_y,v_x)\|_{L^{\infty}}\leq C\|(u,v)\|_{H^3}\leq C\epsilon^2.
\end{align*}
Moreover,by the properties of approximate solution,we derive that for $0<y_1<y_2<1$,
\begin{align*}
	\|\partial_y(u^a-u_e(y))-v_x^a\|_{L^{\infty}([0,2\pi]\times[y_1,y_2])}\le C\|((u^a-u_e(y))_y,v_x^a)\|_{L^{\infty}([0,2\pi]\times[y_1,y_2])}\leq C_{y_1,y_2}\epsilon.
\end{align*}
So $\eqref{estimate on vorticity}$ holds.$\hfill{\square}$

\section{Appendix}

{ \bf Appendix A:  A formal Prandtl–Batchelor selection argument in a periodic channel.}

In this appendix, we derive a  Prandtl-Batchelor theory for the  steady Navier-Stokes equations on a horizontally periodic strip. For Prandtl-Batchelor theory on a simply connected
domain, we refer to \cite{B, K1998, K2000} or Appendix C in \cite{FGLT1}. Also see generalized Prandtl-Batchelor theory on an annulus in Appendix C in \cite{FGLT2}.

Consider the two-dimensional  Navier-Stokes equations (\ref{eq:NS}).
Let $w^\varepsilon=u^\varepsilon_y-v^\varepsilon_x$ be the vorticity, then a direct computation yields
\begin{align*}
(v^\varepsilon,-u^\varepsilon)^T w^\varepsilon+\nabla\Big(p^\varepsilon+\frac12|\mathbf{u}^\varepsilon|^2\Big)=\varepsilon^2\Delta\mathbf{u}^\varepsilon.
\end{align*}
Let $\Phi^\varepsilon(x,y)$ be the stream function which satisfies $\partial_y\Phi^\varepsilon=u^\varepsilon, -\partial_x\Phi^\varepsilon=v^\varepsilon$. Assuming that
\begin{align*}
\Gamma^\varepsilon_c:=\{(x,y)|\Phi^\varepsilon(x,y)=c\}
\end{align*}
is a closed curve for any $c\in \Phi^\varepsilon(\Omega)$.
By integrating the Navier-Stokes equations over the closed curve $\Gamma^\varepsilon_c$ we deduce that
\begin{align*}
\int_{\Gamma^\varepsilon_c}\Big((v^\varepsilon,-u^\varepsilon)^T  w^\varepsilon+\nabla\Big(p^\varepsilon+\frac12|\mathbf{u}^\varepsilon|^2\Big)\Big)d\vec{\tau}=0,
\end{align*}
thus we obtain that
\begin{align*}
\int_{\Gamma^\varepsilon_c}\Delta\mathbf{u}^\varepsilon d\vec{\tau}=0.
\end{align*}
Assume that $\mathbf{u}^\varepsilon\rightarrow \mathbf{u}_e$ in some strong sense and $\mathbf{u}_e$ is the solution of steady Euler equations, then there holds
\begin{align}
\int_{\Gamma_c^e}\Delta\mathbf{u}_ed\vec{\tau}=0, \nonumber
\end{align}
where $\Gamma_c^e$ is the stream line of the Euler flow.

If the inviscid limit of Navier-Stokes flow is a shear flow $(u_e(y),0)^T$,
then there also holds
$$u''_{e}(y)=0.$$

{ \bf Appendix B: Construction of corrector $h(x,y)$ .}

We give a construction of corrector $h(x,y)$ defined in section $\ref{subsection:Approximate solution}$. First, we give a simple lemma which can be refered to 6.1 in Appendix B in \cite{FGLT1}.

\begin{lemma}\label{correctorH}
	Assume that $K(x,y)$ is a $2\pi-$periodic smooth function which satisfies
	$$\int_{0}^{2\pi}K(x,y)dx=0,\forall y\in[0,1];~K(x,0)=K(x,1)=0,$$
	then there exists a $2\pi-$periodic function $h(x,y)$ such that
	$$\partial_xh(x,y)=K(x,y);~h(x,0)=h(x,1)=0;$$
	$$\int_{0}^{2\pi}h(x,y)dx=0,\|\partial_x^j\partial_y^kh\|_2\leq C\|\partial_x^j\partial_y^kK\|_2.$$
\end{lemma}
Then we construct the corrector $h(x,y)$ by the above lemma. Direct computation gives
\begin{align*}
	\partial_x(u_e^a+\bar{u}_{p}^a)+\partial_y(v_e^a+\bar{v}_p^a)&=\chi'(y)(v_p^a-\hat{v}_p^a):=-K(x,y).
\end{align*}
Since $\chi'(y)=0$ near the wall, we have
\begin{align*}
	|K(x,y)|&=|\epsilon^m\chi'(y)(y-1)^{-m}\zeta^{-m}\zeta^{2m}v_p^a+\epsilon^m\chi'(y)y^{-m}\eta^{-m}\eta^{2m}\hat{v}_p^a|\\
	&\leq C_m\epsilon^m\delta(|\zeta|^{-m}+\eta^{-m}).
\end{align*}
Notice that $\chi'(y)=0$ for $y\in[0,1/4]\cup[3/4,1]$ and the properties of $v_p^a,\hat{v}_p^a$, we deduce that $K(x,y)$ satisfies the assumption in Lemma $\ref{correctorH}$,  there exists $h(x,y)$ such that
$$\partial_x(u_e^a+\bar{u}_p^a)+\partial_y(v_e^a+\bar{v}_p^a)=-\partial_xh(x,y),$$
$$h(x,0)=h(x,1)=0,\|\partial_x^j\partial_y^kh\|_2\leq C_{j,k,m}\epsilon^m\delta\rho_m,$$
which indicates that
$$\partial_xu^a+\partial_yv^a=0.$$

\section*{Acknowledgments}

 T.Tao is partially supported by the NSF of China under Grant 12371234 and 12671282. \\

\textbf{Data availability statement}

We do not analyse or generate any datasets, because our work proceeds within a theoretical and mathematical approach.\\

\textbf{Conflict of interest}

On behalf of all authors, the corresponding author states that there is no conflict of interest.\\

\textbf{Disclosure of AI assistance} During the preparation of this manuscript, the authors
used ChatGPT to assist in exploring proof strategies,
checking technical estimates, particularly the $L^2$
estimate for the linearized remainder equations, and
improving the exposition. The authors are responsible
for independently verifying the mathematical arguments
and take full responsibility for the content and
conclusions of the paper.

\renewcommand{\refname}{References} 


\begin{thebibliography}{9}
	
	
	\bibitem {B}G.K.Batchelor, \textit{On steady laminar flow with closed
		streamlines at large Reynolds number}, J.Fluid Mech., 7 (1956), no.1, 177-190.
	
	
	
	
	\bibitem {CWZ}Q.Chen, D.Wu and Z.Zhang, \textit{On the stability of shear
		flows of Prandtl type for the steady Navier-Stokes equations}, Sci. China Math., 66 (2023), no.4, 679-722.
	
	\bibitem{CFLZ}Z.Chen, M.Fei, Z.Lin and J.Zhao, \textit{Prandtl-Batchelor flows with a point vortex on a disk}, arXiv:2410.07759.
	
	\bibitem {DM}A.Dalibard and N.Masmoudi, \textit{Separation for the stationary
		Prandtl equation}, Publ.Math.Inst.Hautes \'{E}tudes Sci., 130 (2019), 187-297.
	
	
	\bibitem{DIT}T.Drivas, S.Iyer and  T.Nguyen, \textit{The Feynman-Lagerstrom criterion for boundary layer}, Arch. Ration. Mech. Anal., 248 (2024), no. 3, Paper No. 55, 41 pp.
	
	
	
	\bibitem {edwards}D.A.Edwards, \textit{Viscous Boundary-Layer Effects in Nearly
		Inviscid Cylindrical Flows}, Nonlinearity, 10 (1997), 277--290.
	
	
	\bibitem{FGLT1}M.Fei, C.Gao, Z.Lin and T.Tao, \textit{Prandtl-Batchelor flows on a disk}, Comm.Math.Phys., 397 (2023), no.3, 1103-1161.
	
	\bibitem{FGLT2}M.Fei, C.Gao, Z.Lin and T.Tao, \textit{Prandtl-Batchelor flows on an annulus}, Adv. Math, 458 (2024), Paper No. 109994, 86 pp.
	
	\bibitem{FPZ} M. Fei, X. Pan and J. Zhao, \textit{Vanishing viscosity limits of the Navier-Stokes equations in a horizontally periodic strip}, Science China Mathematics 69(2026), 2337-2372
	
	\bibitem{Pan2026}X. Pan, C. Zhao and J. Zhao, \textit{Vanishing viscosity limit of the stationary Navier-Stokes equations in an expanding annulus}, arXiv:2610.17754 (2026).
	
	\bibitem {feymann-lagerstrom}R.P.Feynman and P.A.Lagerstrom, \textit{Remarks on
		high Reynolds number flows in finite domains},  Proc. IX International Congress
	on Applied Mechanics, Brussels, 3 (1956), 342--343.
	
	\bibitem{G}
	G.P.Galdi,  \textit{An Introduction to the Mathematical Theory of the Navier-Stokes Equations, Steady-State Problems}, 2nd Edition, Springer Monographs in Mathematics, Springer (2011).
	
	
	\bibitem {GZ}C.Gao and L.Zhang, \textit{On the steady Prandtl boundary layer
		expansions}, Sci. China Math., 66 (2023), no. 9, 1993-2020. 
	
	\bibitem {GZ2}C.Gao and L.Zhang, \textit{Remarks on the steady Prandtl
		boundary layer expansions}, J. Math. Phys., 65 (2024), no. 7, Paper No. 071504, 27 pp.
	
	\bibitem {GZO}C.Gao and C.Zhao, \textit{Sharp asymptotic behavior of the steady pressure-free Prandtl system}, arXiv:2504.11870.
	
	\bibitem {GMM}D.Gerard-Varet and Y.Maekawa, \textit{Sobolev stability of
		Prandtl expansions for the steady Navier-Stokes equations}, Arch. Ration.
	Mech. Anal., 233 (2019), no. 3, 1319-1382.
	
	
	\bibitem {GN}Y.Guo and T.Nguyen, \textit{Prandtl boundary layer expansions of
		steady Navier-Stokes flows over a moving plate}, Ann.PDE, 2017, DOI 10.1007/s40818-016-0020-6.
	
	\bibitem {GI1}Y.Guo and S.Iyer, \textit{Validity of steady Prandtl layer
		expansions}, Comm. Pure Appl. Math., 76 (2023), no. 11, 3150-3232. 
	
	\bibitem {GI2}Y.Guo and S.Iyer, \textit{Steady Prandtl layer expansions with
		external forcing}, Quart. Appl. Math., 81 (2023), no. 2, 375-411. 
	
	\bibitem {GI3}Y.Guo and S.Iyer, \textit{Regularity and expansion for steady
		Prandtl equations}, Comm.Math.Phys., 382 (2021), no.3, 1403-1447.
	
	\bibitem{GuoYang2025}
	Y.Guo and Z.Yang, \textit{Small scales in inviscid limits of steady fluids}, SIAM J. Math. Anal. (2025), arXiv:2409.09604.
	
	\bibitem{Guo2026}Y.Guo, H.Li, X.Ma and M.Zhang, \textit{Validity of three-dimensional axisymmetric steady Prandtl expansion}, arXiv:2611.02492 (2026).
	
	
	
	
	\bibitem {Iyer}S.Iyer, \textit{On global-in-x  stability of Blasius profiles},
	Arch. Ration. Mech. Anal. 237 (2020), no. 2, 951-998.
	
	\bibitem {Iyer1}S.Iyer, \textit{Steady Prandtl boundary layer expansions over
		a rotating disk}, Arch.Ration.Mech.Anal., 224 (2017), no.2, 421-469.
	
	\bibitem {Iyer2}S.Iyer, \textit{Global steady Prandtl boundary layer over a
		moving boundary}, Peking Math. J., I: 2 (2019), no. 2, 155-238; II: 2 (2019),
	no. 3-4, 353-437; III: 3 (2020), no. 1, 47-102.
	
	\bibitem {Iyer3}S.Iyer, \textit{Steady Prandtl boundary layer over a moving
		boundary:nonshear Euler flow}, SIAM, J.Math.Anal., 51 (2019), no.3, 1657-1685.
	
	\bibitem {IM1}S.Iyer and N.Masmoudi, \textit{Boundary layer expansions of
		steady Navier-Stokes equation}, Ars Inven. Anal.(2021), Paper No.6, 47 pp.
	
	\bibitem {IM2}S.Iyer and N.Masmoudi, \textit{Global-in-x stability of steady
		Prandtl expansions for 2D Navier-Stokes flows}, arXiv:2008.12347.
	
	\bibitem {Iyer4}S.Iyer and C.Zhou, \textit{Stationary inviscid limit to shear
		flow}, J. Differential Equations 267 (2019), no.12, 7135-7153.
	
	\bibitem{Iyer2026}
	S.Iyer, \textit{Stability of the favorable Falkner--Skan profiles for the stationary Prandtl equations}, Adv. Math. \textbf{487} (2026), 110749.
	
	\bibitem{JLY} H. Jia, Z. Lei and C. Yuan, \textit{Sharp asymptotic stability of Blasius profile in the steady Prandtl equation}, Adv. Math. 480 (2025), 110533.
	
	\bibitem {K1998}S.C.Kim, \textit{On Prandtl-Batchelor theory of a
		cylindrical eddy: asymptotic study}, SIAM J.Appl.Math., 58 (1998), 1394-1413.
	
	\bibitem {K2000}S.C.Kim, \textit{On Prandtl-Batchelor theory of a
		cylindrical eddy: existence and uniqueness}, Z.Angew.Math.Phys.,
	51 (2000),674-686.
	
	\bibitem {K2003}S.C.Kim, \textit{Asymptotic study of Navier-Stokes flows},
	Trends in Mathematics, information center for Mathematical sciences, 6 (2003), 29-33.
	%
	\bibitem {kim-childress}S.C.Kim and S.Childress, \textit{Vorticity selection with
		multiple eddies in two-dimensional steady flow at high Reynolds number}, SIAM
	J.Appl.Math., 61 (2001), no.5, 1605--1617.
	
	\bibitem {kim-free boundary}S.C.Kim, \textit{A free-boundary problem for Euler
		flows with constant vorticity}, Appl.Math.Lett.,12 (1999),101--104.
	
	\bibitem {kim-thesis}S.C.Kim, \textit{On Prandtl-Batchelor theory of steady flow at
		large Reynolds number}, Ph.D thesis, New York University, 1996.
	
	\bibitem {kim-formula}S.C.Kim, \textit{Batchelor-Wood formula for negative wall
		velocity}, Phys. Fluids, 11 (1999), 1685-1687.
	
	\bibitem{kim-annular} S.C.Kim and J.Y.Lee, \textit{Prandtl-Batchelor theory for an annular domain}, SIAM J. Appl. Math. 86 (2026), 133–159.
	
	
	\bibitem{KKR}
	M.V.Korobkov, K.Pileckas and R.Russo, \textit{Solution of Leray's problem for stationary Navier-Stokes equations in plane and axially symmetric spatial domains}, Ann. of Math. (2)181 (2015), no.2, 769-807.
	
	\bibitem{Leray}
	J.Leray, \textit{Etude de diverses equations integrales non lineaire et de quelques
		problemes que pose lhydrodynamique}, J.Math.Pures Appl., 12 (1933), 1-82.
	
	
	
	
	
	
	
	
	
	
	\bibitem {OK} H.Okamoto, \textit{A variational problem arising in the
		two-dimensional Navier-Stokes equation with vanishing viscosity}, Appl.Math.Lett., 7 (1994), no.1, 29-33.
	
	\bibitem {OK1} H.Okamoto, \textit{Nearly singular two-dimensional Kolmogorov flows for large Reynolds numbers}, J.Dynam.Differential Equations, 8 (1996), no.2, 203-220.
	
	\bibitem {OS} O.A.Oleinik and V.N.Samokhin, \textit{Mathematical models in
		boundary layer theory}, Applied Mathematics and Mathematical Computation, 15.
	Champan and Hall/CRC, Boca Raton, FL, 1999.
	
	
	\bibitem {prandtl} L.Prandtl, 1905. \"{U}ber Fl\"{u}ssigkeitsbewegung bei sehr
	kleiner Reibung. Verhandlungen des III. Internationalen Mathematiker
	Kongresses, Heidelberg, 1904, pp. 484--491, Teubner, Leizig. See Gesammelte
	Abhandlungen II, pp. 575--584.
	
	
	
	\bibitem {renardy}M.Renardy, \textit{On non-existence of steady periodic solutions of
		the Prandtl equations}, J. Fluid Mech., 717 (2013), R7, 1-5.
	
	
	
	
	\bibitem {riley81}N.Riley, \textit{High Reynolds number flows with closed
		streamlines}, J. Eng. Math., 15 (1981), 15-27.
	
	
	\bibitem {van wijngarden}L.Van Wijngaarden, \textit{Prandtl-Batchelor flows
		revisited}, Fluid Dyn. Res., 39 (2007), 267--278.
	
	\bibitem {SWZ}W.Shen, Y.Wang and Z.Zhang, \textit{Boundary layer separation
		and local behavior for the steady Prandtl equation}, Adv. Math., 389 (2021),107896, 25pp.
	
	\bibitem {WZ}Y.Wang and Z.Zhang, \textit{Global $C^{\infty}$ regularity of the
		steady Prandtl equation with favorable pressure gradient},  Ann. Inst. H. Poincar\'{e} Anal. Non Lin\'{e}aire, 38 (2021), no. 6, 1989-2004.
	
	\bibitem {WZ2}Y.Wang and Z.Zhang, \textit{Asymptotic behavior of the steady Prandtl equation}, Math. Ann. 387 (2023), no. 3-4, 1289-1331.
	
	
	
	
	
	\bibitem {W}W.W.Wood, \textit{Boundary layers whose streamlines are closed},
	J.Fluid Mech, 01 (1957), no.2, 77-87.
	
\end{thebibliography}
\end{document}